\documentclass[10pt]{amsart}

\usepackage{amsmath}
\usepackage{amsthm}
\usepackage{stmaryrd}
\numberwithin{equation}{subsection}
\renewcommand{\theequation}{\arabic{equation}}
\usepackage{url}
\usepackage{amssymb}
\usepackage{amsxtra}
\usepackage[all]{xy}
\makeindex

\DeclareMathOperator{\depth}{depth}
\swapnumbers
\newtheorem{prop}[subsection]{Proposition}
\newtheorem{thm}[subsection]{Theorem}
\newtheorem{lem}[subsection]{Lemma}
\newtheorem{cor}[subsection]{Corollary}
\newtheorem{gathlem}[subsection]{Gathering Lemma}
\newtheorem{LNT}[subsection]{Local Norm Theorem}

\theoremstyle{definition}
\newtheorem{defn}[subsection]{Definition}
\newtheorem{eg}[subsection]{Example}

\theoremstyle{remark}
\newtheorem{rmk}[subsection]{Remark}

\DeclareMathOperator{\Nrd}{Nrd}

\newcommand{\Implies}{\Rightarrow}

\renewcommand{\theenumi}{\roman{enumi}}

\newcommand{\F}{\mathbb{F}}

\newcommand{\iso}{\xrightarrow{\sim}}
\newcommand{\ot}{\otimes}
\newcommand{\val}{\lambda}
\newcommand{\s}{\sigma}

\newcommand{\Fb}{\overline{F}}

\newcommand{\Db}{\overline{D}}
\newcommand{\Fx}{F^\times}
\newcommand{\Lx}{L^\times}

\newcommand{\Ox}{\mfo^\times}

\newcommand{\floor}[1]{\left\lfloor{#1}\right\rfloor}

\newcommand{\texp}{\tex}
\newcommand{\qform}[1]{{\left\langle{#1}\right\rangle}}                   % a quadratic form
\newcommand{\pform}[1]{{\dla{#1}\rangle\!\rangle}} % a Pfister form
\newcommand{\qpform}[1]{{\dla{#1}\rrbracket}}   % a characteristic 2 Pfister form

\newcommand{\balpha}{\bar{\alpha}}
\newcommand{\bbeta}{\bar{\beta}}
\newcommand{\bmu}{\bar{\mu}}
\newcommand{\bu}{\bar{u}}
\newcommand{\ba}{\bar{a}}
\newcommand{\bb}{\bar{b}}
\newcommand{\bc}{\bar{c}}
\newcommand{\be}{\bar{\varepsilon}}
\newcommand{\bx}{\bar{x}}
\newcommand{\bpsi}{\bar{\psi}}
\newcommand{\by}{\bar{y}}

\newcommand{\bv}{\bar{v}}
\newcommand{\bchi}{\bar{\chi}}

\newcommand{\dlog}[1]{\frac{\der{#1}}{#1}}
\newcommand{\der}{\mathrm{d}}

\newcommand{\eand}{\quad\text{and}\quad}
\newcommand{\Z}{\mathbb{Z}}
\newcommand{\Q}{\mathbb{Q}}
\newcommand{\Zm}[1]{\Z/{#1}\Z}

\DeclareMathOperator{\gr}{gr}
\DeclareMathOperator{\chr}{char}
\DeclareMathOperator{\chpoly}{chpoly}

\newcommand{\ord}{O}

\DeclareMathOperator{\tex}{texp}  %\newcommand{\tex}{\mathrm{texp}}
\DeclareMathOperator{\nexp}{nexp}
\newcommand{\nex}{\nexp}

\def\lz{\par \vspace{8pt}}
\def\IR{\mathbb{R}}   % Das sind normalerweise Befehle aus rslybo.sty
\def\IF{\mathbb{F}}   % zum Ausdruck der Zahlbereiche. Aus Gewohnheits-
\def\IC{\mathbb{C}}   % gr\"unden werden sie hier synonym f\"ur die ams-Befehle
\def\IO{\mathbb{O}}   % genutzt.

\def\IZ{\mathbb{Z}}
\def\IQ{\mathbb{Q}}

\newcounter{flabel}%
{%
\begin{list}{\mbox{}\hfill(\roman{flabel})\hfill\mbox{}}%
{\itemsep0pt\parsep0pt\topsep0ex%
\settowidth{\labelwidth}{(xiii)}%
\setlength{\leftmargin}{\labelwidth}%
\addtolength{\leftmargin}{.8em}%
\usecounter{flabel}}}%
{\end{list}}
\begin{document}
\setcounter{section}{0} \setcounter{subsection}{-1}
\setcounter{tocdepth}{1}
\sloppy
\newcommand{\Eins}{\mathbf{1}}
\newcommand{\Sp}{\mathcal{S}}
\newcommand{\Alb}{\mathcal{A}}
\def\klz{\par \vspace{2pt}}
\newcommand{\ssp}{\vspace{5pt}}
\newcommand{\bsp}{\vspace{16pt}}
\newcommand{\zp}{\hspace{-4,75pt}.}
\newcommand{\kalg}{k\mathchar45\mathbf{alg}}
\newcommand{\Ralg}{R\mathchar45\mathbf{alg}}
\newcommand{\Salg}{S\mathchar45\mathbf{alg}}
\newcommand{\bfa}{\mathbf{a}}
\newcommand{\bfb}{\mathbf{b}}
\newcommand{\bfh}{\mathbf{h}}
\newcommand{\bfs}{\mathbf{s}}
\newcommand{\bft}{\mathbf{t}}
\newcommand{\bfAut}{\mathbf{Aut}}
\newcommand{\bfT}{\mathbf{T}}
\newcommand{\bfM}{\mathbf{M}}
\newcommand{\bfN}{\mathbf{N}}
\newcommand{\bfD}{\mathbf{D}}
\newcommand{\bfG}{\mathbf{G}}
\newcommand{\bfX}{\mathbf{X}}
\newcommand{\bfY}{\mathbf{Y}}
\newcommand{\bfGL}{\mathbf{GL}}
\newcommand{\bfSL}{\mathbf{SL}}
\newcommand{\bfkap}{\boldsymbol{\kappa}}
\newcommand{\bfSb}{\mathbf{Splbas}}
\newcommand{\clh}{\mathcal{H}}
\newcommand{\clm}{\mathcal{M}}
\newcommand{\clp}{\mathcal{P}}
\newcommand{\mfg}{\mathfrak{g}}
\newcommand{\mfp}{\mathfrak{p}}
\newcommand{\p}{\mathfrak{p}}
\newcommand{\mfq}{\mathfrak{q}}
\newcommand{\q}{\mathfrak{q}}
\newcommand{\mfm}{\mathfrak{m}}
\newcommand{\mfo}{\mathfrak{o}}
\newcommand{\mfgl}{\mathfrak{gl}}
\newcommand{\mfD}{\mathfrak{D}}
\newcommand{\mfL}{\mathfrak{L}}
\newcommand{\mfO}{\mathfrak{O}}
\newcommand{\mfsl}{\mathfrak{sl}}
\newcommand{\hA}{\hat{A}}
\newcommand{\hD}{\hat{D}}
\newcommand{\hL}{\hat{L}}
\newcommand{\hR}{\hat{R}}
\newcommand{\ha}{\hat{a}}
\newcommand{\hb}{\hat{b}}
\newcommand{\oB}{\overline{B}}
\newcommand{\oC}{\overline{C}}
\newcommand{\oa}{\overline{a}}
\newcommand{\ob}{\overline{b}}
\newcommand{\oab}{\overline{ab}}
\newcommand{\odel}{\overline{\delta}}
\newcommand{\on}{\overline{n}}
\newcommand{\ou}{\overline{u}}
\newcommand{\ov}{\overline{v}}
\newcommand{\ow}{\overline{w}}
\newcommand{\ox}{\overline{x}}
\newcommand{\oy}{\overline{y}}
\newcommand{\oxy}{\overline{xy}}
\newcommand{\bfGm}{\mathbf{G}_\mathrm{m}}
\newcommand{\bfGmA}{\mathbf{G}_{\mathrm{m}A}}
\newcommand{\bfGmB}{\mathbf{G}_{\mathrm{m}B}}
\newcommand{\bfGmk}{\mathbf{G}_{\mathrm{m}k}}
\newcommand{\sesanT}{\langle T\rangle_{\mathrm{ses}}}
\newcommand{\sesanIntbT}{\langle b^{-1}Tb \rangle_{\mathrm{ses}}}
\newcommand{\sesancotrT}{\langle\overline{T}^t\rangle_{\mathrm{ses}}}
\newcommand{\sesanmu}{\langle \mu\rangle_{\mathrm{ses}}}
\newcommand{\sesanmumfp}{\langle \mu_{\mfp}\rangle_{\mathrm{ses}}}
\newcommand{\sesanDH}{\langle \det_\Delta\,h \rangle_{\mathrm{ses}}}
\newcommand{\mbalg}{\mathrm{alg}}
\newcommand{\mbCa}{C_\mathrm{a}}
\newcommand{\mbCOa}{C_{0\mathrm{a}}}
\newcommand{\mbCROa}{C_{0R\mathrm{a}}}
\newcommand{\mbla}{\mathrm{lev}_{\mathrm{ass}}}
\newcommand{\mblc}{\mathrm{lev}_{\mathrm{com}}}
\newcommand{\mbMa}{M_\mathrm{a}}
\newcommand{\mbc}{\mathrm{char}}
\newcommand{\mbd}{\mathrm{dim}}
\newcommand{\mbi}{\mathrm{inf}}
\newcommand{\mbm}{\mathrm{min}}
\newcommand{\mbal}{\mathrm{alt}}
\newcommand{\mba}{\mathrm{ass}}
\newcommand{\hass}{h_{\mathrm{ass}}}
\newcommand{\hcom}{h_{\mathrm{com}}}
\newcommand{\hgcom}{\mathrm{hgt}_{\mathrm{com}}}
\newcommand{\hgass}{\mathrm{hgt}_{\mathrm{ass}}}
\newcommand{\mbj}{\mathrm{Jord}}
\newcommand{\mbC}{\mathrm{Cay}}
\newcommand{\mbE}{\mathrm{End}}
\newcommand{\mbH}{\mathrm{Hom}}
\newcommand{\mbCe}{\mathrm{Cent}}
\newcommand{\mbN}{\mathrm{Nuc}}
\newcommand{\Nex}{\mathrm{Nex}}
\newcommand{\mbZ}{\mathrm{Zor}}
\newcommand{\mbS}{\mathrm{Spec}}
\newcommand{\mbP}{\mathrm{Pic}}
\newcommand{\mbPol}{\mathrm{Pol}}
\newcommand{\mbD}{\mathrm{Der}}
\newcommand{\mbID}{\mathrm{InDer}}
\newcommand{\mbMD}{\mathrm{MulDer}}
\newcommand{\mbL}{\mathrm{Lie}}
\newcommand{\mbA}{\mathrm{AssDer}}
\newcommand{\mbSt}{\mathrm{StanDer}}
\newcommand{\mbCo}{\mathrm{ComDer}}
\newcommand{\mbPrr}{\mathrm{Prind}_\mathrm{reg}}
\newcommand{\mbPrs}{\mathrm{Prind}_\mathrm{sing}}
\newcommand{\mbPr}{\mathrm{Prind}}
\newcommand{\mbszo}{^\mathrm{o}}
%%%%%%%%%%%%%%%%%%%%%%%%%%%%%%%%%%%%%%%%%%%%%%%%%%%
\newcommand{\tr}{\mathop \mathrm{tr}\nolimits}
\newcommand{\End}{\mathop \mathrm{End}\nolimits}
\newcommand{\Ker}{\mathop\mathrm{Ker}\nolimits}
\renewcommand{\Im}{\mathop\mathrm{Im}\nolimits}
\newcommand{\Nuc}{\mathop\mathrm{Nuc}\nolimits}
\newcommand{\wid}{\mathop\mathrm{\omega}\nolimits}
\newcommand{\Proof}{\ifvmode\else\par\fi\medskip\noindent \emph{Proof.}\ \ }
\newcommand{\Remark}{\ifvmode\else\par\fi\medskip\noindent \emph{Remark.}\ \ }
\newcommand{\Remarks}{\ifvmode\else\par\fi\medskip\noindent \emph{Remarks.}\ \ }
\renewcommand{\:}{\colon\,}
\newcommand{\Mat}{\mathop\mathrm{Mat}\nolimits}
\newcommand{\sgn}{\mathop \mathrm{sgn}\nolimits}
\newcommand{\rk}{\mathop \mathrm{rk}\nolimits}
\newcommand{\sd}{\mathbf{d}}
\newcommand{\Spec}{\mathrm{Spec}}
\newcommand{\AUT}{\mathbf{Aut}}
\newcommand{\Aut}{\mathrm{Aut}}
\newcommand{\Hom}{\mathop \mathrm{Hom}\nolimits}
\newcommand{\GL}{\mathop \mathrm{GL}\nolimits}
\newcommand{\SL}{\mathop \mathrm{SL}\nolimits}
\newcommand{\Ad}{\mathop \mathrm{Ad}\nolimits}
\newcommand{\Quot}{\mathop \mathrm{Quot}\nolimits}
\newcommand{\Isom}{\mathrm{Isom}}
\newcommand{\ISOM}{\mathbf{Isom}}
\newcommand\set{\mathbf{set}}
\newcommand\grp{\mathbf{grp}}
\newcommand{\la}{\mathbf{\langle}}
\newcommand{\ra}{\mathbf{\rangle}}
\newcommand{\dla}{\langle\!\langle} %\mathbf{\langle\langle}}
\newcommand{\dra}{\rangle\!\rangle} %{\mathbf{\rangle\rangle}}

\long\def\ignore#1\endignore{\relax}
%
%\thispagestyle{fancy}
%\pagestyle{fancy}

%\lhead{\footnotesize \tt \today}

%\chead{}

%\rhead{\footnotesize \tt \jobname.tex}

%\rhead{}

%\lfoot{}

%\cfoot{\thepage}

%\rfoot{}

%\footrulewidth 0.3pt
%\renewcommand{\headrulewidth}{0pt}
% \cfoot{}
%\vspace*{5mm}

%%%%
%
% use title and maketitle so that subjclass is printed out
%
%%%%

\title[Wild Pfister forms]{Wild Pfister forms over Henselian fields, $K$-theory, and conic division algebras}

%Literatur: Osborn [1962], Bruhat-Tits [I, IHES, 1971], Bruhat-Tits
%[Bull.Soc.Math.France 1984]

\author{Skip Garibaldi}
\address{(Garibaldi) Department of Mathematics and Computer Science, Emory University, Atlanta, GA 30322, USA}
\email{skip@mathcs.emory.edu}

\author{Holger P. Petersson}
\address{(Petersson) Fakult\"at f\"ur Mathematik und Informatik, FernUniversit\"at in Hagen, D-58084 Hagen, Germany}
\email{holger.petersson@fernuni-hagen.de}

%\date{\tt Version of \today}

%\dedicatory{\tt Preliminary version}

\begin{abstract}
The epicenter of this paper concerns Pfister quadratic forms over a
field $F$ with a Henselian discrete valuation. All
characteristics are considered but we focus on the most complicated
case where the residue field has characteristic 2 but $F$ does not.
We also prove results about round quadratic forms, composition
algebras, generalizations of composition algebras we call conic
algebras, and central simple associative symbol algebras. Finally we
give relationships between these objects and Kato's filtration on
the Milnor $K$-groups of $F$.
\end{abstract}
\maketitle

\section*{Introduction} \label{s.INTRO}
The theory of quadratic forms over a field $F$ with a
Henselian discrete valuation is well understood in case the
residue field $\Fb$ has characteristic different from 2 thanks to
Springer \cite{MR0070664}.  But when $\Fb$ is imperfect of
characteristic 2, the theorems are much more complicated---see,
e.g., \cite{Tietze} and \cite{MR868606}---reflecting perhaps the
well-known fact that quadratic forms are not determined by
valuation-theoretic data, as illustrated below in Example
\ref{e.PAIRS}.  However, more can be said when one focuses on
Pfister quadratic forms over $F$ and more generally round forms, see
Part \ref{hensel.part} below.

In Part \ref{p.CONKAT} we change our focus to Kato's filtration on
the mod-$p$ Milnor $K$-theory of  a Henselian discretely valued
field of characteristic zero where the residue field has
characteristic $p$.  (Again, if $\Fb$ has characteristic
different from $p$, the mod-$p$ Galois cohomology and Milnor
$K$-theory of $F$ are easily described in terms of $\Fb$, see
\cite[pp.~17--19]{GMS} and \cite[7.1.10]{GilleSz}.)  We give
translations between valuation-theoretic properties of Pfister
forms, octonion algebras, central simple associative algebras of
prime degree (really, symbol algebras), and cyclic field extensions
of prime degree over $F$ on the one hand and properties of the
corresponding symbols in Milnor $K$-theory on the other.

Along the way, we prove some results that are of independent
interest, which we now highlight.  Part \ref{comp.part} treats
quadratic forms and composition algebras over arbitrary fields.  It
includes a Skolem-Noether Theorem for purely inseparable subfields
of composition algebras (Th.~\ref{t.SKONOE}) and a result on
factoring quadratic forms (Prop.~\ref{p.PFISUB}).  We also give a
new family of examples of what we call conic division algebras,
which are roughly speaking division algebras where every element
satisfies a polynomial of degree 2, see Example \ref{e.CODIAL}.
More precisely, we show that---contrary to what is known, e.g., over
the reals---a Pfister quadratic form of characteristic $2$ is
anisotropic if and only if it is the norm of such a division algebra
(Cor.~\ref{c.PFICODI}).

Part \ref{hensel.part} focuses on round quadratic forms and
composition algebras over a field $F$ with a 2-Henselian discrete
valuation with residue field of characteristic 2.  From this part,
our Local Norm Theorem \ref{t.HENOTH} has already been applied in \cite{We09}.  We also
relate Tignol's height $\omega$ from \cite{Tignol:wild} with
Saltman's level $\hgcom$ from \cite{MR589083}, with a nonassociative
version of Saltman's level that we denote by $\hgass$, and with
valuation-theoretic properties of composition algebras, see
Th.~\ref{t.HEWI} and Cor.~\ref{texp.cor}. We show that composition
division algebras over $F$ having pre-assigned valuation data,
subject to a few obvious constraints, always exist
(Cor.~\ref{c.PREASINV}), though they are far from unique up to
isomorphism (Example~\ref{e.PAIRS}).

Part \ref{p.CONKAT} gives a $K$-theoretic proof of the Local Norm
Theorem (Th.~\ref{LNTK}); it has an easy proof but stronger
hypotheses than the version in Part \ref{hensel.part}. Finally, our
Gathering Lemma \ref{gathlem} is independent of the rest of the
paper and says that one may rewrite symbols in a convenient form.

\smallskip

Let us now discuss what happens ``under the hood".  The basic
technical result in Part \ref{comp.part} is a non-orthogonal
analogue of the classical Cayley-Dickson construction for algebras
of degree 2.  It is used to prove the Skolem-Noether Theorem
mentioned above as well as to construct the examples of conic
division algebras.

In Part \ref{hensel.part}, we proceed to a more arithmetic
set-up by considering a base field $F$, discretely valued by a
normalized discrete valuation $\lambda\:F \to \IZ \cup \{\infty\}$,
which is $2$-Henselian in the sense that it satisfies Hensel's lemma
for quadratic polynomials. Our main goal is to understand
composition algebras and Pfister quadratic forms over $F$. For this purpose,
the results of \cite{MR51:635}---where the base field
was assumed to be \emph{complete} rather than Henselian---carry over
to this more general (and also more natural) setting virtually
unchanged; we use them here here without further ado. Moreover, we
will be mostly concerned with ``wild'' composition algebras over $F$
(see \ref{ss.TAWI} below for the precise definition of this term in
a more general context) since a complete description of the ``tame''
ones in terms of data living over the residue field of $F$ has been
given in \cite{MR51:635}. (Alternatively---and from a different perspective---one has a good description of the tame part of the Witt group of $F$ from \cite{ET:springer}.)
The approach adopted here owes much to the work of Kato
\cite[\S~1]{Kato:gen1}, Saltman \cite{MR589083} and particularly
Tignol \cite{Tignol:wild} on wild associative division algebras of
degree the residual characteristic $p > 0$ of their (possibly
non-discrete) Henselian base field. Moreover, our approach is not
confined to composition algebras but, at least to a certain extent,
works more generally for (non-singular) pointed quadratic spaces
that are round (e.g., Pfister) and anisotropic. We attach valuation
data to these spaces, among which not so much the usual ones
(ramification index (\ref{ss.VADAL}~(b)) and pointed quadratic
residue space (\ref{ss.VADAL}~(c))), but wildness-detecting
invariants like the trace exponent (\ref{ss.TITE}) play a
significant role. After imitating the quadratic defect
\cite[63A]{MR0152507} for round and anisotropic pointed quadratic
spaces (\ref{ss.NE}) and extending the local square theorem
\cite[63:1]{MR0152507} to this more general setting
(Thm.~\ref{t.HENOTH}), we proceed to investigate the behavior of our
valuation data when passing from a wild, round and anisotropic
pointed quadratic space $P$ having ramification index $1$ as input
to the output $Q := \dla\mu\dra \otimes P$, for any non-zero scalar
$\mu \in F$ (Section~\ref{s.VADAEN}). In all cases except one, the
output, assuming it is anisotropic, will again be a wild pointed
quadratic space. Remarkably, the description of the exceptional case
(Thm.~\ref{t.CHAGORE}), where the input is assumed to be Pfister and
the output turns out to be tame, when specialized to composition
algebras, relies critically on the non-orthogonal Cayley-Dickson
construction encountered in the first part of the paper
(Cor.~\ref{c.CONTAM}). This connection is due to the fact that our
approach also lends itself to the study of what we call
$\lambda$-normed and $\lambda$-valued conic algebras
(Section~\ref{s.LANOVA}), the latter forming a class of conic
division algebras over $F$ that generalize ordinary composition
algebras and turn out to exist in all dimensions $2^n$, $n =
0,1,2,\dots$, once $F$ has been chosen appropriately
(Examples~\ref{e.LAPR},\ref{e.LAU}). There is yet another unusual
feature of the exceptional case: though exclusively belonging to the
theory of quadratic forms (albeit in an arithmetic setting), it can
be resolved here only by appealing to elementary properties of
flexible conic algebras (Thm.~\ref{t.CONCHA}).
The second part of the paper concludes with extending
Tignol's notion of height \cite{Tignol:wild}, which agrees with
Saltman's notion of level \cite{MR589083}, to composition division
algebras over $F$ and relating them to the valuation data introduced
before (Thm.~\ref{t.HEWI}).

In the third part of the paper, we consider the case where $F$ has
characteristic zero and a primitive $p$-th root of unity and has a
Henselian discrete valuation with residue field of characteristic
$p$.  In that setting, Kato, Bloch, and Gabber gave a description of
the mod-$p$ Milnor $K$-groups $k_q(F)$ in \cite{Kato:gen1,
Kato:gen2, Kato:cdv, BlochKato}; we use \cite{CT:kato} as a
convenient reference.  We relate properties of a symbol in $k_q(F)$
with valuation-theoretic properties of the corresponding algebra,
see Prop.~\ref{dict} and Th.~\ref{texp.depth}.

%\pagebreak

\tableofcontents

%%%%%%%%%%%%%%%%%%%%%%%%%%%%%%%%%%%%%%%
\part{Base fields of characteristic $2$} \label{comp.part}

\section{Standard properties of conic algebras}
\label{s.COALARFI} Although composition algebras are our main
concern in this paper, quite a few of our results remain valid
under far less restrictive conditions. The appropriate framework
for some of these conditions is provided by the category of conic
algebras. They are the subject of the present section.

We begin by fixing some terminological and notational conventions
about non-associative algebras in general and about quadratic
forms. For the time being, we let $k$ be a field of arbitrary
characteristic. Only later on (Sections
\ref{s.PFIBILQUAD}$\--$\ref{s.CODI}) will we confine ourselves to
base fields of characteristic $2$.

\subsection{Algebras.} \label{ss.ALTA} Non-associative (= not
necessarily associative) algebras\index{algebra} play a dominant
role in the present investigation. For brevity, they will often be
referred to simply as algebras (over $k$) or as $k$-algebras. A
good reference for the standard vocabulary is \cite{MR668355}.

Left and right multiplication of a $k$-algebra $A$ will be denoted
by $x \mapsto L_x$ and $x \mapsto R_x$,
respectively\index{algebra!left, right multiplication}. $A$ is
called unital if it has an identity (or unit) element, denoted by
$1_A$. A subalgebra of $A$ is called unital if it contains the
identity element of $A$. Algebra homomorphisms are called unital
if they preserve identity elements.
Commutator\index{algebra!commutator} and
associator\index{algebra!associator} of $A$ will be denoted by
$[x,y] = xy - yx$ and $[x,y,z] = (xy)z - x(yz)$, respectively. If
$A$ is unital, then
\[
\mbN(A) := \big\{x \in A\mid [A,A,x] = [A,x,A] = [x,A,A] =
\{0\}\big\}
\]
is a unital associative subalgebra of $A$, called its
nucleus\index{algebra!nucleus}, and
\[
\mbCe(A) := \big\{x \in \mbN(A)\mid [A,x] = \{0\}\big\}
\]
is a unital commutative associative subalgebra of $A$, called its
centre\index{algebra!centre}. We say $A$ is central (resp.~has
trivial nucleus) if $\mbCe(A) = k1_A$ (resp.~$\mbN(A) = k1_A$).

A $k$-algebra $A$ is called flexible\index{algebra!flexible} if it
satisfies the flexible law
\begin{align}
\label{FLEX}  xyx := (xy)x = x(yx).
\end{align}
$A$ is said to be alternative\index{algebra!alternative} if the
associator is an alternating (trilinear) function of its
arguments. This means that $A$ is flexible and satisfies the left
and right alternative laws
\begin{align}
\label{ALT} x(xy) = x^2y, \quad (yx)x = yx^2.
\end{align}
Furthermore, the left, middle and right Moufang identities
\begin{align}
\label{MOU} x\big(y(xz)\big) = (xyx)z, \quad x(yz)x = (xy)(zx),
\quad \big((zx)y\big)x = z(xyx)
\end{align}
hold.

\subsection{Quadratic forms.} \label{ss.QUAFO}Our main reference
in this paper for the algebraic theory of quadratic forms is
\cite{MR2427530}, although our notation will occasionally be
different and we sometimes work in infinite dimensions. Let $V$ be
vector space over $k$, possibly infinite-dimensional. Deviating
from the notation used in \cite{MR2427530}, we write the polar
form of a quadratic form $q\:V \to k$, also called the bilinear
form associated with $q$ or its bilinearization, as $\partial q$,
so $\partial q\:V \times V \to k$\index{$\partial q$} is the
symmetric bilinear form given by
\begin{align*}
\partial q(x,y) := q(x + y) - q(x) - q(y) &&(x,y \in V).
\end{align*}
Most of the time we simplify notation and write $q(x,y) :=
\partial q(x,y)$ if there is no danger of confusion. The quadratic
form $q$ is said to be \emph{non-singular}\index{quadratic form,
non-singular} if it has finite dimension and its polar form is
non-degenerate in the usual sense, i.e., for any $x \in V$, the
relations $\partial q(x,y) = 0$ for all $y \in V$ imply $x = 0$.
Recall that non-singular quadratic forms have even dimension if
the characteristic is $2$ \cite[Chap.~2, Remark~7.22]{MR2427530}.

\subsection{Conic algebras.} \label{ss.QUALG} We consider a class
of non-associative algebras that most authors refer to as
quadratic \cite{MR668355} or algebras of degree $2$
\cite{MR769825}. In order to avoid confusion with Bourbaki's
notion of a quadratic algebra \cite{MR0354207}, we adopt a
different terminology. A $k$-algebra $C$ is said to be
\emph{conic}\index{conic algebra} if it has an identity element
$1_C \ne 0$ and there exists a quadratic form $n\:C \to k$ with
$x^2 - t(x)x + n(x)1_C = 0$ for all $x \in C$, where $t$ is
defined by $t := \partial n(1_C,-)\:C \to k$  and hence is a
linear form. The quadratic form $n$ is uniquely determined by
these conditions and is called the \emph{norm} of $C$\index{conic
algebra!norm}, written as $n_C$. We call $t_C := t =
\partial n_C(1_C,-)$ the \emph{trace}\index{conic algebra!trace}
of $C$ and have
\begin{align}
\label{QUE} x^2 - t_C(x)x + n_C(x)1_C = 0, \quad n_C(1_C) = 1,
\quad t_C(1_C) = 2 &&(x \in C).
\end{align}
Finally, the linear map\index{$\iota_C$}
\begin{align}
\label{CONJ} \iota_C\:C \longrightarrow C, \quad x \longmapsto
\iota_C(x) := x^\ast := t_C(x)1_C - x,
\end{align}
called the \emph{conjugation} of $C$\index{conic
algebra!conjugation}, has period $2$ and is characterized by the
condition
\begin{align}
\label{CENCON} 1_C^\ast = 1_C, \quad xx^\ast = n_C(x)1_C &&(x \in
C).
\end{align}
The property of an algebra to be conic is inherited by unital
subalgebras. Injective unital homomorphisms of conic algebras are
automatically norm preserving. A conic algebra $C$ over $k$ is said
to be \emph{non-degenerate}\index{conic algebra!non-degenerate} if
the polar form $\partial n_C$ has this property. Thus
finite-dimensional conic algebras are non-degenerate iff their norms
are non-singular as quadratic forms. Orthogonal complementation in
$C$ always refers to $\partial n_C$. We say $C$ \emph{is simple as
an algebra with conjugation} if only the trivial (two-sided) ideals
$I \subseteq C$ satisfy $I^\ast = I$.

\subsection{Invertibility in conic algebras.} \label{ss.INCO} Let
$C$ be a conic algebra over $k$. By (\ref{ss.QUALG}.\ref{QUE}),
the unital subalgebra of $C$ generated by an element $a \in C$,
written as $k[a]$, is commutative associative and spanned by
$1_C,a$ as a vector space over $k$; in particular it has dimension
at most $2$. We say $a$ is \emph{invertible} in $C$\index{conic
algebra!invertible elements} if this is so in $k[a]$, i.e., if
there exists an element $a^{-1} \in k[a]$ (necessarily unique and
called the \emph{inverse} of $a$ in $C$) such that $aa^{-1} =
1_C$. For $a$ to be invertible in $C$ it is necessary and
sufficient that $n_C(a) \neq 0$, in which case $a^{-1} =
n_C(a)^{-1}a^\ast$. The set of invertible elements in $C$ will
always be denoted by $C^\times$.

As usual, a non-associative $k$-algebra $A$ is called a
\emph{division algebra}\index{algebra!division} if for all $a,b
\in A$, $a \neq 0$, the equations $ax = b$, $ya = b$ can be solved
uniquely in $A$. The quest for conic division algebras is an
important topic in the present investigation. The following
necessary criterion, though trivial, turns out to be useful.

\subsection{Proposition.} \label{p.NODIV} \emph{The norm of a
conic division algebra is anisotropic.} \hfill $\square$ \lz

\noindent The converse of this proposition does not hold (cf.
\ref{ss.CODI}). For $\mbc(k) \neq 2$, conditions that are
necessary and sufficient for a conic algebra to be division have
been given by Osborn \cite[Thm.~3]{MR0140550}.

\subsection{Inseparable field extensions.} \label{ss.INSFI} Exotic
examples of conic algebras arise in connection with
inseparability. Suppose $k$ has characteristic $2$ and let $K/k$
be a purely inseparable field extension of exponent at most $1$,
so $K^2 \subseteq k$. Then $K$ is a conic $k$-algebra with $n_K(u)
= u^2$ for all $u \in K$, $\partial n_K = 0$, $t_K = 0$ and
$\iota_K = \Eins_K$. In particular, inseparable field extensions
of exponent at most $1$ over $k$ are degenerate, hence singular,
conic division algebras.

\subsection{Composition algebras.} \label{ss.CONCOM} Composition
algebras form the most important class of conic algebras.
Convenient references, including base fields of characteristic
$2$, are \cite{MR2000a:16031,MR1763974}, although
\cite{MR2000a:16031} introduces a slightly more general notion. An
algebra $C$ over $k$ is said to be a \emph{composition
algebra}\index{composition algebra} if it is non-zero, contains a
unit element and carries a non-singular quadratic form $n\:C \to
k$ that \emph{permits composition}: $n(xy) = n(x)n(y)$ for all
$x,y \in C$. Composition algebras are automatically conic. In
fact, the only quadratic form on $C$ permitting composition is the
norm of $C$ in its capacity as a conic algebra.

\subsection{Basic properties of composition algebras.} \label{ss.PROCOM}
Composition algebras exist only in dimensions $1,2,4,8$ and are
alternative. They are associative iff their dimension is at most
$4$, and commutative iff their dimension is at most $2$. The base
field is a composition algebra if and only if it has
characteristic different from $2$. Composition algebras of
dimension $2$ are the same as quadratic \'etale algebras.
Composition algebras of dimension $4$ (resp.~$8$) are called
quaternion (resp.~octonion or Cayley) algebras. Two composition
algebras are isomorphic if and only if their norms are isometric
(as quadratic forms). The conjugation of a composition algebra $C$
over $k$ is an algebra involution.

What is denied to arbitrary conic algebras holds true for
composition algebras:

\subsection{Norm criterion for division algebras.}
\label{ss.NOCRI} \emph{A composition algebra $C$ over $k$ is a
division algebra if and only if its norm is anisotropic.}
Otherwise its norm is hyperbolic, in which case we say $C$ is
split. Up to isomorphism, split composition algebras are uniquely
determined by their dimension, and their structure is explicitly
known.

\subsection{The Cayley-Dickson construction.} \label{ss.CDCO}\index{Cayley-Dickson construction}
The main tool for dealing with conic algebras in general and
composition algebras in particular is the Cayley-Dickson
construction. Its inputs are a conic algebra $B$ and a non-zero
scalar $\mu \in k$. Its output is a conic algebra $C :=
\mbC(B,\mu)$ that is given on the vector space direct sum $C = B
\oplus Bj$ of two copies of $B$ by the multiplication
\begin{align}
\label{CDMU} (u_1 + v_1j)(u_2 + v_2j) := (u_1u_2 + \mu v_2^\ast
v_1) + (v_2u_1 + v_1u_2^\ast)j &&(u_i,v_i \in B, i = 1,2).
\end{align}
Norm, polarized norm, trace and conjugation of $C$ are related to
the corresponding data of $B$ by the formulas
\begin{align}
\label{CDNO} n_C(u + vj) =\,\,&n_B(u) - \mu n_B(v), \\
\label{CDPN} n_C(u_1+ v_1j,u_2 + v_2j) =\,\,&n_B(u_1,u_2) - \mu
n_B(v_1,v_2), \\
\label{CDTR} t_C(u + vj) =\,\,&t_B(u), \\
\label{CDCO} (u + vj)^\ast =\,\,&u^\ast - vj
\end{align}
for all $u,v,u_i,v_i \in B$, $i = 1,2$. Note that $B$ embeds into
$C$ as a unital conic subalgebra through the first summand; we
always identify $B \subseteq C$ accordingly. The Cayley-Dickson
construction $\mbC(B,\mu)$ is clearly functorial in $B$, under
injective unital homomorphisms.

It is a basic fact that $C$ is a composition algebra iff $B$ is an
associative composition algebra. Conversely, we have the following
embedding property, which fails for arbitrary conic algebras, cf.
Example~\ref{e.LAPRESP} below.

\subsection{Embedding property.} \label{ss.EMPRO} Any proper
composition subalgebra $B$ of a composition algebra $C$ over $k$
is associative and \emph{admits a scalar $\mu \in k^\times$ such
that the inclusion $B \hookrightarrow C$ extends to an embedding
$\mbC(B,\mu) \to C$ of conic algebras.} More precisely, $\mu \in
k^\times$ satisfies this condition iff $\mu = -n_C(y)$ for some $y
\in B^\perp \cap C^\times$. \lz

\subsection{The Cayley-Dickson process.} \label{ss.CDP}\index{Cayley-Dickson process}
Let $B$ be a conic $k$-algebra. Using non-zero scalars
$\mu_1,\dots,\mu_n \in k^\times$ ($n \geq 1$), we write
inductively
\begin{align*}
C := \mbC(B;\mu_1,\dots,\mu_n) :=
\mbC\big(\mbC(B;\mu_1,\dots,\mu_{n-1}),\mu_n\big)
\end{align*}
for the corresponding iterated Cayley-Dickson construction
starting from $B$. It is a conic $k$-algebra of dimension
$2^n\mbd_k(B)$. We say $C$ arises from $B$ and the
$\mu_1,\dots,\mu_n$ by means of the \emph{Cayley-Dickson process}.
The norm of $C$ is given by
\begin{align}
\label{ENBC} n_C = \dla\mu_1,\dots,\mu_n\dra \otimes n_B.
\end{align}
Here are the most important special cases of the Cayley-Dickson
process. $\ssp$ \\
\emph{Case~$1$.} $B = k$, $\mbc(k) \neq 2$. \\
Then $n_C = \dla\mu_1,\dots,\mu_n\dra$ is an $n$-Pfister quadratic
form. $C$ is a composition algebra iff $n \leq 3$. $\ssp$ \\
\emph{Case~$2$.} $B = k$, $\mbc(k) = 2$. \\
Then $n_C = \dla\mu_1,\dots,\mu_n\dra_q$ is a quasi-Pfister
(quadratic) form \cite[\S~10, p.~56]{MR2427530}. Moreover, $n_C$ is
anisotropic iff $C = k(\sqrt{\mu_1},\dots,\sqrt{\mu_n})$ is an
extension field of $k$, necessarily purely inseparable of exponent
$1$, hence never a composition algebra. $\ssp$ \\
\emph{Case~$3$.} $B$ is a quadratic \'etale $k$-algebra.  \\
Then $n_B = \dla\mu\rrbracket$ for some $\mu \in k$
\cite[Example~9.4]{MR2427530} and \eqref{ENBC} shows that
\[
n_C = \dla\mu_1,\dots,\mu_n,\mu\rrbracket
\]
is an $(n+1)$-Pfister quadratic form over $k$. Moreover, $C$ is a
composition algebra iff $n \leq 2$. $\ssp$ \\
Composition algebras other than the base field itself always
contain quadratic \'etale subalgebras. Hence, by Cases~$1,3$ above
and by the embedding property \ref{ss.EMPRO}, they may all be
obtained from each one of these, even from the base field itself
if the characteristic is not $2$, by the Cayley-Dickson process.
The preceding discussion also shows that all Pfister and all
quasi-Pfister quadratic forms are the norms of appropriate conic
algebras.

\subsection{Inseparable subfields.} \label{ss.INSU} Let $C$ be a
composition \emph{division} algebra over $k$. A unital subalgebra
of $C$ is either a composition (division) algebra itself or an
inseparable extension field of $k$; in the latter case, $k$ has
characteristic $2$ and the extension is purely inseparable of
exponent at most $1$ \cite{MR0152555}. The extent to which this
case actually occurs may be described somewhat more generally as
follows.

Suppose $k$ has characteristic $2$, $B$ is a conic $k$-algebra and
$\mu_1,\dots,\mu_n \in k^\times$ are such that the norm of
\begin{align}
C := \mbC(B;\mu_1,\dots,\mu_n)
\end{align}
is anisotropic, so all non-zero elements of $C$ are invertible
(\ref{ss.INCO}). Then, by Case~$2$ of \ref{ss.CDP},
\[
K := \mbC(k;\mu_1,\dots,\mu_n) \subseteq C
\]
is a purely inseparable subfield of degree $2^n$ and exponent $1$.

Specializing this observation to $n = 2$ and $B$ quadratic \'etale
over $k$, we conclude \emph{that every octonion division algebra
over a field of characteristic $2$ contains an inseparable
subfield of degree $4$}.

\section{Flexible and alternative conic algebras.}
\label{s.FLEXALT} This section is devoted to some elementary
properties of flexible and alternative conic algebras. In
particular, we derive expansion formulas for the norm of
commutators and associators that turn out to be especially useful
in subsequent applications.

Phrased with appropriate care, most of the results obtained here
remain valid over any commutative associative ring of scalars. For
simplicity, however, we continue to work over a field $k$ of
arbitrary characteristic. We fix a conic algebra $C$ over $k$ and
occasionally adopt the abbreviations $1 = 1_C$, $n = n_C$, $t =
t_C$.

\subsection{Identities in arbitrary conic algebras.}
\label{ss.IDCON} The following identities, some of which have been
recorded before, are assembled here  for the convenience of the
reader and either hold by definition or are straightforward to
check.
\begin{align}
\label{NUNT} n_C(1_C) =\,\,&1_C, \\
\label{TUNT} t_C(1_C) =\,\,&2, \\
\label{NOLI} t_C(x) =\,\,&n_C(1_C,x), \\
\label{QUAD} x^2 =\,\,&t_C(x)x - n_C(x)1_C, \\
\label{QUALI} x \circ y := xy + yx =\,\,&t_C(x)y + t_C(y)x -
n_C(x,y)1_C, \\
\label{CJ} x^\ast =\,\,&t_C(x)1_C - x, \\
\label{SCIN} xx^\ast =\,\,&n_C(x)1_C, \\
\label{NJ} n_C(x^\ast) =\,\,&n_C(x).
\end{align}

\subsection{Identities in flexible conic algebras.}
\label{ss.FLECON} We now assume that $C$ is flexible. By McCrimmon
\cite[3.4,Thm.~3.5]{MR769825}, this implies the following
relations:
\begin{align}
\label{ADM} n_C(xy,x) =  n_C(x)t_C(y) =\,\,&n_C(yx,x), \\
\label{NAS} n_C(x,zy^\ast) = n_C(xy,z) = \,\,&n_C(y,x^\ast z), \\
\label{TEEN} n_C(x,y) = t_C(xy^\ast) = \,\,&t_C(x)t_C(y) - t_C(xy), \\
\label{TAS} t_C(xy) = t_C(yx), \quad t_C(xyz)
:=\,\,&t_C\big((xy)z\big) = t_C\big(x(yz)\big).
\end{align}
Moreover, the conjugation is an algebra involution of $C$, so we
have $(xy)^\ast = y^\ast x^\ast$ for all $x,y \in C$. Dealing with
flexible conic algebras has the additional advantage that this
property is preserved under the Cayley-Dickson construction
\cite[Thm.~6.8]{MR769825}.

\subsection{Remark.} \label{r.EFLENAS} By \cite[3.4]{MR769825},
each one of the four(!) identities in
(\ref{ss.FLECON}.\ref{ADM}),(\ref{ss.FLECON}.\ref{NAS}) is
actually \emph{equivalent} to $C$ being flexible. \lz

\noindent The norm of a flexible conic algebra will in general not
permit composition. But we have at least the following result.

\subsection{Proposition.} \label{p.NOCOMT} \emph{Let $C$ be a
flexible conic algebra over $k$. Then}
\begin{align}
\label{NOCO} n_C(xy) =\,\,&n_C(yx), \\
\label{NOCOMT} n_C([x,y]) =\,\,&4n_C(xy) - t_C(x)^2n_C(y) -
t_C(y)^2n_C(x) + t_C(xy)t_C(xy^\ast)
\end{align}
\emph{for all $x,y \in C$.}

\Proof Expanding the expression $n(x \circ y)$ by means of
(\ref{ss.IDCON}.\ref{QUALI}),(\ref{ss.IDCON}.\ref{NOLI}) yields
\[
n(x \circ y) = t(x)^2n(y) + t(y)^2n(x) + n(x,y)^2 -
t(x)t(y)n(x,y),
\]
where flexibility allows us to invoke
(\ref{ss.FLECON}.\ref{TEEN}); we obtain
\begin{align}
\label{ENCI} n(x \circ y) = t(x)^2n(y) + t(y)^2n(x) -
t(xy)t(xy^\ast).
\end{align}
Now let $\varepsilon = \pm 1$. Then
\[
n(xy + \varepsilon yx) = n(xy) + \varepsilon n(xy,yx) + n(yx) =
n(xy) + (1 - 2\varepsilon)n(yx) + \varepsilon n(x \circ y,yx),
\]
and combining (\ref{ss.IDCON}.\ref{QUALI}) with
(\ref{ss.FLECON}.\ref{ADM})(\ref{ss.FLECON}.\ref{TEEN}),(\ref{ss.FLECON}.\ref{TAS}),
we conclude
\begin{align}
\label{ENEP} n(xy + \varepsilon yx) = n(xy) + (1 -
2\varepsilon)n(yx) + \,\,&\varepsilon\big(t(x)^2n(y) + t(y)^2n(x)
- t(xy)t(xy^\ast)\big).
\end{align}
Comparing \eqref{ENCI} and \eqref{ENEP} for $\varepsilon = 1$
yields \eqref{NOCO}, while \eqref{NOCO} and \eqref{ENEP} for
$\varepsilon = -1$ yield \eqref{NOCOMT}. \hfill $\square$

\subsection{Proposition.} \label{p.ADSIM}\emph{Let $C$ be a non-degenerate
conic algebra over $k$.} $\ssp$ \\
(a) \emph{$C$ is simple as an algebra with conjugation} (cf.~\ref{ss.QUALG}). $\ssp$ \\
(b) \emph{If $\iota_C$ is an algebra involution of $C$, in
particular, if $C$ is flexible, then $C$ is either simple or split
quadratic \'etale.}

\Proof (a) Let $I \subseteq C$ be an ideal with $I^\ast = I$. For
$x \in I$, $y \in C$ we linearize (\ref{ss.IDCON}.\ref{SCIN}) and
obtain $n_C(x,y)1_C = xy^\ast + yx^\ast \in I$. Then either $I =
C$ or $n_C(x,y) = 0$ for all $x \in I$, $y \in C$, forcing $I =
\{0\}$ by non-degeneracy.

(b) Assuming $C$ is not simple, we must show it has dimension $2$.
Since $(C,\iota_C)$ is simple as an algebra with involution by
(a), there exists a $k$-algebra $A$ such that $(C,\iota_C) \cong
(A^{\mathrm{op}} \oplus A,\varepsilon)$ as algebras with
involution, $\varepsilon$ being the exchange involution of
$A^{\mathrm{op}} \oplus A$. Then $A$, embedded diagonally into
$A^{\mathrm{op}} \oplus A$, and
\begin{align*}
H := H(C,\iota_C) := \{x \in C\mid x = x^\ast\} \subseteq C
\end{align*}
identify canonically as vector spaces over $k$. In particular, not
only the dimension of $H$ but also its codimension in $C$ agree
with the dimension of $A$. Applying (\ref{ss.IDCON}.\ref{CJ}), we
conclude that $H$ has dimension $1$ for $\mbc(k) \neq 2$ and
codimension $1$ for $\mbc(k) = 2$. In both cases, $C$ must be
$2$-dimensional. \hfill $\square$

\subsection{Proposition.} \label{p.DICAN} \emph{Let $C$ be a
flexible conic algebra over $k$ whose norm is anisotropic. Then
either $C$ is non-degenerate or $\partial n_C = 0$.}

\Proof Since, by \ref{ss.INCO}, all non-zero elements of $C$ are
invertible, $C$ is a simple algebra. On the other hand,
(\ref{ss.FLECON}.\ref{NAS}) shows that $I := C^\perp \subseteq C$
is an ideal. If $I = \{0\}$, then $C$ is non-degenerate. If $I =
C$, then $\partial n_C = 0$. \hfill $\square$

\subsection{Identities in conic alternative algebras.}
\label{ss.CONALT} If $C$ is a conic alternative algebra, then by
\cite[p.~97]{MR769825} its norm permits composition:
\begin{align}
\label{CP} n_C(xy) = n_C(x)_Cn(y).
\end{align}
Linearizing \eqref{CP}, we obtain
\begin{align}
\label{CPLL} n_C(x_1y,x_2y) =\,\,&n_C(x_1,x_2)n_C(y), \\
\label{CPRL} n_C(xy_1,xy_2) =\,\,&n_C(x)n_C(y_1,y_2), \\
\label{CPTL} n_C(x_1y_1,x_2y_2) + n_C(x_1y_2,x_2y_1)
=\,\,&n_C(x_1,x_2)n_C(y_1,y_2).
\end{align}
Moreover, by \cite[Prop.~3.9]{MR769825} and
(\ref{ss.FLECON}.\ref{TEEN}),
\begin{align}
\label{UOP} xyx = n_C(x,y^\ast)x - n_C(x)y^\ast = t_C(xy)x -
n_C(x)y^\ast.
\end{align}

\subsection{Theorem.} \label{t.NOAS} \emph{Let $C$ be a
conic alternative algebra over $k$. Then}
\begin{align}
\label{NOAS} n_C([x_1,x_2,x_3]) =\,\,&4n_C(x_1)n_C(x_2)n_C(x_3) -
\sum\,t_C(x_i)^2n_C(x_j)n_C(x_l) + \\
\,\,&\sum\,t_C(x_ix_j)t_C(x_ix_j^\ast)n_C(x_l) -
t_C(x_1x_2)t_C(x_2x_3)t_C(x_3x_1) + \notag\\
\,\,&t_C(x_1x_2x_3)t_C(x_2x_1x_3) \notag
\end{align}
\emph{for all $x_1,x_2,x_3 \in C$, where both summations on the right of
\eqref{NOAS} are taken over the cyclic permutations $(ijl)$ of
$(123)$.}

\Proof Expanding
\[
n([x_1,x_2,x_3]) = n((x_1x_2)x_3) - n((x_1x_2)x_3,x_1(x_2x_3)) +
n(x_1(x_2x_3))
\]
and applying (\ref{ss.CONALT}.\ref{CP}), we conclude
\begin{align}
\label{NASL} n([x_1,x_2,x_3]) = 2n(x_1)n(x_2)n(x_3) -
n\big((x_1x_2)x_3,x_1(x_2x_3)\big).
\end{align}
Turning to the second summand on the right of \eqref{NASL}, we
obtain, by (\ref{ss.CONALT}.\ref{CPTL}),
\[
n\big((x_1x_2)x_3,x_1(x_2x_3)\big) = n(x_1x_2,x_1)n(x_3,x_2x_3) -
n\big((x_1x_2)(x_2x_3),x_1x_3\big),
\]
where applying (\ref{ss.FLECON}.\ref{ADM}) to the first summand on
the right yields
\begin{align}
\label{NAST} n\big((x_1x_2)x_3,x_1(x_2x_3)\big) =
t(x_2)^2n(x_3)n(x_1) - n\big((x_1x_2)(x_2x_3),x_1x_3\big).
\end{align}
Manipulating the expression $(x_1x_2)(x_2x_3)$ by means of
(\ref{ss.IDCON}.\ref{QUALI}) and the Moufang identities
(\ref{ss.ALTA}.\ref{MOU}), we obtain
\begin{align*}
(x_1x_2)(x_2x_3) =\,\,&(x_1x_2) \circ (x_2x_3) - (x_2x_3)(x_1x_2)
\\
=\,\,&t(x_1x_2)x_2x_3 + t(x_2x_3)x_1x_2 - n(x_1x_2,x_2x_3)1 -
x_2(x_3x_1)x_2, %\\
%=\,\,&t(x_1x_2)x_2x_3 + t(x_2x_3)x_1x_2 - n(x_1x_2,x_2x_3)1 - \\
%\,\,&n\big(x_2,(x_3x_1)^\ast\big)x_2 + n(x_2)(x_3x_1)^\ast.
\end{align*}
where (\ref{ss.CONALT}.\ref{UOP}),(\ref{ss.FLECON}.\ref{TAS})
yield
\begin{align}
\label{FOURP} (x_1x_2)(x_2x_3) =\,\,&t(x_1x_2)x_2x_3 +
t(x_2x_3)x_1x_2
- n(x_1x_2,x_2x_3)1 - \\
\,\,&t(x_1x_2x_3)x_2 + n(x_2)(x_3x_1)^\ast. \notag
\end{align}
Here we use (\ref{ss.FLECON}.\ref{NAS}),(\ref{ss.IDCON}.\ref{CJ})
to compute
\begin{align*}
n(x_1x_2,x_2x_3) =\,\,&n(x_1,x_2x_3x_2^\ast) = t(x_2)n(x_1,x_2x_3)
- n(x_1,x_2x_3x_2) %\\
%=\,\,&t(x_2)n(x_1,x_2x_3) - n\big(x_1,n(x_2,x_3^\ast)x_2\big) +
%n\big(x_1,n(x_2)x_3^\ast\big),
\end{align*}
and (\ref{ss.FLECON}.\ref{TEEN}),(\ref{ss.CONALT}.\ref{UOP}) give
\begin{align*}
n(x_1x_2,x_2x_3) =\,\,&t(x_2)n(x_1,x_2x_3) - t(x_2x_3)n(x_1,x_2) +
t(x_3x_1)n(x_2) \\
=\,\,&t(x_1)t(x_2)t(x_2x_3) - t(x_2)t(x_1x_2x_3) -
t(x_1)t(x_2)t(x_2x_3) + \\
\,\,&t(x_1x_2)t(x_2x_3) + t(x_3x_1)n(x_2),
\end{align*}
hence
\begin{align*}
n(x_1x_2,x_2x_3) = t(x_1x_2)t(x_2x_3) - t(x_2)t(x_1x_2x_3) +
t(x_3x_1)n(x_2).
\end{align*}
Inserting this into \eqref{FOURP}, and \eqref{FOURP} into the
second term on the right of \eqref{NAST}, we conclude
\begin{align*}
n\big((x_1x_2)(x_2x_3),x_1x_3\big) =\,\,&t(x_1x_2)n(x_2x_3,x_1x_3)
+ t(x_2x_3)n(x_1x_2,x_1x_3) - \\
\,\,&t(x_1x_2)t(x_2x_3)t(x_3x_1) + t(x_2)t(x_3x_1)t(x_1x_2x_3) -
\\
\,\,&t(x_3x_1)^2n(x_2) - t(x_1x_2x_3)n(x_2,x_1x_3) + \\
\,\,&n\big((x_3x_1)^\ast,x_1x_3\big)n(x_2) \\
=\,\,&t(x_1x_2)t(x_1x_2^\ast)n(x_3) + t(x_2x_3)t(x_2x_3^\ast)n(x_1) - \\
\,\,&t(x_1x_2)t(x_2x_3)t(x_3x_1) + t(x_2)t(x_3x_1)t(x_1x_2x_3) -
\\
\,\,&t(x_3x_1)^2n(x_2) - t(x_2)t(x_3x_1)t(x_1x_2x_3) + \\
\,\,&t(x_1x_2x_3)t(x_2x_1x_3) + t(x_3x_1^2x_3)n(x_2),
\end{align*}
where we may use
(\ref{ss.IDCON}.\ref{QUAD}),(\ref{ss.FLECON}.\ref{TEEN}),(\ref{ss.FLECON}.\ref{NAS})
to expand
\begin{align*}
t(x_3x_1^2x_3)n(x_2) - t(x_3x_1)^2n(x_2) =\,\,&t(x_3^2x_1^2)n(x_2)
-t(x_3x_1)^2n(x_2) \\
=\,\,&t\big([t(x_3)x_3 - n(x_3)1][t(x_1)x_1 - n(x_1)1]\big)n(x_2)
- \\
\,\,&t(x_3x_1)^2n(x_2) \\
=\,\,&t(x_3)t(x_1)t(x_3x_1)n(x_2) - t(x_1)^2n(x_2)n(x_3) - \\
\,\,&t(x_3)^2n(x_1)n(x_2) + 2n(x_1)n(x_2)n(x_3) - \\
\,\,&t(x_3x_1)^2n(x_2) \\
=\,\,&t(x_3x_1)t(x_3x_1^\ast)n(x_2) - t(x_1)^2n(x_2)n(x_3) - \\
\,\,&t(x_3)^2n(x_1)n(x_2) + 2n(x_1)n(x_2)n(x_3).
\end{align*}
Inserting the resulting expression
\begin{align*}
n\big((x_1x_2)(x_2x_3),x_1x_3\big)
=\,\,&\sum\,t(x_ix_j)t(x_ix_j^\ast)n(x_l) -
t(x_1x_2)t(x_2x_3)t(x_3x_1) + \\
\,\,&t(x_1x_2x_3)t(x_2x_1x_3) - t(x_1)^2n(x_2)n(x_3) - \\
\,\,&t(x_3)^2n(x_1)n(x_2) + 2n(x_1)n(x_2)n(x_3)
\end{align*}
into \eqref{NAST} and \eqref{NAST} into \eqref{NASL}, the theorem
follows. \hfill $\square$

\Remark The associator of a conic alternative algebra being
alternating, its norm must be totally symmetric in all three
variables. Since the expression $t_C(xy^\ast)$  is symmetric in
$x,y \in C$ by (\ref{ss.FLECON}.\ref{TEEN}), this fact is in
agreement with the right-hand side of (\ref{t.NOAS}.\ref{NOAS}).

%%%%%%%%%%%%%%%%%%%%%%%%%%%%%%%%%%%%%%%%%%%%%%%%%%
\section{Interlude: Pfister bilinear and quadratic forms in
characteristic $2$} \label{s.PFIBILQUAD}

Working over an arbitrary field $k$ of characteristic $2$, the main
purpose of this section is to collect a few results on Pfister
bilinear and Pfister quadratic forms that are hardly new but, at
least in their present form, apparently not in the literature. The
final result, Prop.~\ref{p.PFISUB}, concerns quadratic forms in
arbitrary characteristic and may be amusing even for experts.

Recall the following from, for example, \cite[8.5(iv)]{HL}. If one is
given an anisotropic $n$-Pfister bilinear form $b \cong
\pform{\alpha_1, \ldots, \alpha_n}$, then $K := k(\sqrt{\alpha_1},
\ldots, \sqrt{\alpha_n})$ is a field of dimension $2^n$ over $k$
such that $K^2 \subseteq k$.  Conversely, given such an extension
$K/k$, we can view $K$ as a $k$-vector space endowed with a
$k$-valued quadratic form $q \!: x \mapsto x^2$, and one checks that
$q$ is isomorphic to the quadratic form $v \mapsto b(v,v)$.  This
defines a bijection between isomorphism classes of extensions $K/k$
of dimension $2^n$ with $K^2 \subseteq k$ and anisotropic
$n$-quasi-Pfister quadratic forms \cite[10.4]{MR2427530}.

\subsection{Unital linear forms} \label{ss.ULF}\index{unital linear form}
We now refine this bijection. Fix an extension $K/k$ as in the
previous paragraph, and a $k$-linear form $s \!: K \rightarrow k$
such that $s(1_K) = 1$, i.e., $s$ is a retraction of the inclusion
$k \hookrightarrow K$.  (We say that $s$ is \emph{unital}.)  Define
a symmetric bilinear
\[
b_{K,s}\:K \times K \longrightarrow k \quad \text{via} \quad
b_{K,s}(u,v) := s(uv)
\]
for $u, v \in K$; it is the transfer $s_*\qform{1}$ in the notation
of \cite[\S20.A]{MR2427530}. Moreover, it is anisotropic (hence
non-degenerate) since $b_{K,s}(u,u) = u^2$ for $u \in K$.

We remark that if both $s$ and $t$ are unital linear forms on $K$,
then (by the nondegeneracy of $b_{K,s}$) there is some $u \in
K^\times$ such that $b_{K,s}(u, \text{--}) = t$, i.e., $s(uv) =
t(v)$ for all $v \in K$.

%Fix a unital $k$-algebra $K$ and a $k$-linear map $s \!: K \rightarrow k$ that is \emph{unital}, meaning that $s(1_K) = 1_k$.
%For the rest of this section, the degree of $K/k$ is
%assumed to be finite, so $[K:k] = 2^n$ for some integer $n \geq
%0$. Then $b_{K,s}$ is non-degenerate and $t\:K \to k$
%is a unital linear form if and only if there exists an element $u
%\in K$, necessarily unique, such that $s(u) = 1$ and $t(v) =
%s(uv)$ for all $v \in K$. Unital linear forms connect with Pfister
%bilinear forms in a natural way.

\subsection{Pfister bilinear forms} \label{p.BKSL}
Fix a finite extension $K/k$ such that $K^2 \subseteq k$ as assumed
above.  We compute $b_{K,s}$ for a particular $s$. Fix a 2-basis
$\sqrt{\alpha_1}, \ldots, \sqrt{\alpha_n}$ of $K$.  The monomials
$\sqrt{\alpha_1}^{i_1} \cdots \sqrt{\alpha_n}^{i_n}$ with $i_j \in
\{ 0, 1 \}$ are a $k$-basis for $K$ and we define $s$ to be 1 on
$1_k$ and $0$ on the other monomials.  One sees immediately that
$b_{K,s}$ is isomorphic to the Pfister bilinear form
$\pform{\alpha_1, \ldots, \alpha_n}$.

In fact, the $s$ constructed above is the general case.  Suppose we
are given a unital linear form $s \!: K \rightarrow k$; we will
construct a 2-basis so that $s$ is as in the previous paragraph.
Start with $A = \emptyset$ and repeat the following loop: If $k(A) =
K$, then we are done.  Otherwise, $[K:k(A)]$ is at least 2, hence
there is some $a \neq 0$ in the intersection of the $k$-subspaces
$\Ker(s)$ and $k(A)^\perp$, orthogonal complementation relative to
$b_{K,s}$. As $a$ is not in $k(A)$, $A \cup \{ a \}$ is $p$-free by
\cite[\S{V.13.1}, Prop.~3]{Bou:alg2}.  We replace $A$ by $A \cup \{
a \}$ and repeat.

In this way, we have proved: If $[K:k] = 2^n$, then $b_{K,s}$ is an
anisotropic $n$-Pfister bilinear form. Conversely, if we are given
an anisotropic symmetric bilinear form $\pform{\alpha_1, \ldots,
\alpha_n}$ over $k$, then $\sqrt{\alpha_1}, \ldots, \sqrt{\alpha_n}$
is a 2-basis for a purely inseparable extension $K/k$ as in the
coarser correspondence recalled at the beginning of the section.

\smallskip
The preceding considerations can be made more precise by looking at
the category of pairs $(K,s)$, where $K/k$ is a finite purely
inseparable field extension of exponent at most $1$, $s\:K \to k$ is
a unital linear form, and morphisms are (automatically injective)
$k$-homomorphisms of field extensions preserving unital linear
forms.  The map $(K, s) \mapsto b_{K,s}$ defines a functor from this
category into the category of anisotropic Pfister bilinear forms
over $k$ where morphisms are (automatically injective) isometries of
bilinear forms.

\begin{prop} \label{t.CHARPFIBIL}
The functor just defined is an equivalence of
categories.
\end{prop}

\begin{proof}
Given two pairs $(K,s)$, $(K^\prime,s^\prime)$ of finite purely
inseparable field extensions of exponent at most $1$ over $k$ with
unital linear forms, we need only show that any isometry
$\varphi\:b_{K,s} \to b_{K^\prime,s^\prime}$ is, in fact, a field
homomorphism preserving unital linear forms. We have
$s^\prime(\varphi(u)\varphi(v)) = s(uv)$ for all $u,v \in K$, hence
$\varphi(u)^2 = s^\prime(\varphi(u)^2) = s(u^2) = u^2$ by unitality
of $s,s^\prime$, which implies $\varphi(uv)^2 = (uv)^2 = u^2v^2 =
(\varphi(u)\varphi(v))^2$ and therefore $\varphi(uv) =
\varphi(u)\varphi(v)$ since we are in characteristic $2$. Thus
$\varphi$ is a $k$-homomorphism of fields preserving unital linear
forms in view of $s^\prime(\varphi(u)) =
s^\prime(\varphi(1_K)\varphi(u)) = s(1_Ku) = s(u)$.
\end{proof}

\subsection{The passage to Pfister quadratic forms.}
\label{ss.BILQUA} Let $\alpha \in k$ be a scalar and $s\:K \to k$ a
unital linear form.  We define the quadratic form $q_{K;\alpha,s}$
to be the transfer $s_*(\qpform{\alpha} \ot K)$, where, as usual,
$\qpform{\alpha}$ stands for the binary quadratic form given on $k
\oplus kj$ by the matrix $\left(
\begin{smallmatrix}
1 & 1 \\
0 & \alpha
\end{smallmatrix}
\right)$, so
\[
\qpform{\alpha}(\beta + \gamma j) = \beta^2 + \beta\gamma + \alpha
\gamma^2
\]
for $\beta,\gamma \in k$. By the projection formula
\cite[p.~84]{MR2427530}, $q_{K;\alpha,s}$ is isomorphic to $b_{K,s}
\ot \qpform{\alpha}$.  More concretely, on the vector space direct
sum $K \oplus Kj$ of two copies of $K$ over $k$, we can define
$q_{K;\alpha,s}$ by the formula
\begin{align}
\label{QUKA} q_{K;\alpha,s}(u + vj) := u^2 + s(uv) + \alpha v^2 &&
(u,v \in K).
\end{align}
%Example
%\ref{p.BKSL} gives:
%\subsection{Proposition.} \label{c.BAPFI} \emph{Let $a =
%(a_1,\dots,a_n) \in K^n$ be a $2$-basis of $K/k$, put $\alpha_i =
%a_i^2$ for $1 \leq i \leq n$, and let $\alpha_{n+1} \in k$. Then}
%\[
%q_{K;\alpha_{n+1},s_a} \cong
%\langle\langle\alpha_1,\dots,\alpha_n,\alpha_{n+1}\rrbracket.
%\]
%\begin{flushright}
%\vspace{-10pt} \hfill $\square$
%\end{flushright}

\begin{eg} \label{c.BAPFI}
For $K$ and $s$ as in the first paragraph of \ref{p.BKSL} and
$\alpha_{n+1} \in k$, the form $q_{K;\alpha_{n+1},s}$ is isomorphic
to $\qpform{\alpha_1, \ldots, \alpha_n, \alpha_{n+1}}$.
\end{eg}

\subsection{Theorem.} \label{c.QUEPFI} \emph{For $\alpha \in k$
and $s\:K \to k$ a unital linear form, $q_{K;\alpha,s}$ is an
$(n+1)$-Pfister quadratic form over $k$. Conversely, every}
anisotropic \emph{$(n+1)$-Pfister quadratic form over $k$ is
isomorphic to $q_{K;\alpha,s}$ for some purely inseparable field
extension $K/k$ of exponent at most $1$ and degree $2^n$, some
$\alpha \in k$ and some unital linear form $s\:K \to k$.}

%\Proof
%The first part follows immediately from
%Props.~\ref{p.FILIN},~\ref{c.BAPFI}. Conversely, suppose $q \cong
%\dla \alpha_1,\dots,\alpha_n,\alpha_{n+1}\rrbracket$ with
%$\alpha_1,\dots,\alpha_n \in k^\times$, $\alpha_{n+1} \in k$ is an
%anisotropic $(n + 1)$-Pfister quadratic form over $k$. Then the
%$n$-Pfister bilinear form $\dla \alpha_1,\dots,\alpha_n \dra$ is
%anisotropic as well, so Thm.~\ref{t.CHARPFIBIL} and
%(\ref{ss.BILQUA}.\ref{QUDEC}) apply. \hfill $\square$ \lz

\begin{proof}
Combine \ref{p.BKSL} and Example \ref{c.BAPFI}.
\end{proof}

\noindent The Pfister quadratic forms $q_{K;\alpha,s}$ of
\ref{ss.BILQUA} need not be anisotropic, nor can isometries between
two of them be described as easily as in the case of Pfister
bilinear forms (cf. Prop.~\ref{t.CHARPFIBIL}).

\subsection{Proposition.} \label{p.ISOM} \emph{Given scalars
$\alpha,\beta \in k$ and unital linear forms $s,t\:K \to k$, the
following conditions are equivalent.}
\begin{itemize}
\item [(i)] \emph{The Pfister quadratic forms $q_{K;\alpha,s}$ and
$q_{K;\beta,t}$ are isometric.} \item [(ii)] \emph{There exist
elements $u_0,v_0 \in K$ such that}
\begin{align}
\label{MUS} \beta = u_0^2 + s(u_0v_0) + \alpha v_0^2, \quad t(u) =
s(uv_0) &&\text{($u \in K$)}.
\end{align}
\end{itemize}
\Proof We identify $K \subseteq K \oplus Kj$ canonically through
the first summand. $\ssp$ \\
(i) $\Longrightarrow$ (ii). If $q := q_{K;\alpha,s}$ and $q^\prime
:= q_{K;\beta,t}$ are isometric, Witt's theorem
\cite[Thm.~8.3]{MR2427530} yields a bijective $k$-linear isometry
$\varphi\:K \oplus Kj \overset{\sim} \to K \oplus Kj$ from
$q^\prime$ to $q$ (so $q \circ \varphi = q^\prime$) that is the
identity on $K$. Then there are elements $u_0,v_0 \in K$ such that
\begin{align}
\label{VIE} \varphi(j) = u_0 + v_0j,
\end{align}
and there are $k$-linear maps $f,g\:K \to K$ such that
\begin{align}
\label{VIUV} \varphi(u + vj) = \big(u + f(v)\big) + g(v)j
&&\text{($u,v \in K$)}.
\end{align}
Combining \eqref{VIE}, \eqref{VIUV}, we conclude
\begin{align}
\label{FGUN} f(1_K) = u_0, \quad g(1_K) = v_0,
\end{align}
while evaluating $q \circ \varphi = q^\prime$ at $u + vj \in K
\oplus Kj$ with the aid of \eqref{VIUV} and \eqref{QUKA} yields
\[
t(uv) + \beta v^2 = f(v)^2 + s\big(ug(v)\big) +
s\big(f(v)g(v)\big) + \alpha g(v)^2
\]
for all $u,v\in K$. Setting $u = 0$, $v = 1_K$ and observing
\eqref{FGUN}, we obtain the first equation of \eqref{MUS}. The
second one now follows by setting $v = 1_K$ again but keeping $u
\in K$ arbitrary.

(ii) $\Longrightarrow$ (i). $s(v_0) = t(1_K) = 1$ implies $v_0 \in
K^\times$, and a straightforward verification using \eqref{MUS}
shows that the map
\[
K \oplus Kj \overset{\sim} \longrightarrow K \oplus Kj, \quad u +
vj \longmapsto (u + u_0v) + (v_0v)j,
\]
is a bijective isometry from $q^\prime$ to $q$. \hfill $\square$
\subsection{Corollary.} \label{c.FIX} \emph{Let $s\:K \to k$ be a
unital linear form. Given $\beta \in k$ and a unital linear form
$t\:K \to k$, there exists an element $\alpha \in k$ such that}
\[
q_{K;\beta,t} \cong q_{K;\alpha,s}.
\]
\begin{proof} By \ref{ss.ULF}, some $v_0 \in K$ satisfies $t(u) =
s(uv_0)$ for all $u \in K$. Now Prop.~\ref{p.ISOM} applies.
\end{proof}

\subsection{The Artin-Schreier map.} \label{ss.AS}\index{Artin Schreier map}
Let $s\:K \to k$ be a unital linear form. Then
\[
\wp_{K,s}\:K \longrightarrow k, \quad u \longmapsto \wp_{K,s}(u)
:= u^2 + s(u)
\]
is called the \emph{Artin-Schreier map} of $K/k$ relative to $s$.
It is obviously additive and becomes the usual Artin-Schreier map
(in characteristic $p = 2$), simply written as $\wp$, for $K = k$.
Given another unital linear form $t\:K \to k$, we conclude from
\ref{ss.ULF} that there is a unique element $v_0 \in K^\times$
satisfying $s(v_0) = 1$ and $t(u) = s(uv_0)$ for all $u \in K$.
Hence
\begin{align*}
\wp_{K,t}(u) = v_0^2\wp_{K,s}(uv_0^{-1}) &&\text{($u \in K$)}
\end{align*}
since the left-hand side is equal to $u^2 + s(uv_0) =
v_0^2[(uv_0^{-1})^2 + s(uv_0^{-1})]$.

With $K^\prime := \Ker(s)$, we obtain the orthogonal splitting $K
= k1_K \perp K^\prime$ relative to $b_{K,s}$, so the symmetric
bilinear form $b^\prime_{K,s} := b_{K,s}\vert_{K^\prime \times
K^\prime}$ up to isometry is uniquely determined by $b_{K,s} \cong
\langle 1 \rangle \perp b^\prime_{K,s}$. For $\alpha \in k$,
$u^\prime \in K^\prime$ we obtain
\begin{align}
\label{WPAN} \wp_{K,s}(\alpha 1_K + u^\prime) = \wp(\alpha) +
b^\prime_{K,s}(u^\prime,u^\prime) = \wp(\alpha) + u^{\prime 2}.
\end{align}

\subsection{Corollary.} \label{c.ISAS} \emph{Let $s\:K \to k$ be a
unital linear form and $\alpha,\beta \in k$.} $\ssp$ \\
(a) \emph{$q_{K;\alpha,s}$ is isotropic if and only if $\alpha \in
\Im(\wp_{K,s})$.} $\ssp$ \\
(b) \emph{$q_{K;\alpha,s}$ and $q_{K;\beta,s}$ are isometric if
and only if $\alpha \equiv \beta \bmod \Im(\wp_{K,s})$.}

\Proof (a) If $q := q_{K;\alpha,s}$ is isotropic,
(\ref{ss.BILQUA}.\ref{QUKA}) shows $u^2 + s(uv) + \alpha v^2 = 0$
for some $u,v \in K$ not both zero, which implies $v \neq 0$ and
then $\alpha = (uv^{-1})^2 + s(uv^{-1}) \in \Im(\wp_{K,s})$.
Conversely, $\alpha \in \Im(\wp_{K,s})$ implies $\alpha = u^2 +
s(u)$ for some $u \in K$, forcing $q(u + j) = 0$ by
(\ref{ss.BILQUA}.\ref{QUKA}).

(b) By Prop.~\ref{p.ISOM}, the two quadratic forms are isometric
iff (\ref{p.ISOM}.\ref{MUS}) holds with $t = s$. But this forces
$v_0 = 1_K$, and (\ref{p.ISOM}.\ref{MUS}) becomes equivalent to
$\alpha \equiv \beta \bmod \Im(\wp_{K,s})$. \hfill $\square$ \lz

\subsection{Remark.} \label{r.CLASCRI} (a) For a symmetric
bilinear form $b$ on a vector space $V$ over $k$, we adopt the
usual notation from the theory of quadratic forms by writing
$\tilde{D}(b) := \{b(v,v)\vert v \in V\}$
\cite[Def.~1.12]{MR2427530} as an additive subgroup of $k$ (recall
$\mbc(k) = 2$). With this notation, and observing
(\ref{ss.AS}.\ref{WPAN}), Cor.~\ref{c.ISAS}~(a) amounts to saying
that the quadratic form $q_{K;\alpha,s}$ is isotropic if and only
if $\alpha$ belongs to $\Im(\wp) + \tilde{D}(b^\prime_{K,s}$);
written in this way, Cor.~\ref{c.ISAS}~(a) agrees with
\cite[Lemma~9.11]{MR2427530}. $\ssp$ \\
(b) The definition of the quadratic forms $q_{K;\alpha,s}$ in
\ref{ss.BILQUA} as well as of the Artin-Schreier map in
\ref{ss.AS} make sense also if $K$ has infinite degree over $k$,
and the isomorphism $q_{K;\alpha,s} \cong b_{K,s} \ot \qpform{\alpha}$ as well as Cor.~\ref{c.ISAS}~(a)
continue to hold under these more general circumstances. \lz

\noindent We close this section with an observation that will play
a useful role in the discussion of types for composition algebras
over $2$-Henselian fields, cf. in particular \ref{r.COSLO} below.
For the remainder of this section, we dispense ourselves from the
overall restriction that $k$ have characteristic $2$.

\subsection{Proposition.} \label{p.PFISUB} \emph{Let $q$ be a Pfister
quadratic form over a field $k$ of arbitrary characteristic and
suppose $q_1,q_2$ are Pfister quadratic subforms of $q$ with
$\dim(q_1) \leq \dim(q_2)$. Then there are Pfister bilinear forms
$b_1,b_2$ over $k$ such that
\[
b_2 \otimes b_1 \otimes q_1 \cong q \cong b_2 \otimes q_2.
\]
In particular, if $\dim(q_1 ) = \dim(q_2)$, then $b \otimes q_1
\cong q \cong b \otimes q_2$ for some Pfister bilinear form $b$
over $k$.}

\Proof We may assume that $q$ is anisotropic. Writing $q\:V \to
k$, $q_i = q\vert_{V_i}\:V_i \to k$ for $i = 1,2$, with a vector
space $V$ of dimension $n$ over $k$ and subspaces $V_i \subseteq
V$ of dimension $n_i$ ($0 \leq n_1 \leq n_2 \leq n$), we then
argue by induction on $r := n - n_1$. If $r = 0$, we put $b_1 :=
b_2 := \langle 1 \rangle$. Now suppose $r > 0$. Since the $q_i$
represent $1$, there are $e_i \in V_i$ with $q(e_i) = q_i(e_i) =
1$. Now Witt's theorem \cite[Thm.~8.3]{MR2427530} yields an
orthogonal transformation $f \in O(V,q)$ with $f(e_1) = e_2$, and
replacing $q_1$ by $q_1 \circ f^{-1}\vert_{f(V_1)}\:f(V_1) \to k$
if necessary, we may assume $e_1 = e_2 \in V_1 \cap V_2$. Assuming
$n_2 < n$, and taking orthogonal complements relative to $\partial
q$, we obtain
\begin{align*}
\dim_k(V_1^\perp \cap V_2^\perp) =\,\,&\dim_k\big((V_1 +
V_2)^\perp\big) = 2^n - \dim_k(V_1 + V_2) \\
=\,\,&2^n - \dim_k(V_1) - \dim_k(V_2) + \dim_k(V_1 \cap V_2) > 2^n
- 2^{n_1} - 2^{n_2} \\
\geq\,\,&2^n - 2^{n_2+1} \geq 0.
\end{align*}
Hence we can find a non-zero vector $j \in V_1^\perp$ that also
belongs to $V_2^\perp$ if $n_2 < n$. Setting $\mu := -q(j) \in
k^\times$, we conclude from \cite[Lemma~23.1]{MR2427530} that
$q_1^\prime := \dla \mu \dra \otimes q_1$ is a subform of $q$, and
the same holds true for $q_2^\prime := \dla \mu \dra \otimes q_2$
unless $n_2 = n$, in which case we put $q_2^\prime := q_2 = q$. In
any event, the induction hypothesis yields Pfister bilinear forms
$b_1^\prime,b_2^\prime$ over $k$ such that $b_2^\prime$ is a
subform of $b_1^\prime$, $b_2^\prime = \langle 1 \rangle$ for $n_2
= n$ and $b_1^\prime \otimes q_1^\prime \cong q \cong
b_2^\prime\otimes q_2^\prime$. Hence $b_1 := b_1^\prime \otimes
\dla \mu \dra$,
\begin{align*}
b_2 :=
\begin{cases}
\langle 1 \rangle &\text{for $n_2 = n$,} \\
b_2^\prime \otimes \dla \mu \dra &\text{for $n_2 < n$}
\end{cases}
\end{align*}
are Pfister bilinear forms that satisfy $b_1 \otimes q_1 \cong
b_1^\prime \otimes \dla \mu \dra \otimes q_1 \cong b_1^\prime
\otimes q_1^\prime \cong q \cong b_2 \otimes q_2$, completing the
induction step since $b_2$ obviously is a subform of $b_1$. \hfill
$\square$ \lz
%%%%%%%%%%%%%%%%%%%%%%%%%%%%%%%%%%%%%%%%
\section{A non-orthogonal Cayley-Dickson
construction}\index{Cayley-Dickson construction!non-orthogonal}
\label{s.NOORCD} By the embedding property \ref{ss.EMPRO}, the
Cayley-Dickson construction may be regarded as a tool to recover
the structure of a composition algebra $C$ having dimension $2^n$,
$n = 1,2,3$, from a composition subalgebra $B \subseteq C$ of
dimension $2^{n-1}$: there exists a scalar $\mu \in k^\times$ with
$C \cong \mbC(B,\mu)$. In the present section, a similar
construction will be developed achieving the same objective when
$B$ is replaced by an inseparable subfield of degree $2^{n-1}$;
concerning the existence of such subfields, see \ref{ss.INSU}. The
construction we are going to present is less general but more
intrinsic than the one recently investigated by Pumpl\"un
\cite{MR2223466}. Throughout this section, we work over a field
$k$ \emph{of characteristic $2$.}

\subsection{Unital linear forms and conic algebras.} \label{ss.ULIF}
Let $K/k$ be a purely inseparable field extension of exponent at
most $1$, $C$ a conic alternative $k$-algebra containing $K$ as a
unital subalgebra and let $l \in C$. Since $K$ has trivial
conjugation by \ref{ss.INSFI} and $C$ satisfies
(\ref{ss.FLECON}.\ref{TEEN}), the relation
\begin{align}
\label{ULIF} s_l(u) := n_C(u,l) = t_C(ul)
\end{align}
holds for all $u \in K$ and defines a linear form $s_l\:K \to k$,
which is unital in the sense of \ref{ss.ULF} if $t_C(l) = 1$. In
this case, $s_l$ is called the unital linear form on $K$
\emph{associated with} $l$. \lz

\noindent The following proposition paves the way for the
non-orthogonal Cayley-Dickson construction we have in mind.

\subsection{Proposition.} \label{p.NOCOAL} \emph{Let $C$ be a
conic alternative algebra over $k$ and $K \subseteq C$ a purely
inseparable subfield of exponent at most $1$. Suppose $l \in C$
satisfies $t_C(l) = 1$, put $\mu := n_C(l) \in k$ and write $s :=
s_l$ for the unital linear form on $K$ associated with $l$. Then
$B := K + Kl \subseteq C$ is the subalgebra of $C$ generated by
$K$ and $l$. More precisely, the relations}
\begin{align}
\label{LJR} (vl)u =\,\,&s(u)v +uv + (uv)l, \\
\label{LRJ} u(vl) =\,\,&s(uv)1_C + s(u)v + s(v)u + uv + (uv)l, \\
\label{LJRJ} (v_1l)(v_2l) =\,\,&s(v_1v_2)1_C + s(v_1)v_2 \,+
s(v_2)v_1 + (1 + \mu)v_1v_2 + \\
&\big(s(v_1v_2)1_C + s(v_1)v_2 + v_1v_2\big)l \notag
\end{align}
\emph{hold for all $u,v,v_1,v_2 \in K$.}

\Proof The assertion about $B$ will follow once we have
established the relations \eqref{LJR}$\--$\eqref{LJRJ}. To do so,
we first prove
\begin{align}
\label{SYMU} u(vl) + (vl)u =\,\,&s(uv)1_C + s(v)u.
\end{align}
Since $K$ has trace zero, we combine the definition of $s$ with
(\ref{ss.IDCON}.\ref{QUALI}),(\ref{ss.FLECON}.\ref{NAS}) and
obtain
\begin{align*}
u(vl) + (vl)u =\,\,&t_C(u) vl + t_C(vl)u - n_C(u,vl)1_C \\
=\,\,&s(v)u - n_C(v^\ast u,l)1_C = s(v)u + n_C(uv,l)1_C
\\
=\,\,&s(uv)1_C + s(v)u,
\end{align*}
giving \eqref{SYMU}. In order to establish \eqref{LJR}, it
suffices to show
\begin{align}
\label{HLJR} (vx)u = t_C(ux)v + t_C(x)uv + (uv)x && (u,v \in K, x
\in C)
\end{align}
since this implies \eqref{LJR} for $x = l$. The assertion being
obvious for $u = v = 1$, we may assume $K \neq k$, forcing $k$ to
be infinite. By Zariski density, we may therefore assume that $x$
is invertible. Then the Moufang identities
(\ref{ss.ALTA}.\ref{MOU}) combine with
(\ref{ss.CONALT}.\ref{UOP}), (\ref{ss.IDCON}.\ref{QUAD}), the
right alternative law (\ref{ss.ALTA}.\ref{ALT}) and
(\ref{ss.ULIF}.\ref{ULIF}) to imply
\begin{align*}
\big((vx)u\big)x =\,\,&v(xux) = v\big(t_C(xu)x -
n_C(x)u^\ast\big) = t_C(ux)vx + n_C(x)uv \\
=\,\,&t_C(ux)vx + (uv)\big(t_C(x)x + x^2\big) = \big(t_C(ux)v +
t_C(x)uv + (uv)x\big)x,
\end{align*}
hence \eqref{HLJR}. Combining \eqref{SYMU} and \eqref{LJR}, we now
obtain \eqref{LRJ}. Finally, making use of the Moufang identities
again and of \eqref{LJR},\eqref{SYMU},(\ref{ss.ULIF}.\ref{ULIF}),
we conclude
\begin{align*}
(v_1l)(v_2l) =\,\,&(v_1l + lv_1)(v_2l) + (lv_1)(v_2l) =
\big(s(v_1)1_C +
s(1_C)v_1\big)(v_2l) + l(v_1v_2)l \\
=\,\,&s(v_1)v_2l + v_1(v_2l) + t_C\big((v_1v_2)l\big)l -
n_C(l)(v_1v_2)^\ast \\
=\,\,&s(v_1)v_2l + v_1(v_2l) + s(v_1v_2)l + \mu v_1v_2.
\end{align*}
Applying \eqref{LRJ} to the second summand on the right gives
\eqref{LJRJ}. \hfill $\square$

\subsection{Embedding inseparable field extensions into conic
algebras.} \label{ss.ENLA} We now look at the converse of the
situation described in Prop.~\ref{p.NOCOAL}. Let $K/k$ be a purely
inseparable field extension of exponent at most $1$. Suppose we
are given a unital linear form $s:K \rightarrow k$ and a scalar
$\mu \in k$. Inspired by the relations
(\ref{p.NOCOAL}.\ref{LJR})$\--$(\ref{p.NOCOAL}.\ref{LJRJ}), we now
observe the obvious fact that the vector space direct sum $K
\oplus Kj$ of two copies of $K$ carries a unique unital
non-associative $k$-algebra structure
\[
C := \mbC(K;\mu,s)
\]
into which $K$ embeds as a unital subalgebra through the first
summand such that the relations
\begin{align}
\label{FLJR} (vj)u = \,\,&\big(s(u)v + uv\big) + (uv)j, \\
\label{FLRJ} u(vj) = \,\,&\big(s(uv)1_K  + s(u)v + s(v)u +
uv\big) + (uv)j, \\
\label{FLJRJ} (v_1j)(v_2j) = \,\,&\big(s(v_1v_2)1_K  +
s(v_1)v_2 + s(v_2)v_1 + (1 + \mu)v_1v_2\big) + \\
\,\,&\big(s(v_1v_2)1_K + s(v_1)v_2 + v_1v_2\big)j \notag
\end{align}
hold for all $u,v,v_1,v_2 \in K$. Adding \eqref{FLJR} to
\eqref{FLRJ}, we obtain
\begin{align}
\label{FSYMU} u \circ (vj) = u(vj) + (vj)u = s(uv)1_K + s(v)u.
\end{align}

\subsection{Proposition.} \label{p.FOCOAL} \emph{With the notations and
assumptions of} \ref{ss.ENLA}, \emph{$C = \mbC(K;\mu,s)$ is a
non-degenerate flexible conic $k$-algebra with norm, polarized
norm, trace, conjugation respectively given by the formulas}
\begin{align}
\label{FNO} n_C(u + vj) =\,\,& u^2 + s(uv) + \mu v^2, \\
\label{FPONO} n_C(u_1 + v_1j,u_2 + v_2j) =\,\,&s(u_1v_2) +
s(u_2v_1), \\
\label{FTR} t_C(u + vj) =\,\,&s(v), \\
\label{FCONJ} (u + vj)^\ast = \,\,&s(v)1_K + u + vj
\end{align}
\emph{for all $u,u_1,u_2,v,v_1,v_2 \in K$. Moreover, $n_C =
q_{K;\mu,s}$} (\emph{cf.} \ref{ss.BILQUA}), \emph{and this is an
$(n+1)$-Pfister quadratic form if $K$ has finite degree $2^n$ over
$k$.}

\Proof The right-hand side of \eqref{FNO} defines a quadratic form
$n:C \rightarrow k$ whose polarization is given by the right-hand
side of \eqref{FPONO}. In particular, setting $t :=
\partial n(1_C,\--)$, we obtain $t(u + vj) = s(v)$ and then, using
(\ref{ss.ENLA}.\ref{FLJRJ}),(\ref{ss.ENLA}.\ref{FSYMU}),
\begin{align*}
(u + vj)^2 =\,\,&u^2 + u \circ (vj) + (vj)^2 = u^2 + s(uv)1_C
+ s(v)u + v^2 + (1+ \mu)v^2 + s(v)vj \\
=\,\,&s(v)(u + vj) + \big(u^2 + s(uv) + \mu v^2\big)1_C
\\
=\,\,&t(u + vj)(u + vj) -n(u + vj)1_C.
\end{align*}
Thus $C$ is indeed a conic $k$-algebra, and
\eqref{FNO}$\--$\eqref{FCONJ} hold. Since $K$ is a field and $s$
is unital, forcing $b_{K,s}$ to be non-degenerate by \ref{ss.ULF},
we conclude from \eqref{FPONO} that $\partial n_C$ is
non-degenerate as well. The final statement of the proposition
follows from comparing \eqref{FNO} with
(\ref{ss.BILQUA}.\ref{QUKA}) and applying Thm.~\ref{c.QUEPFI}. It
therefore remains to show that $C$ is flexible. We do so by
invoking Remark~\ref{r.EFLENAS} and verifying the first relation
of (\ref{ss.FLECON}.\ref{ADM}). Letting $u,v,w \in K$ be arbitrary
and setting $x = u + vj$, we may assume $y = w$ or $y = wj$.
Leaving the former case to the reader, we apply
\eqref{FNO},\eqref{FTR},(\ref{ss.ENLA}.\ref{FLJR}$\--$\ref{FLJRJ})
and compute
\begin{align*}
n_C(&xy,x) - n_C(x)t_C(y) = n_C\big(u(wj) +
(vj)(wj),u + vj\big) - n_C(u + vj)t_C(wj) \\
=\,\,&n_C\Big(s(uw)1_K  + s(u)w + s(w)u + uw
+ s(vw)1_K  + s(v)w + s(w)v \\
\,\,&+ (1 + \mu)vw + \big(uw + s(vw)1_K + s(v)w + vw\big)j,u +
vj\Big) +
\big(u^2 + s(uv) + \mu v^2\big)s(w) \\
=\,\,&u^2s(w) + s(u)s(vw) + s(v)s(uw) +
s(uvw) + s(v)s(uw) + s(u)s(vw) \\
\,\,&+ s(w)s(uv) + s(uvw) + s(v)s(vw) +
s(v)s(vw) + v^2s(w) + v^2s(w) \\
\,\,&+ \mu v^2s(w) + u^2s(w) + s(w)s(uv) + \mu
v^2s(w) \\
=\,\,&0,
\end{align*}
which completes the proof. \hfill $\square$ \lz

\noindent Prop.~\ref{p.NOCOAL} not only serves to illuminate the
intuitive background of the non-orthogonal Cayley-Dickson
construction presented in \ref{ss.ENLA} and Prop.~\ref{p.FOCOAL}.
It also allows for the following application.

\subsection{Proposition.} \label{p.ISCOAL} \emph{Let $C$ be a
conic alternative algebra over $k$ and $K \subseteq C$ a purely
inseparable subfield of exponent at most $1$. Suppose $l \in C$
satisfies $t_C(l) = 1$, put $\mu := n_C(l) \in k$ and write $s :=
s_l$ for the unital linear form on $K$ associated with $l$. Then
there is a unique homomorphism}
\[
\varphi:\mbC(K;\mu,s) \longrightarrow C
\]
\emph{extending the identity of $K$ and satisfying $\varphi(j) =
l$. Moreover, $\varphi$ is injective, and its image is the
subalgebra of $C$ generated by $K$ and $l$.} \Proof The uniqueness
assertion being obvious, define $\varphi:\mbC(K;\mu,\delta)
\rightarrow C$ by $\varphi(u + vj) := u + vl$ for $u,v \in K$.
Then $\varphi$ is a $k$-linear map whose image by
Prop.~\ref{p.NOCOAL} agrees with the subalgebra of $C$ generated
by $K$ and $l$. But $\varphi$ is also a homomorphism of algebras,
which follows by comparing
(\ref{ss.ENLA}.\ref{LJR})$\--$(\ref{ss.ENLA}.\ref{LJRJ}) with the
corresponding relations
(\ref{p.NOCOAL}.\ref{FLJR})$\--$(\ref{p.NOCOAL}.\ref{FLJRJ}). It
remains to show that $\varphi$ is injective, i.e., that $u + vl =
0$ for $u,v \in K$ implies $u = v = 0$. Otherwise, we would have
$v \neq 0$, which implies $0 = v^{-1}(u + vl) = v^{-1}u + l$ and
applying $t_C$ yields a contradiction. \hfill $\square$ \lz
\noindent It is a natural question to ask for conditions that are
necessary and sufficient for a non-orthogonal Cayley-Dickson
construction as in \ref{ss.ENLA} to be a composition algebra. The
answer is given by the following result.

\subsection{Theorem.} \label{t.NOCOMP} \emph{Let $K/k$ be a purely
inseparable field extension of exponent at most $1$, $s\:K \to k$
a unital linear form and $\mu \in k$ a scalar. We put $C :=
\mbC(K;\mu,s)$.} $\ssp$ \\
(a) \emph{The map $f_C\:K^4 \to k$ defined by}
\begin{align*}
f_C(u_1,u_2,u_3,u_4) := s(u_1u_2u_3u_4) + \sum\,s(u_i)s(u_ju_lu_m)
+ \sum\,s(u_iu_j)s(u_lu_4)
\end{align*}
\emph{for $u_1,u_2,u_3,u_4 \in K$, where the first (resp.~second)
sum on the right is taken over all cyclic permutations $ijlm$
(resp.~$ijl$) of $1234$ (resp.~$123$), is an alternating
quadri-linear map. Setting $\dot K := K/k1_K$, $f_C$ induces
canonically an alternating quadri-linear map $\dot f_C\:\dot K^4
\to k$. Moreover, the relation}
\begin{align}
\label{CPALT} n_C(xy) = n_C(x)n_C(y) + f_C(u_1,u_2,v_1,v_2)
\end{align}
\emph{holds for all $x = u_1 + v_1j, y = u_2 + v_2j \in C$ with
$u_i,v_i \in K$, $i = 1,2$.} $\ssp$ \\
(b) \emph{The following conditions are equivalent.}
\begin{itemize}
\item [(i)] \emph{$C$ is a composition algebra.}

\item [(ii)] $f_C = 0.$

\item [(iii)] $[K:k] \leq 4$.
\end{itemize}

\Proof (a) Setting $f := f_C$, which is obviously quadri-linear,
it is straightforward to check that it is alternating as well
(hence symmetric since we are in characteristic $2$) and satisfies
$f(1_K,u_2,u_3,u_4) = 0$ for all $u_i \in K$, $i = 1,2,3$. This
proves existence of $\dot f$ with the desired properties. It
remains to establish \eqref{CPALT}. Subtracting the first summand
on the right from the left, and expanding the resulting expression
in the obvious way, we conclude that it decomposes into the sum of
terms
\begin{align*}
n_C\big((vj)u\big) -\,\,&n_C(vj)n_C(u), \\
n_C\big(u(vj)\big) -\,\,&n_C(u)n_C(vj), \\
n_C\big((v_1j)(v_2j)\big) -\,\,&n_C(v_1j)n_C(v_2j), \\
n_C\big(u_1u_2,u_1(v_2j)\big) -\,\,&n_C(u_1)n_C(u_2,v_2j), \\
n_C(u_1u_2,(v_1j)u_2)\big) -\,\,&n_C(u_1,v_1j)n_C(u_2), \\
n_C\big(u_1(v_2j),(v_1j)(v_2j)\big) -\,\,&n_C(u_1,v_1j)n_C(v_2j), \\
n_C\big((v_1j)u_2,(v_1j)(v_2j)\big) -\,\,&n_C(v_1j)n_C(u_2,v_2j), \\
n_C\big(u_1u_2,(v_1j)(v_2j)\big) +
n_C\big(u_1(v_2j),(v_1j)u_2\big) -\,\,&n_C(u_1,v_1j)n_C(u_2,v_2j).
\end{align*}
A tedious but routine computation, involving
(\ref{p.NOCOAL}.\ref{LJR})$\--$(\ref{p.NOCOAL}.\ref{LJRJ}) and
(\ref{p.FOCOAL}.\ref{FNO}),(\ref{p.FOCOAL}.\ref{FPONO}) shows that
all but the very last one of these expressions are equal to zero.
Hence \eqref{CPALT} follows from
\begin{align*}
n_C&\big(u_1u_2,(v_1j)(v_2j)\big) +
n_C\big(u_1(v_2j),(v_1j)u_2\big) \\
=\,\,&n_C\Big(u_1u_2,\big(s(v_1v_2)1_K + s(v_1)v_2 +
s(v_2)v_1 + (1 + \mu)v_1v_2\big) + \\
\,\,&\big(s(v_1v_2)1_K + s(v_1)v_2 + v_1v_2\big)j\Big) +
\\
\,\,&n_C\Big(\big(s(u_1v_2)1_K + s(u_1)v_2 + s(v_2)u_1 +
u_1v_2\big) + (u_1v_2)j,\big(s(u_2)v_1 +
u_2v_1\big) + (u_2v_1)j\Big) \\
=\,\,&s(u_1u_2)s(v_1v_2) + s(v_1)s(u_1u_2v_2)
+ s(u_1u_2v_1v_2) + s(u_1v_2)s(u_2v_1) + \\
\,\,&s(u_1)s(u_2v_1v_2) + s(v_2)s(u_1u_2v_1) + s(u_1u_2v_1v_2) +
s(u_2)s(u_1v_1v_2) + s(u_1u_2v_1v_2) \\
=\,\,&s(u_1u_2v_1v_2) + s(u_1)s(u_2v_1v_2) + s(u_2)s(v_1v_2u_1) +
s(v_1)s(v_2u_1u_2) + s(v_2)s(u_1u_2v_1) + \\
\,\,&s(u_1u_2)s(v_1v_2) + s(u_2v_1)s(u_1v_2) +
s(v_1u_1)s(u_2v_2) + s(u_1v_1)s(u_2v_2) \\
=\,\,&f(u_1,u_2,v_1,v_2) + n_C(u_1,v_1j)n_C(u_2,v_2j).
\end{align*}
(b) (i) $\Longrightarrow$ (iii) follows from the fact that
composition algebras  have dimension at most $8$ and that the
dimension of $C$ is twice the degree of $K/k$. \\
(iii) $\Longrightarrow$ (ii) follows from the fact that $\dot K$
has dimension at most $3$, so \emph{any} alternating
quadri-linear linear map on $\dot K$ must be zero. \\
(ii) $\Longrightarrow$ (i) follows immediately from \eqref{CPALT}
since $\partial n_C$ by Prop.~\ref{p.FOCOAL} is non-degenerate.
\hfill $\square$ \lz \Remark Thm.~\ref{t.NOCOMP}~(b) can be proved
without recourse to any a priori knowledge of composition
algebras. To do so, it suffices to show directly that $[K:k] > 4$
implies $f_C \neq 0$, which can be done quite easily. We omit the
details. \lz

%%%%%%%%%%%%%%%%%%%%%%%%%%%%%%%%%%%%%%%%%%%%%%%%%%
\section{The Skolem-Noether theorem for inseparable subfields} \label{s.SKONO}
The non-orthogonal Cayley-Dickson construction introduced in the
preceding section will be applied in two ways. Recalling from
\cite{MR1763974} that every isomorphism between composition
subalgebras of a composition algebra $C$ can be extended to an
automorphism of $C$, the aim of our first application will be to
derive an analogous result where the composition subalgebras are
replaced by inseparable subfields. We begin by exploiting more
fully our description of Pfister quadratic forms presented in
Section \ref{s.PFIBILQUAD} within the framework of composition
algebras.

Throughout we continue to work over an arbitrary field $k$ of
characteristic $2$.

\subsection{Comparing non-orthogonal Cayley-Dickson
constructions.} \label{ss.CONOOR} Let $K/k$ be a purely
inseparable field extension of exponent at most $1$ and degree
$2^{n-1}$, $1 \leq n \leq 3$. Suppose we are given scalars
$\mu,\mu^\prime \in k$ and unital linear forms $s,s^\prime\:K \to
k$. By Thm.~\ref{t.NOCOMP}, $C := \mbC(K;\mu,s)$ and $C^\prime :=
\mbC(K;\mu^\prime,s^\prime)$ are composition algebras over $k$;
the notational conventions of \ref{ss.ENLA} will remain in force
for $C$ and will be extended to $C^\prime = K \oplus Kj^\prime$ in
the obvious manner. By a \emph{$K$-isomorphism} from $C$ to
$C^\prime$ we mean an isomorphism that induces the identity on
$K$; we say $C$ and $C^\prime$ are \emph{$K$-isomorphic} if a
$K$-isomorphism from $C$ to $C^\prime$ exists.

\subsection{Proposition.} \label{p.CHAPA} \emph{With the notations and
assumptions of} \ref{ss.CONOOR}, \emph{we have:} $\ssp$ \\
(a) \emph{$C$ is split if and only if $\mu \in \Im(\wp_{K,s})$.}
$\ssp$ \\
(b) \emph{$C$ and $C^\prime$ are $K$-isomorphic if and only if
there exist $u_0,v_0 \in K$ such that}
\begin{align}
\label{CHAPAM} \mu^\prime = u_0^2 + s(u_0v_0) + \mu v_0^2, \quad
s^\prime(u) = s(uv_0) &&(u \in K).
\end{align}
\Proof (a) $C$ is split iff $n_C$ is isotropic iff $\mu \in
\Im(\wp_{K,s})$ by Cor.~\ref{c.ISAS}~(a) and Prop.~\ref{p.FOCOAL}.

(b) If $C$ and $C^\prime$ are $K$-isomorphic, their norms are
isometric, so by Prop.~\ref{p.ISOM}, some $u_0,v_0 \in K$ satisfy
\eqref{CHAPAM}. Conversely, let this be so. Setting $l := u_0 +
v_0j \in C$ and combining
(\ref{p.FOCOAL}.\ref{FNO})$\--$(\ref{p.FOCOAL}.\ref{FTR})with
\eqref{CHAPAM}, we conclude $t_C(l) = s(v_0) = s^\prime(1_K) = 1$,
$n_C(l) = u_0^2 + s(u_0v_0) + \mu v_0^2 = \mu^\prime$ and
$n_C(u,l) = s(uv_0) = s^\prime(u)$ for all $u \in K$, so
Prop.~\ref{p.ISCOAL} yields a unique $K$-isomorphism $C^\prime
\overset{\sim} \to  C$ sending $j^\prime$ to $l$. \hfill $\square$

\subsection{Remark.} \label{r.GENSET} While the set-up described
in \ref{ss.CONOOR} above extends to the case of allowing purely
inseparable extensions of exponent $1$ and arbitrary degree
$2^{n-1}$, $n \geq 1$, in the obvious manner (replacing
composition algebras by flexible conic ones in the process), our
methods of proof become unsustainable in this generality. A
typical example is provided by the proof of
Prop.~\ref{p.CHAPA}~(b), where a $K$-isomorphism from $C^\prime$
to $C$ sending $j^\prime$ to $l$ in general does not exist unless
$n \leq 3$. Indeed, assuming that every element $l \in C$ of trace
$1$ allows a $K$-isomorphism
\[
\varphi\:\mbC\big(K;n_C(l),\partial n_C(\--,l)\big) \overset{\sim}
\longrightarrow C
\]
with $\varphi(j^\prime) = l$, one computes the expression
$\varphi((vj^\prime)u) - \varphi(vj^\prime)\varphi(u) = 0$ and
arrives at the conclusion that the trilinear map defined by
\begin{align}
\label{TRIALT} H(u_1,u_2,u_3) :=\,\,&s(u_1u_2u_3)1_K +
\sum\Big(s(u_iu_j)u_l + s(u_i)u_ju_l\Big) + u_1u_2u_3
\end{align}
for $u_1,u_2,u_3 \in K$ is zero, the sum on the right being
extended over all cyclic permutations $ijl$ of $123$. Since $H$ is
obviously alternating and satisfies the relations $H(1_K,u_2,u_3)
= H(u_1,u_2,u_1u_2) = 0$ for all $u_1,u_2,u_3 \in K$, it is easy
to check (see the corresponding argument in the proof of
Thm.~\ref{t.NOCOMP}) that $H$ vanishes for $n \leq 3$ (as it
should). But for $n > 3$, we may choose $u_1,u_2,u_3 \in K$ to be
$2$-independent over $k$ \cite[V~\S13.2, Thm.~2]{Bou:alg2}, forcing
the right-hand side of \eqref{TRIALT} to be a $k$-linear
combination of linear independent vectors over $k$ \cite[V~\S13.2,
Prop.~1]{Bou:alg2}, whence $H(u_1,u_2,u_3) \neq 0$.

\subsection{Proposition.} \label{t.REMDEL} \emph{Let $K/k$ be a
purely inseparable field extension of exponent at most $1$ and
degree $2^{n-1}$, $1 \leq n \leq 3$. Furthermore, let $s\:K \to k$
be a unital linear form.} $\ssp$ \\
(a) \emph{If $C$ is a composition algebra of dimension $2^n$ over
$k$ containing $K$ as a unital subalgebra, then there exists a
scalar $\mu \in k$ such that $C$ and $\mbC(K;\mu,s)$ are
$K$-isomorphic.} $\ssp$ \\
(b) \emph{For $\mu,\mu^\prime \in k$, the following conditions are
equivalent.}
\begin{itemize}
\item[(i)] \emph{$\mbC(K;\mu,s)$ and $\mbC(K;\mu^\prime,s)$ are
$K$-isomorphic.}
\item [(ii)] \emph{$\mbC(K;\mu,s) \cong \mbC(K;\mu^\prime,s)$.}
\item [(iii)] $\mu \equiv \mu^\prime \bmod \Im(\wp_{K,s})$.
\end{itemize}
\Proof (a) Pick any $l \in C$ of trace $1$. Then $C = K \oplus
Kl$, and Prop.~\ref{p.ISCOAL} yields a $K$-isomorphism $C
\overset{\sim} \to \mbC(K;\mu^\prime,s^\prime)$, for some
$\mu^\prime \in k$ and some unital linear form $s^\prime\:K \to
k$. Hence there is a unique $v_0 \in K^\times$ such that
$s^\prime(u) = s(uv_0)$ for all $u \in K$. Setting $\mu :=
v_0^{-2}\mu^\prime + \wp_{K,s}(u_0v_0^{-1})$ and consulting
Prop.~\ref{p.CHAPA}~(b), we obtain a $K$-isomorphism
$\mbC(K;\mu^\prime,s^\prime) \overset{\sim} \to \mbC(K;\mu,s)$.

(b) While (i) $\Rightarrow$ (ii) is obvious, (ii) $\Rightarrow$
(iii) follows from Cor.~\ref{c.ISAS}~(b). It remains to show (iii)
$\Rightarrow$ (i). But (iii) implies $\mu^\prime = \mu +
\wp_{K,s}(u_0)$ for some $u_0 \in K$, so
(\ref{p.CHAPA}.\ref{CHAPAM}) holds with $v_0 =1_K$. This gives
(i). \hfill $\square$ \lz

\noindent Every element of an octonion algebra over a field is
contained in a suitable quaternion subalgebra
\cite[Prop.~1.6.4]{MR1763974}. However, it doesn't seem entirely
obvious that, if the octonion algebra is split, the quaternion
subalgebra can be chosen to be split as well. But, in fact, it
can:
\subsection{Proposition.} \label{c.EMBQUAT} \emph{Let $C$ be a split
octonion algebra over an arbitrary field $F$. Then every element
of $C$ is contained in a split quaternion subalgebra of $C$.}
\Proof Let $x \in C$. We may assume that $R := k[x]$ has dimension
$2$ over $k$. There are three cases \cite[III \S2 Prop.~3]{MR0354207}. $\ssp$ \\
\emph{Case $1$.} $R$ is the algebra of dual numbers. In
particular, it contains zero divisors. Hence so does every
quaternion subalgebra of $C$ containing $x$, which therefore must
be split. $\ssp$ \\
\emph{Case $2$.} $R$ is (quadratic) \'etale. Since $C$ up to
isomorphism is uniquely determined by splitness, it may be
obtained from $R$ by the Cayley-Dickson process as $C \cong
\mbC(R;1,1)$, which implies that $x$ is contained in the split
quaternion subalgebra $\mbC(R,1)$ of $C$. $\ssp$ \\
\emph{Case $3$.} We are left with the most delicate possibility
that $F = k$ has characteristic $2$ and $K :=R$ is a purely
inseparable field extension of $k$ having exponent at most $1$.
Let $c \in C$ be an idempotent different from $0,1_C$, which
exists since $C$ is split, and write $B$ for the subalgebra of $C$
generated by $K$ and $c$. Since $c$ has trace $1$,
Prop.~\ref{p.ISCOAL} yields a scalar $\mu \in k$, a unital linear
form $s\:K \to k$ and an isomorphism from $B^\prime :=
\mbC(K;\mu,s)$ onto $B$. But $B^\prime$ is a quaternion algebra by
Thm.~\ref{t.NOCOMP} whereas $B$, containing $c$, has zero
divisors. Hence $B$ is a split quaternion subalgebra of $C$
containing $x$. \hfill $\square$ \lz
\noindent And finally, we need a purely technical result.
\subsection{Lemma.} \label{l.BIGINS} \emph{Let $C$ be an octonion
division algebra over $k$ and suppose $\varphi\:K_1 \overset{\sim}
\to K_2$ is an isomorphism of inseparable quadratic subfields
$K_1,K_2 \subseteq C$. Then there exist inseparable subfields
$L_1,L_2 \subseteq C$ of degree $4$ over $k$ such $K_i \subseteq
L_i$ for $i = 1,2$ and $\varphi$ extends to an isomorphism
$\psi\:L_1 \overset{\sim} \to L_2$.}

\Proof Given any elements $y \in C$, $x_i \in K_i \setminus k 1_C$
for $i = 1,2$, denote by $L_i$ the subalgebra of $C$ generated by
$K_i$ and $y$. Since $C$ has no zero divisors, $L_i$ is either a
composition algebra or an inseparable field extension, the latter
possibility being equivalent to the trace of $C$ vanishing
identically on $L_i$. In any event, the dimension of $L_i$ is
either $2$ or $4$. Moreover, $L_i$ is spanned by $1_C,x_i,y,x_iy$
as a vector space over $k$. Summing up, $L_i/k$ is therefore an
inseparable field extension of degree $4$ if and only if $y \notin
K_i$ satisfies the condition $t_C(y) = t_C(x_iy) = 0$. To choose
$y$ appropriately, we now write $C^0$ for the space of trace zero
elements in $C$ and consider the hyperplane intersection $V := C^0
\cap x_1C^0 \cap x_2C^0 \subseteq C$, which is a subspace of
dimension at least $5$. Since a group cannot be the union of two
proper subgroups, we conclude $V \nsubseteq K_1 \cup K_2$ and,
accordingly, pick a $y \in V$ that neither belongs to $K_1$ nor to
$K_2$. Then $y$ has trace zero and $y = x_iz_i$ for some $z_i \in
C^0$, so $x_iy = x_i^2z_i = n_C(x_i)z_i \in k z_i$ has trace zero
as well, and by the above, $L_i/k$ is an inseparable field
extension of degree $4 $ containing $K_i$. Moreover, $K_1 \cong
K_2$ implies $K_1^2 =K_2^2$, hence
\[
L_1^2 = K_1^2 + K_1^2y^2 = K_2^2 + K_2^2y^2 = L_2^2,
\]
so there is an isomorphism $\psi:L_1 \overset{\sim} \to L_2$,
which necessarily extends $\varphi$ since the only $k$-embedding
$K_1 \to L_2$ is the one induced by $\varphi$. \hfill $\square$
\lz After these preparations, we can now establish the
Skolem-Noether theorem for inseparable subfields of composition
algebras.

\subsection{Theorem.} \label{t.SKONOE} (Skolem-Noether) \emph{Let
$C$ be a composition algebra over $k$. Then every isomorphism
between inseparable subfields of $C$ can be extended to an
automorphism of $C$.}

\Proof Write $\mbox{dim}_k(C) = 2^n$, $0 \leq n \leq 3$ and let
$\varphi\:K_1 \overset{\sim} \to K_2$ be an isomorphism of
inseparable subfields $K_1,K_2 \subseteq C$ having degree
$2^{n^\prime}$, $0 \leq n^\prime < n$, over $k$. We may assume
$n^\prime > 0$ and first reduce to the case $n^\prime = n - 1$. To
do so, suppose the theorem holds for $n^\prime = n - 1$ and let
$n^\prime < n - 1$. Then $n^\prime = 1$, $n = 3$, so $C$ is an
octonion algebra containing $K_1,K_2$ as inseparable quadratic
subfields. If $C$ is split, there are split quaternion subalgebras
$B_i \subseteq C$ containing $K_i$ for $i = 1,2$
(Prop.~\ref{c.EMBQUAT}). But $B_1,B_2$ are isomorphic, hence
conjugate under the automorphism group of $C$, by the classical
Skolem-Noether theorem for composition algebras. Hence up to
conjugation by automorphisms of $C$, we may assume $B_1 = B_2 =:
B$. But then $\varphi$ extends to an automorphism of $B$, which in
turn extends to an automorphism of $C$. We are left with the case
that $C$ is a division algebra. By Lemma~\ref{l.BIGINS}, there are
inseparable subfields $K_i \subseteq L_i \subseteq C$ of degree
$4$ ($i = 1,2$) such that $\varphi$ extends to an isomorphism
$\psi\:L_1 \overset{\sim} \to L_2$. But $\psi$ in turn extends to
an automorphism of $C$, completing the reduction to the case
$n^\prime = n - 1$. From now on we assume $n^\prime = n - 1$ and
fix a unital linear form $s_2:K_2 \to k$. Then $s_1 := s_2 \circ
\varphi\:K_1 \to k$ is a unital linear form as well. For $i =
1,2$, Prop.~\ref{t.REMDEL} (a) yields scalars $\mu_i \in k$ and
$K_i$-isomorphisms $\psi_i\:\mbC(K_i;\mu_i,s_i) \overset{\sim} \to
C$. Now observe that the non-orthogonal Cayley-Dickson
construction of \ref{ss.ENLA} is functorial in the parameters
involved. Hence $\varphi$ determines canonically an isomorphism
\[
\psi := \mbC(\varphi)\:\mbC(K_1;\mu_1,s_1) \overset{\sim}
\longrightarrow \mbC(K_2;\mu_1,s_1 \circ \varphi^{-1}) =
\mbC(K_2;\mu_1,s_2).
\]
Putting things together, we thus obtain an isomorphism
\[
\psi_2^{-1} \circ \psi_1 \circ \psi^{-1}:\mbC(K_2;\mu_1,s_2)
\overset{\sim} \longrightarrow \mbC(K_2;\mu_2,s_2).
\]
Applying Prop.~\ref{t.REMDEL} (b), we see that there exists a
$K_2$-isomorphism
\[
\chi\:\mbC(K_2;\mu_2,s_2) \overset{\sim} \to \mbC(K_2;\mu_1,s_2),
\]
giving rise to the automorphism $\phi := \psi_2 \circ \chi^{-1}
\circ \psi \circ \psi_1^{-1}$ of $C$, which extends $\varphi$
since $\psi_2 = \chi^{-1} = \Eins$ on $K_2$ and $\psi = \varphi$,
$\psi_1^{-1} = \Eins$ on $K_1$. \hfill $\square$

\Remark Let $C$ be an octonion algebra over $k$ and $K \subseteq
C$ an inseparable subfield of degree $4$. Changing scalars to the
algebraic closure, $\bar k$, of $k$, $K \otimes_k \bar k$ becomes
a unital $\bar k$-subalgebra of $\bar C := C \otimes_k \bar k$
containing a $3$-dimensional subalgebra $N$ that consists entirely
of nilpotent elements. Hence $N \subseteq \bar C$ is a Borel
subalgebra in the sense of \cite{MR52:8209}. The fact that all
Borel subalgebras of $\bar C$ are conjugate under its automorphism
group \cite[\S~2,~1.]{MR52:8209} corresponds nicely with the
Skolem-Noether theorem. \lz

\noindent It is a natural question to ask how our non-orthogonal
Cayley-Dickson construction can be converted into the classical
orthogonal one. When dealing within the framework of composition
(and not of arbitrary conic) algebras, here is a simple answer.

\subsection{Orthogonalizing the non-orthogonal Cayley-Dickson
construction.} \label{ss.NOOROR} Let $K/k$ be a purely inseparable
field extension of exponent at most $1$ and degree $2^n$, $n =
1,2$. Suppose we are given an intermediate subfield $k \subseteq
K^\prime \subseteq K$ of degree $2^{n-1}$, a scalar $\mu \in k$
and a unital linear form $s\:K \to k$. Then $s^\prime :=
s\vert_{K^\prime}\:K^\prime \to k$ is a unital linear form on
$K^\prime$ and $B := \mbC(K^\prime;\mu,s^\prime)$ is a composition
subalgebra of $C := \mbC(K;\mu,s)$. Moreover,
(\ref{p.FOCOAL}.\ref{FPONO}) implies
\begin{align}
\label{PERP} B^\perp = K^{\prime\perp} \oplus K^{\prime\perp} j,
\end{align}
orthogonal complementation in $K$ (resp.~$C$) being taken relative
to $b_{K,s}$ (resp.~$\partial n_C$). From \eqref{PERP} we conclude
\begin{align}
\label{CEBEU} C = \mbC(B,-n_C(u)) = \mbC(B;u^2)
\end{align}
for any non-zero element $u \in K^{\prime\perp}$.

To be more specific, let $C$ be an octonion algebra over $k$ and
$K \subseteq C$ an inseparable subfield of degree $4$. Pick a
$2$-basis $a = (a_1,a_2)$ of $K/k$. By Prop.~\ref{t.REMDEL}~(a),
there exists a scalar $\mu \in k$ with $C \cong \mbC(K;\mu,s_a)$.
On the other hand,
\[
L = \mbC(k;\mu,\Eins_k) = k[\bft]/(\bft^2 + \bft + \mu)
\]
is a quadratic \'etale $k$-algebra, and \eqref{CEBEU} yields
\[
C = \mbC(L;a_1^2,a_2^2)
\]
as an ordinary Cayley-Dickson process starting from $L$.

%%%%%%%%%%%%%%%%%%%%%%%%%%%%%%%%%%%%%%%%%%%%%
\section{Conic division algebras in characteristic $2$.} \label{s.CODI}

We now turn to a second application of the non-orthogonal
Cayley-Dickson construction which consists in finding new examples
of conic division algebras. A few comments on the historical
context seem to be in order.

\subsection{Conic division algebras over arbitrary fields.}
\label{ss.CODI} Among all conic division algebras, it is the
composition division algebras that are particularly well
understood and particularly easy to construct: by \ref{ss.NOCRI},
it suffices to ensure that their norms be anisotropic. Of course,
composition division algebras exist only in dimensions $1,2,4,8$,
as do all non-associative division algebras over the reals, by the
Bott-Kervaire-Milnor theorem \cite[Kap.~10,\S~2]{Numbers}; in
particular, the Cayley-Dickson process~\ref{ss.CDP} leads to conic
algebras
\begin{align*}
\mbC(\IR;\mu_1,\dots,\mu_n), \quad \mu_1 = \cdots = \mu_n = -1
&&(n \in \IZ,\,\,n \geq 1)
\end{align*}
over the reals whose norms are positive definite (hence anisotropic)
but which fail to be division algebras unless $n \leq 3$. Hence it
is natural to ask for examples of conic division algebras in
dimensions other than $1,2,4,8$, over fields other than the reals.

From the point of view of non-associative algebras, conic division
algebras that are not central, like purely inseparable field
extensions of characteristic $2$ and exponent $1$ as in
\ref{ss.INSFI}, are not particularly interesting. Over appropriate
fields of characteristic not $2$, the first examples of central
conic division algebras in all dimensions  $2^n$, $n = 0,1,2\dots$
are apparently due to Brown
\cite[pp.~421-422]{MR0215891}. They all arise from the
base field, an iterated Laurent series field in finitely many
variables, by the Cayley-Dickson process;
generalizations of these examples will be discussed in
Example~\ref{e.LAPR} below. Other examples of central conic division
algebras in dimension $16$, using a refinement of the Cayley-Dickson
construction, have been exhibited by Becker
\cite[Satz~16]{MR0320099}. Examples of central
\emph{commutative} conic division algebras in characteristic $2$ and
all dimensions $2^n$, $n \geq 0$, have been constructed by Albert
\cite[Thm.~2]{MR0047027}. These algebras are closely related to
purely inseparable field extensions of exponent $1$ since their
norms bilinearize to zero and hence degenerate when extending
scalars to the algebraic closure.

%Using our non-orthogonal Cayley-Dickson construction, we
%can do better than that by exhibiting flexible central conic
%division algebras whose norms, in addition to the properties
%listed above, are non-singular quadratic forms (hence continue to
%stay that way under arbitrary base change).

In view of the preceding results one may ask whether the dimension
of a finite-dimensional conic division algebra is always a power of
$2$. Though a feeble result along these lines has been obtained by
Petersson \cite{MR46:9132}, the answer to this question doesn't seem
to be known.

In this paper, two classes of conic division algebras
in all dimensions $2^n$, $n = 0,1,2\dots$, will be constructed. The
examples of the first class, to be discussed in the present section,
depend strongly on the non-orthogonal Cayley-Dickson construction,
hence exist only in characteristic $2$ but differ from Albert's by
being central and highly non-commutative and by allowing
\emph{arbitrary} anisotropic Pfister quadratic forms as their norms,
which in particular remain non-singular under all scalar extensions.
The second class of examples will be discussed in \ref{e.LAPR} and
\ref{e.LAU} below.

\subsection{Notations and conventions.} \label{ss.NOCO} For the
remainder of this section, we fix a base field $k$ of
characteristic $2$, a purely inseparable field extension $K/k$ of
exponent at most $1$, a scalar $\mu \in k$ and a unital linear
form $s\:K \to k$ to consider the non-orthogonal Cayley-Dickson
construction $C :=\mbC(K;\mu,s)$ as in \ref{ss.ENLA}. We put
$[K:k] = 2^n$, $n = 0,1,2,\dots$ and explicitly allow the
possibility $n = \infty$, i.e., that $K$ has infinite degree over
$k$. The following proposition paves the way for the application
we have in mind.

\subsection{Proposition.} \label{p.GENOOR} \emph{With the
notations and conventions of} \ref{ss.NOCO}, \emph{the following
assertions hold.} $\ssp$ \\
(a) \emph{$C$ is locally finite-dimensional.} $\ssp$ \\
(b) \emph{For $n \geq 2$, every element of $C$ belongs to an
octonion subalgebra of $C$.} $\ssp$ \\
(c) \emph{$C$ is central simple for $n \geq 1$ and has trivial
nucleus for $n \geq 2$.} $\ssp$
%(d) \emph{$n_C(xy) = n_C(yx)$ for all $x,y \in C$.}

\Proof (a) It suffices to note that finitely many elements $x_i =
u_i + v_ij \in C$ with $u_i,v_i \in K$ $(1 \leq i \leq m)$ are
contained in $\mbC(K^\prime;\mu,s\vert_{K^\prime})$ ($K^\prime
:=k(u_1,v_1,\dots,u_m,v_m)$), which is a subalgebra of $C$ having
dimension at most $2^{2m+1}$.

(b) Let $x = u + vj \in C$, $u,v \in K$. Then there is a subfield
$K^\prime \subseteq K$ containing $u,v$ and having degree $4$ over
$k$, so $C^\prime := \mbC(K^\prime;\mu,s\vert_{K^\prime})
\subseteq C$ by Thm.~\ref{t.NOCOMP} is an octonion subalgebra
containing $x$.

(c) Standard properties of composition algebras allow us to assume
$n > 2$. Combining Props.~\ref{p.ADSIM},~\ref{p.FOCOAL} we see
that $C$ is simple. Let $x \in \mbN(C)$ and apply (b) to pick an
octonion subalgebra $C^\prime \subseteq C$ containing $x$. Since
$C^\prime$ has trivial nucleus \cite[Prop.~1.9.2]{MR1763974}, we
conclude $x \in k1_C$, so $C$ has trivial nucleus as well; in
particular, $C$ is central. \hfill $\square$ \lz

%(d) follows immediately from Thm.~\ref{t.NOCOMP}~(a).

\noindent Referring the reader to our version of the
Artin-Schreier map (\ref{ss.AS}, Remark~\ref{r.CLASCRI}), we can
now state the main result of this section.

\subsection{Theorem.} \label{t.ADCODI} \emph{With the notations
and conventions of} \ref{ss.NOCO}, \emph{the following conditions
are equivalent.}
\begin{itemize}
\item [(i)] \emph{$C$ is a division algebra.}

\item [(ii)] \emph{$n_C$ is anisotropic.}

\item [(iii)] \emph{$\mu \notin \Im(\wp_{K,s})$.}
\end{itemize}

\Proof The implication (i) $\Rightarrow$ (ii) follows from
\ref{p.NODIV}, while (ii) and (iii) are equivalent by
Cor.~\ref{c.ISAS}~(a) and Remark~\ref{r.CLASCRI}~(b). It remains
to prove \\
(ii) $\Longrightarrow$ (i). Since $C$ is locally
finite-dimensional by Prop.~\ref{p.GENOOR}~(a), it suffices to
show that there are no zero divisors, so suppose $x_1,x_2 \in C$
satisfy $x_1x_2 = 0$. By (ii) and (\ref{p.NOCOMT}.\ref{NOCO}),
this implies $x_2x_1 = 0$, and from (\ref{ss.IDCON}.\ref{QUALI})
we conclude $0 = x_1x_2 + x_2x_1 = t_C(x_1)x_2 + t_C(x_2)x_1 -
n_C(x_1,x_2)1_C$. If $t_C(x_1) \neq 0$, this yields $x_2 \in
k[x_1]$, hence $x_2 = 0$ since $n_C$ being anisotropic implies
that $k[x_1]$ is a field. By symmetry, we may therefore assume
$t_C(x_1) = t_C(x_2) = n_C(x_1,x_2) = 0$. Write $x_i = u_i + v_ij$
with $u_i,v_i \in K$ for $i = 1,2$. Then $s(v_i) = 0$ by
(\ref{p.FOCOAL}.\ref{FTR}), $s(u_1v_2) = s(u_2v_1)$ by
(\ref{p.FOCOAL}.\ref{FPONO}), and if $v_1 = 0$ or $v_2 = 0$,
Thm.~\ref{t.NOCOMP}~(a) yields $n_C(x_1)n_C(x_2) = n_C(x_1x_2) =
0$, hence $x_1 = 0$ or $x_2 = 0$. We are thus reduced to the case
\begin{align}
\label{REFO} v_1 \neq 0 \neq v_2, \quad s(v_1) = s(v_2) = 0, \quad
s(u_1v_2) = s(u_2v_1).
\end{align}
Next we use (\ref{ss.ENLA}.\ref{FLJR}$\--$\ref{FLJRJ}) and
\eqref{REFO} to expand $(u_1 + v_1j)(u_2 + v_2j) = 0$. A short
computation gives
\begin{align}
u_1u_2 + s(u_1v_2)1_K + s(u_1)v_2 + u_1v_2 + s(u_2)v_1 &+ u_2v_1 +
s(v_1v_2)1_K + (1 + \mu)v_1v_2 = 0, \notag \\
\label{ZERUV} u_1v_2 + u_2v_1 + s(v_1v_2)1_K &+ v_1v_2 = 0.
\end{align}
Adding these two relations, we obtain
\begin{align}
\label{ZERUU} u_1u_2 + s(u_1v_2)1_K + s(u_1)v_2 + s(u_2)v_1 + \mu
v_1v_2 = 0
\end{align}
and note that \eqref{ZERUV},\eqref{ZERUU} are symmetric in the
indices $1,2$ by \eqref{REFO}. We now claim:
\begin{itemize}
\item [($\ast$)] The subfield $K^\prime$ of $K/k$ generated by
$u_1,u_2,v_1,v_2$ is spanned by
\begin{align}
\label{LINSP} 1_K,u_1,u_2,v_1,v_2,u_1v_2,u_2v_1
\end{align}
as a vector space over $k$
\end{itemize}
Suppose for the time being that this claim has been proved. Then
the field extension $K^\prime/k$ has degree at most $7$. Being
purely inseparable at the same time, it has, in fact, degree at
most $4$. By Thm.~\ref{t.NOCOMP}~(b), $C^\prime :=
\mbC(K^\prime;\mu,s\vert_{K^\prime}) \subseteq C$ is therefore a
composition subalgebra containing $x_1,x_2$ and inheriting its
anisotropic norm from $C$. Thus $C^\prime$ is a division algebra,
and $x_1x_2 = 0$ implies $x_1 = 0$ or $x_2 = 0$, as desired.

We are thus reduced to showing ($\ast$). Writing $V$ for the
linear span of the vectors in \eqref{LINSP}, it suffices to show
that $V \subseteq K$ is a $k$-subalgebra, i.e., that the product
of any two distinct elements in \eqref{LINSP} belongs to $V$.
Since $v_1v_2 \in V$ by \eqref{ZERUV}, hence $u_1u_2 \in V$ by
\eqref{ZERUU}, this will follow once we have shown that
\begin{align*}
v_p^2u_qv_q, \quad u_pu_qv_p, \quad u_pv_pv_q &&(\{p,q\} =
\{1,2\})
\end{align*}
all belong to $V$. By symmetry, we may assume $p = 1$, $q = 2$.
Multiplying \eqref{ZERUV} by $u_2v_1$, we obtain
\[
u_1u_2v_1v_2 + u_2^2v_1^2 + s(v_1v_2)u_2v_1 + v_1^2u_2v_2 = 0,
\]
forcing $v_1^2u_2v_2 \equiv u_1u_2v_1v_2 \bmod V$. But multiplying
\eqref{ZERUU} by $v_1v_2$ implies $u_1u_2v_1v_2 \in V$, so we have
$v_1^2u_2v_2 \in V$ as well. Moreover, multiplying \eqref{ZERUV}
first by $u_1$, then by $v_1$, yields
\begin{align*}
u_1u_2v_1 = u_1^2v_2 + s(v_1v_2)u_1 + u_1v_1v_2 \equiv v_1^2u_2 +
s(v_1v_2)v_1 + v_1^2v_2 \equiv 0 \bmod V.
\end{align*}
Hence also $u_1v_1v_2 = u_1^2v_2 + u_1u_2v_1 + s(v_1v_2)u_1 \in
V$, which completes the proof. \hfill $\square$

\subsection{Corollary.} \label{c.PFICODI} \emph{For a Pfister
quadratic form $q$ over a field of characteristic $2$ to be the
norm of a conic division algebra it is necessary and sufficient
that $q$ be anisotropic.} \hfill $\square$

\subsection{Examples.} \label{e.CODIAL} Letting $k$ be the field
of rational functions in countably many variables
$\bft_1,\bft_2,\bft_3,\dots$ over any field of characteristic $2$,
e.g., over $\IF_2$, the $(n + 1)$-Pfister quadratic forms $\dla
\bft_1,\dots,\bft_n,\bft_{n+1}\rrbracket$, $n \geq 0$, by standard
arguments are easily seen to be anisotropic, and we obtain central
flexible conic division algebras over $k$ in all dimensions $2^n$,
$n = 0,1,2,\dots$.

%%%%%%%%%%%%%%%%%%%%%%%%%%%%%%%%%%%%%%%%%%%%%%%%%%
\part{$2$-Henselian base fields} \label{hensel.part}
%We now pass from the purely algebraic set-up of the previous
%sections to a more arithmetic one by working over a field $F$ that
%is \emph{Henselian} with respect to a discrete valuation. Our
%ultimate goal will be to understand composition algebras over $F$.
%To this end, we first note

%The methods we are going to develop in order to achieve the
%aforementioned objectives are not confined to composition
%algebras; up to a point, they work more generally for pointed
%quadratic spaces that are round and anisotropic. Moreover, our
%approach lends itself to the study of what we call
%$\lambda$-composition algebras, a class of conic division algebras
%over $F$ generalizing ordinary composition algebras and, as we
%shall see, existing in all dimensions $2^n$, $n = 0,1,2,\dots$ .

%%%%%%%%%%%%%%%%%%%%%%%%%%%%%%%%%%%%%%%%%%%%%%%%%%
\section{Pointed quadratic spaces over $2$-Henselian fields.} \label{s.COALHE}
In this section, we recast the conceptual foundations for the
study of quadratic forms over Henselian fields in the setting of
pointed quadratic spaces. Our subsequent considerations also fit
naturally into the valuation theory of Jordan division rings
\cite{MR50:422} when specialized to the Jordan algebras of pointed
quadratic spaces over Henselian fields.

\subsection{Round quadratic forms.} \label{ss.ROQUA} For the time
being, we work over a field $k$ that is completely arbitrary. We
recall from \cite[\S~9A, p.~52]{MR2427530} that a
finite-dimensional quadratic form $q$ over $k$ is said to be
\emph{round} if all its non-zero values are precisely its
similarity factors; in particular, they form a subgroup of
$k^\times$. The most important examples of round quadratic forms
are Pfister forms and quasi-Pfister forms \cite[Cor.~9.9,
Cor.~10.13]{MR2427530}.

\subsection{Pointed quadratic spaces over arbitrary
fields.}\index{pointed quadratic space} \label{ss.POQUA} Adopting
the terminology of Weiss \cite[Def.~1.1]{MR2177056}, by a
\emph{pointed quadratic space} over $k$ we mean a triple $Q =
(V,q,e)$ consisting of a finite-dimensional vector space $V$ over
$k$, a quadratic form $q\:V \to k$ and a vector $e \in V$ which is
a \emph{base point} for $q$ in the sense that $q(e) = 1$.
Morphisms of pointed quadratic spaces are isometries preserving
base points. $Q$ is said to be \emph{non-singular} (resp.
\emph{anisotropic}, \emph{round}, \emph{Pfister}, $\dots$) if $q$
is. Given a conic algebra $C$ over $k$, we obtain in $Q_C :=
(C,n_C,1_C)$ a pointed quadratic space and every pointed quadratic
space arises in this manner (Loos \cite{Loos11}, see also
Rosemeier \cite{Rosemeier02}). Notationally, we do not always
distinguish carefully between $C$ and $Q_C$. If $Q = (V,q,e)$ is
any pointed quadratic space over $k$, we therefore find it
convenient to put $V_Q = V$ as a vector space over $k$ and to call
$n_Q := q$ the \emph{norm}\index{pointed quadratic space!norm},
$1_Q := e$ the \emph{unit element}\index{pointed quadratic
space!unit element}, $t_Q := \partial n_Q(1_Q,\--)$ the
\emph{trace}\index{pointed quadratic space!trace}, $\iota_Q\:Q \to
Q,\,x \mapsto x^\ast := t_Q(x)1_Q - x$ the
\emph{conjugation}\index{pointed quadratic space!conjugtion} of
$Q$. We also put $V_Q^\times := \{x \in Q \mid n_Q(x) \neq 0\}$.

When dealing with quadratic forms representing $1$, insisting on
pointedness is not so much a matter of necessity but one of
convenience, making the language of non-associative algebras the
natural mode of communication. Just as for composition algebras,
Witt's theorem \cite[Thm.~8.3]{MR2427530} implies that pointed
quadratic spaces are classified by their norms:

\subsection{Proposition.} \label{p.PONO} \emph{For non-singular
pointed quadratic spaces over $k$ to be isomorphic it is necessary
and sufficient that their underlying quadratic forms be isometric.}
\hfill $\square$

\subsection{Enlargements.} \label{ss.ENLAR} Let $P$ be a pointed
quadratic space over $k$ and $\mu \in k^\times$. Writing $V_P
\oplus V_Pj$ for the direct sum of two copies of the vector space
$V_P$ over $k$ as in \ref{ss.CDCO}, and identifying $V_P \subseteq
V_P \oplus V_Pj$ as a subspace through the first summand,
\begin{align*}
Q := \dla\mu\dra \otimes P := (V_Q,n_Q,1_Q),\quad V_Q := V_P
\oplus V_Pj, \quad n_Q := \dla\mu\dra \otimes n_P, \quad 1_Q :=
1_P
\end{align*}
with
\begin{align}
\label{ENQU} \big(\dla\mu\dra \otimes n_P\big)(u + vj) = n_P(u) -
\mu n_P(v) &&(u,v \in V_P)
\end{align}
is again a pointed quadratic space over $k$ whose trace and
conjugation are given by the formulas
\begin{align}
\label{TEQU} t_Q(u + vj) =\,\,& t_P(u), \\
\label{CEQU} (u + vj)^\ast =\,\,&u^\ast - vj
\end{align}
for all $u,v \in V_P$. Moreover, if $P$ is round (resp. Pfister),
so is $Q$. Finally, a comparison with the Cayley-Dickson
construction \ref{ss.CDCO} shows $Q_{\mbC(B,\mu)} = \dla\mu\dra
\otimes Q_B$ for any conic $k$-algebra $B$ and any $\mu \in
k^\times$. The following two statements are standard facts about
round quadratic forms, translated into the setting of pointed
quadratic spaces; we refer to \cite[Prop.~9.8,
Lemma~23.1]{MR2427530} for details.

\subsection{Proposition.} \label{p.CDPOQUA} \emph{Let $P$ be
a non-singular round pointed quadratic space over $k$.} \\
(a) \emph{For $\mu \in k^\times$, the following conditions are
equivalent.}
\begin{itemize}
\item [(i)] \emph{$\dla\mu\dra \otimes P$ is isotropic.}

\item [(ii)] \emph{$\mu \in n_P(V_P^\times)$.}

\item [(iii)] \emph{$\dla\mu\dra \otimes P$ is hyperbolic.}
\end{itemize}
(b) \emph{Let $\mu_1,\mu_2 \in k^\times$. Then $\dla\mu_1\dra
\otimes P \cong \dla\mu_2\dra \otimes P$ if and only if $\mu_1 =
\mu_2n_P(u)$ for some $u \in P^\times$.} \hfill $\square$

\subsection{Proposition.} \label{p.EMPROQUA} (Embedding property)
\emph{Let $Q$ be a pointed Pfister quadratic space over $k$ and $P
\subset Q$ a proper pointed Pfister quadratic subspace. Then the
inclusion $P \hookrightarrow Q$ extends to an embedding from
$\dla\mu\dra \otimes P$ to $Q$, for some $\mu \in k^\times$.} \hfill
$\square$

\subsection{$2$-Henselian fields.} \label{ss.HEFI}
Let $F$ be a field of arbitrary characteristic that is endowed
with a normalized discrete valuation $\lambda$, so $\lambda\:F \to
\IZ_\infty := \IZ \cup \{\infty\}$  \index{$\Z_\infty$} is a surjective map satisfying
the following conditions, for all $\alpha,\beta \in F$.
\begin{align*}
\text{$\lambda$ is definite:} \quad &\lambda(\alpha) = \infty
\Longleftrightarrow \alpha = 0. \\
\text{$\lambda$ is sub-additive:} \quad &\lambda(\alpha + \beta)
\geq \mbm\,\{\lambda(\alpha),\lambda(\beta)\}. \\
\text{$\lambda$ is multiplicative:} \quad &\lambda(\alpha\beta) =
\lambda(\alpha) + \lambda(\beta).
\end{align*}

As convenient references for the theory of valuations we mention
\cite{MR2215492,MR2183496}, and particularly \cite{MR1105534} for
the discrete case. We write $\mfo \subseteq F$ for the valuation
ring of $F$ relative to $\lambda$, $\mfp \subseteq \mfo$ for its
valuation ideal and $\bar F := \mfo/\mfp$ for the residue field of
$F$. The natural map from $\mfo$ to $\bar F$ will always be
indicated by $\alpha \mapsto \bar \alpha$. Throughout the
remainder of this paper, we fix a prime element $\pi \in \mfo$.
The quantity
\begin{align}
\label{EEF} e_F := \lambda(2\cdot 1_F),
\end{align}
which is either a non-negative integer or $\infty$, will play an
important role in the sequel. If $\bar F$ has characteristic $2$,
then $e_F > 0$ agrees with what is usually called the
\emph{absolute ramification index}\index{absolute ramification index $e_F$} of $F$.

Due to the quadratic character of the gadgets we are interested in
(composition and conic algebras, pointed quadratic spaces),
requiring $F$ to be Henselian (with respect to $\lambda$) is too
strong a condition. It actually suffices to assume that $F$ be
\emph{$2$-Henselian} in the sense of Dress \cite{MR0206804} or \cite[\S4.2]{MR2183496}, i.e.,
that $F$ satisfies the following two equivalent conditions
\cite[Satz~1]{MR0206804}.
\begin{itemize}
\item [(i)] For any quadratic field extension $K/F$, there is a
unique extension of $\lambda$ to a discrete valuation $\lambda_K$
of $K$ taking values in $\IQ_\infty = \IQ \cup \{\infty\}$.\index{$\Q_\infty$}

\item [(ii)] For all $\alpha_0,\alpha_1,\alpha_2 \in \mfo$ with
$\alpha_0 \in \mfp$, $\alpha_1 \notin \mfp$, the polynomial
\[
\alpha_0\bft^2 + \alpha_1\bft + \alpha_2 \in F[\bft]
\]
is reducible.
\end{itemize}
In this case, the extension $\lambda_K$ of $\lambda$ in (i) is
given by
\begin{align}
\label{EXK} \lambda_K(u) = \frac{1}{2}\lambda\big(N_{K/F}(u)\big)
&& (u \in K).
\end{align}

From now on, $F$ is assumed to be a fixed $2$-Henselian field with
respect to a normalized discrete valuation $\lambda$. For
simplicity, all algebras, quadratic forms etc. over $F$ are
assumed to be finite-dimensional. The characteristic of $F$ is
arbitrary.
%For simplicity, we consider only algebras, quadratic forms etc. that
%are finite-dimensional over $F$.
%\lz

\subsection{Quadratic forms over $2$-Henselian fields.}
\label{ss.QUAHE} Let $q\:V \to F$ be a quadratic form over $F$ and
suppose $q$ is anisotropic. Following Springer \cite{MR0070664}
(in the case of a complete rather than Henselian valuation),
\begin{align}
\label{QUATRI} \lambda\big(q(x + y)\big) \geq
\mbm\,\{\lambda\big(q(x)\big),\lambda\big(q(y)\big)\} &&(x,y \in
V).
\end{align}
For convenience, we give the easy proof of this inequality. It
evidently suffices to show
\begin{align}
\label{POLTRI} \lambda\big(q(x,y)\big) \geq
\mbm\,\{\lambda\big(q(x)\big),\lambda\big(q(y)\big)\}
\end{align}
for all $x,y \in V$. Suppose there are $x,y \in V$ such that
\eqref{POLTRI} does not hold. Then $q(x,y) \neq 0$, and the
polynomial
\[
\frac{1}{q(x,y)}q(\bft x + y) = \frac{q(x)}{q(x,y)}\bft^2 + \bft +
\frac{q(y)}{q(x,y)} \in F[\bft]
\]
has no zero in $F$ since $q$ is anisotropic, but is reducible by
\ref{ss.HEFI}~(ii), a contradiction. \hfill $\square$

Relation \eqref{POLTRI} may be strengthened to
\begin{align}
\label{SPOLTRI} 2\lambda\big(q(x,y)\big) \geq \lambda\big(q(x)\big)
+ \lambda\big(q(y)\big) &&(x,y \in V),
\end{align}
either by appealing to a general result of Bruhat-Tits
\cite[Thm.~10.1.15]{MR0327923}, or by using an ad-hoc argument from
the valuation theory of Jordan rings \cite{MR50:422} in disguise:
for $x,y \in V$, $x \neq 0$, we obtain
\[
Q_xy := q(x,y)x - q(x)y = -q(x)\tau_x(y),
\]
where $\tau_x\:V \to V$ is the reflection in the hyperplane
perpendicular to $x$. In particular, $\tau_x$ leaves $q$
invariant, which implies $\lambda(q(Q_xy)) = 2\lambda(q(x)) +
\lambda(q(y))$, and since $q(x,y)x = Q_xy + q(x)y$, we obtain
\eqref{SPOLTRI}.

%\subsection{The general set-up.} \label{ss.GENSET} One could argue
%with some justification that the appropriate framework for some of
%our subsequent results is provided by \emph{pointed quadratic
%forms} over $F$, i.e., by pairs $(q,e)$ consisting of a quadratic
%form $q\:V \to F$ and a \emph{base point} $e \in V$ for $q$ in the
%sense that $q(e) = 1$; it would also be necessary to require that
%$q$ be round and anisotropic. Instead, we prefer to work with
%non-degenerate conic algebras over $F$ whose norms are round and
%anisotropic. There are two reasons for this. The first reason is
%one of convenience: in this way, the terminology we are going to
%introduce fits more naturally into the context of valuation theory
%and, in particular, doesn't have to be modified when passing later
%on to the setting of ($\lambda$-)composition algebras, where
%pointed quadratic forms over $F$ that are round and anisotropic
%will definitely not be enough. The second reason is that we do not
%lose much in generality by our approach since, over any field,
%\emph{every pointed quadratic form is the norm of some conic
%algebra} \cite{Rosemeier02,Loos11}.

\subsection{Valuation data for pointed quadratic spaces.}\index{pointed quadratic space!valuation data}
\label{ss.VADAL} For the rest of this section, we fix a pointed
quadratic space $Q$ over $F$ which is round and anisotropic. $\ssp$ \\
(a) The map
\begin{align}
\label{EXC} \lambda_Q\:V_Q \longrightarrow \IQ_\infty, \quad x
\longmapsto \lambda_Q(x) := \frac{1}{2}\lambda\big(n_Q(x)\big),
\end{align}
is a norm of $V_Q$ as an $F$-vector space in the sense of
Bruhat-Tits \cite[1.1]{MR788969}, that is, the following relations
hold for all $\alpha \in F$, $x,y \in V_Q$.
\begin{align}
\label{LADEF} \text{$\lambda_Q$ is definite:} \quad &\lambda_Q(x)
= \infty \Longleftrightarrow x = 0, \\
\label{LASU} \text{$\lambda_Q$ is sub-additive:} \quad
&\lambda_Q(x + y) \geq \mbm\,\{\lambda_Q(x),\lambda_Q(y)\}, \\
\label{LASC} \text{$\lambda_Q$ is scalar-compatible:} \quad
&\lambda_Q(\alpha x) = \lambda(\alpha) + \lambda_Q(x),
\end{align}
where \eqref{LASU} is a consequence of
(\ref{ss.QUAHE}.\ref{QUATRI}). Moreover,
\begin{align}
\label{LAST} \lambda_Q(x^\ast) =\,\,&\lambda_Q(x), \\
\label{LANP} \lambda\big(n_Q(x,y)\big) \geq\,\,&\lambda_Q(x) +
\lambda_Q(y), \\
\label{LAT} \lambda\big(t_Q(x)\big) \geq\,\,&\lambda_Q(x),
\end{align}
for all $x,y \in V_Q$, where \eqref{LAST} follows from conjugation
invariance of $n_Q$, \eqref{LANP} from
(\ref{ss.QUAHE}.\ref{SPOLTRI}), and \eqref{LAT} from \eqref{LANP}
for $y = 1_Q$. $\ssp$ \\
(b) $n_Q$ being round, $\Gamma_Q := \lambda_Q(V_Q^\times)$ is an
additive subgroup of $\IQ$ for which \eqref{EXC} implies
\begin{align}
\label{RAIN} \IZ \subseteq \Gamma_Q \subseteq \frac{1}{2}\IZ,
\quad e_{Q/F} := [\Gamma_Q:\IZ] \in \{1,2\}, \quad \Gamma_Q =
\frac{1}{e_{Q/F}}\IZ.
\end{align}
We call $e_{Q/F}$ the \emph{ramification index}\index{pointed quadratic space!ramification index}
of $Q$. $\ssp$ \\
(c) We put
\begin{align}
\label{OCE} \mfo_Q :=\,\,&\{x \in V_Q \mid \lambda_Q(x) \geq 0\}, \\
\label{PECE} \mfp_Q := \{x \in V_Q \mid \lambda_Q(x) > 0\}
=\,\,&\{x \in V_Q \mid \lambda_Q(x) \geq \frac{1}{e_{Q/F}}\}
\subseteq \mfo_Q,
\end{align}
which are both full $\mfo$-lattices in $V_Q$, and
\begin{align}
\label{OCET} \mfo_Q^\times := \mfo_Q \setminus \mfp_Q = \{x \in
V_Q \mid \lambda_Q(x) = 0\},
\end{align}
which is just a subset of $\mfo_Q$ containing $1_Q$. By abuse of
language, the elements of $\mfo_Q^\times$ are called
\emph{units}\index{pointed quadratic space!unit} of $\mfo_Q$. An
element $\Pi \in V_Q^\times$ such that $\lambda_Q(\Pi) > 0$
generates the infinite cyclic group $\Gamma_Q$ belongs to $\mfp_Q$
and is called a \emph{prime element}\index{pointed quadratic
space!prime element} of $\mfo_Q$. Writing $x \mapsto \bar x$ for
the natural map from $\mfo_Q$ to $\mfo_Q/\mfp_Q$ and setting
\[
\bar Q := (V_{\bar Q},n_{\bar Q},1_{\bar Q}), \quad V_{\bar Q} :=
\mfo_Q/\mfp_Q, \quad 1_{\bar Q} := \overline{1_Q},
\]
where $n_{\bar Q}\:\bar Q \longrightarrow \bar F, \quad \bar x
\longmapsto n_{\bar Q}(\bar x) := \overline{n_Q(x)}$ is the first
residue form of $n_Q$, we obtain a pointed quadratic space over
$\bar F$, called the \emph{pointed quadratic residue
space}\index{pointed quadratic space!pointed quadratic residue
space} of $Q$, which is round and anisotropic. Here only roundness
of $n_{\bar Q}$ demands a proof, so let $u \in \mfo_Q^\times$.
Since $n_Q$ is round, there exists a linear bijection $f\:V_Q \to
V_Q$ with $n_Q(f(x)) = n_Q(u)n_Q(x)$ for all $x \in V_Q$. Hence
$f$ stabilizes $\mfo_Q$ as well as $\mfp_Q$ and thus canonically
induces a similarity transformation $V_{\bar Q} \to V_{\bar Q}$
relative to $n_{\bar Q}$ with multiplier $n_{\bar Q}(\bar u)$.

We call
\begin{align}
\label{RESDEG} f_{Q/F} := \mbd_{\bar F}(V_{\bar Q})
\end{align}
the \emph{residue degree}\index{pointed quadratic space!residue
degree} of $Q$. For convenience, we collect some of the above
properties of $\bar Q$ in the following proposition, which is
really stating the obvious.

\begin{prop} \label{p.PORESQUA} {Norm, trace and
conjugation of the pointed quadratic residue space $\bar Q$ of $Q$
are given by the formulas}
\begin{align}
\label{PORESNO} n_{\bar Q}(\bar x) =\,\,&\overline{n_Q(x)}, \\
\label{PORESTR} t_{\bar Q}(\bar x) =\,\,&\overline{t_Q(x)}, \\
\label{PORESCO} \bar x^\ast =\,\,&\overline{x^\ast}
\end{align}
{for all $x \in \mfo_Q$. Moreover, $n_{\bar Q}$ is round and
anisotropic.} \hfill $\square$
\end{prop}

Here is the easiest example that is not totally trivial.

\begin{eg} \label{PORES.eg}
Let $Q = (V, q, e)$ be such that $V$ has basis $e, j$ and $q(ue +
vj) = u^2 - \mu v^2$ for some $\mu \in \mfo^\times$ and all $u, v
\in F$.  (This is the case $P = \qform{1}$ of \ref{ss.ENLAR}; $q
\cong \pform{\mu}$.)  If $\bmu$ is not a square in $\Fb$, then
$\mfo_Q$ and $\mfp_Q$ are the $\mfo$- and $\mfp$-spans of $e, j$
respectively and the pointed quadratic residue space of $Q$ is
naturally identified with the quadratic extension $\Fb(\sqrt{\bmu})$
with quadratic form the squaring map $z \mapsto z^2$ and base point
$1 \in \Fb$.
\end{eg}

\subsection{Lemma.} \label{l.COVA} \emph{The map}
\[
\IZ \times \left\{y \in V_Q^\times \left| \ 0 \leq \lambda_Q(y) \leq 1 -
\frac{1}{e_{Q/F}} \right. \right\} \longrightarrow V_Q^\times, \quad (m,y)
\longmapsto \pi^my,
\]
\emph{is surjective.}

\Proof For $x \in V_Q^\times$, some $n \in \IZ$ has $\lambda_Q(x)
= \frac{n}{e_{Q/F}}$ by (\ref{ss.VADAL}.\ref{RAIN}), and writing
$n = me_{Q/F} + r$, $m,r \in \IZ$, $0 \leq r \leq e_{Q/F} - 1$,
the element $y := \pi^{-m}x \in V_Q^\times$ satisfies $0 \leq
\lambda_Q(y) \leq 1 - \frac{1}{e_{Q/F}}$.
\begin{flushright}
\hfill $\square$
\end{flushright}

\subsection{Proposition.} \label{p.EFEN} $e_{Q/F}f_{Q/F} = \mbd_F(V_Q)$.

\Proof If $e_{Q/F} = 1$, then $\mfp_Q = \pi \mfo_Q$, and the
dimension of $V_{\bar Q} = \mfo_Q \otimes_\mfo \bar F$ over $\bar
F$ agrees with the dimension of $V_Q$ over $F$. We may therefore
assume $e_{Q/F} = 2$. Let $\Pi$ be a prime element of $\mfo_Q$ and
$f\:V_Q \to V_Q$ a norm similarity with multiplier $n_Q(\Pi)$.
Then $\lambda_Q(f(x)) = \frac{1}{2} + \lambda_Q(x)$ for all $x \in
V_Q^\times$, which implies $f(\mfo_Q) = \mfp_Q$, $f(\mfp_Q) = \pi
\mfo_Q$, and $f$ induces canonically an $\bar F$-linear bijection
from $V_{\bar Q}$ onto $\mfp_Q/\pi \mfo_Q$. Combined with the
filtration $\mfo_Q \otimes_\mfo \bar F = \mfo_Q/\pi\mfo_Q \supset
\mfp_Q/\pi\mfo_Q \supset \{0\}$ and the isomorphism
$(\mfo_Q/\pi\mfo_Q)/(\mfp_Q/\pi\mfo_Q) \cong V_{\bar Q}$ as vector
spaces over $\bar F$, this implies $\mbd_F(V_Q) =2\mbd_{\bar
F}(V_{\bar Q)}$, as desired. \hfill $\square$

\Remark Prop.~\ref{p.EFEN} becomes \emph{false} if $Q$ is not
assumed to be round, for example if $e_{Q/F} = 2$ and $\mbd_F(V_Q)$
is odd, see \cite[Satz~6.3]{MR50:422} for generalization.

\subsection{Connecting with conic algebras.} \label{ss.CONCON} Let
$C$ be a finite-dimensional conic algebra over $F$ and suppose its
norm is round and anisotropic; for example, $C$ could itself be a
composition algebra, or it could arise from a composition algebra
by means of the Cayley-Dickson process. Applying as we may the
preceding considerations to $Q_C$, the pointed quadratic space
associated with $C$ via \ref{ss.POQUA}, we systematically adhere
to the following convention: all notation and terminology
developed up to now and later on for pointed quadratic spaces will
be applied without further comment to $C$ in place of $Q_C$,
modifying subscripts accordingly whenever possible. For example,
\begin{align}
\label{EXCO} \lambda_C\:C \longrightarrow \IQ_\infty, \quad x
\longmapsto \lambda_C(x) := \frac{1}{2}\lambda\big(n_C(x)\big),
\end{align}
is a norm of $C$ as an $F$-vector space, $e_{C/F} := e_{Q_C/F}$ is
the \emph{ramification index}\index{conic algebra!ramification
index}, $f_{C/F} := f_{Q_C/F}$ is the \emph{residue
degree}\index{conic algebra!residue degree} and $\bar C:=
\overline{Q_C}$ is the \emph{pointed quadratic residue
space}\index{conic algebra!pointed quadratic residue space} of $C$.
But note that $\bar C$ in general is \emph{not} a conic algebra over
$\bar F$ in a natural way unless $C$ is a composition algebra; see
\ref{ss.TAWI} and Section~\ref{s.LANOVA} below for further
discussion.

\subsection{Tame and wild pointed quadratic spaces.} \label{ss.TAWI} If
$t_{\bar Q}$ is non-zero, then $Q$ is said to be
\emph{tame}.\index{pointed quadratic space!tame}\index{tame pointed quadratic space} %
Otherwise, i.e.,
if $t_{\bar Q} = 0$, then $Q$ is said to be
\emph{wild}.\index{pointed quadratic space!wild}\index{wild pointed quadratic space} %
For $Q$ to be wild it is clearly necessary that $\bar F$ have
characteristic $2$. Applying \cite[Prop.~1]{MR51:635}, and bearing
in mind the conventions of \ref{ss.CONCON}, a composition division
algebra $C$ over $F$ is tame (resp. wild) iff $\bar C$ is a
composition algebra (resp. a purely inseparable field extension of
exponent at most $1$) over $\bar F$. Extending the terminology of
\cite{MR51:635} to the present more general set-up, we call $Q$
\emph{unramified} (resp. \emph{ramified}) if $Q$ is tame with
$e_{Q/F} = 1$ (resp. $e_{Q/F} = 2$). For $Q$ to be unramified it is
necessary and sufficient that the quadratic form $n_Q$ have good
reduction with respect to $\lambda$ in the sense of Knebusch
\cite{MR0412101,{Kn07}}. This will follow from
Prop.~\ref{p.PROTE}~(c) below.

The preceding definitions seem to assign a distinguished role to
the base point of a pointed quadratic space. But this is not so as
will be seen in Prop.~\ref{p.PROTE} below.

%%%%%%%%%%%%%%%%%%%%%%%%%%%%%%%%%%%%%%%%%%%%%%%%%%
\section{Trace and norm exponent.} \label{s.VADA}
This section serves a double purpose. Working with a fixed
non-singular, round and anisotropic pointed quadratic space $Q$
over a $2$-Henselian field $F$, we attach wildness-detecting
invariants to $Q$. Moreover, we present a device measuring how far
a scalar in $\mfo^\times$ is removed from being the norm of an
appropriate element in $V_Q$.

\subsection{Trace ideal and trace exponent.} \label{ss.TITE} Since
$\mfo_Q \subseteq V_Q$ is a full $\mfo$-lattice, its image under
the trace of $Q$ by non-singularity is a non-zero ideal in $\mfo$,
called the \emph{trace ideal}\index{pointed quadratic space!trace
ideal} of $Q$. But $\mfo$ is a discrete valuation ring, so there
is a unique integer $\tex(Q) \geq 0$ such that
\begin{align}
\label{TIP} t_Q(\mfo_Q) = \mfp^{\tex(Q)}.
\end{align}
We call $\tex(Q)$ the \emph{trace exponent}\index{pointed
quadratic space!trace exponent}\index{trace exponent ($\texp$)} of $Q$ and have
\begin{align}
\label{TEXMIN} \tex(Q) = \mbm\,\{\lambda\big(t_Q(x)\big) \mid x
\in \mfo_Q\}.
\end{align}
The image of $\mfo_Q \otimes_\mfo \mfo_Q$ under the linear map $x
\otimes y \mapsto n_Q(x,y)$ is an ideal in $\mfo$ denoted by
$\partial n_Q(\mfo_Q \otimes_\mfo \mfo_Q)$. The following result
relates the ideals just defined to one another but also to wild
and tame pointed quadratic spaces.

\subsection{Proposition.} \label{p.PROTE} (a) \emph{$Q$ is wild
if and only if $\tex(Q) > 0$.} $\ssp$ \\
(b) \emph{If $P \subseteq Q$ is pointed quadratic subspace that is
round and non-singular, then $\tex(P) \geq \tex(Q)$.} $\ssp$ \\
(c) \emph{$t_Q(\mfo_Q) = \partial n_Q(\mfo_Q \otimes_\mfo
\mfo_Q)$.} $\ssp$ \\
(d) \emph{$Q$ is tame if and only if $\bar Q$ is non-singular.}

\Proof (a) and (b) are obvious. Before proving (c),(d), we let $x
\in Q^\times$ and, by using a Jordan isotopy argument in disguise,
pass to $Q^x : = (V_Q,n_Q(x)^{-1}n_Q,x)$, which is a non-singular
quadratic space over $F$. Since $n_Q$ is round, the norms of $Q$
and $Q^x$ are isometric, forcing $Q$ and $Q^x$ to be isomorphic as
pointed quadratic spaces by Prop.~\ref{p.PONO}. This implies
\begin{align}
\label{ISEX} t_Q(\mfo_Q) = t_{Q^x}(\mfo_{Q^x}) =\,\,&
\big\{n_Q\big(x,n_Q(x)^{-1}y\big) \mid y \in V_Q,\,\,\lambda_Q(y)
\geq \lambda_Q(x)\big\} \\
=\,\,&\{n_Q(x,y) \mid y \in V_Q,\,\,\lambda_Q(y) \geq
-\lambda_Q(x)\}. \notag
\end{align}
(c) The left-hand side is clearly contained in the right, so it
suffices to show $n_Q(x,y) \in t_Q(\mfo_Q)$ for all $x,y \in
\mfo_Q$, $x \neq 0$. But this follows from \eqref{ISEX} and
$\lambda_Q(x) \geq 0$. \\
(d) Non-singularity of $\bar Q$ is clearly sufficient for $Q$ to
be tame. Conversely, suppose $Q$ is tame and let $x \in
\mfo_Q^\times$. Then \eqref{ISEX} produces an element $y \in
\mfo_Q$ with $n_Q(x,y) = 1 \in t_Q(\mfo_Q)$, hence $n_{\bar
Q}(\bar x,\bar y) \neq 0$. \hfill $\square$

\subsection{Trace generators.} \label{ss.TRAG} An element $w_0 \in
\mfo_Q$ where the minimum in (\ref{ss.TITE}.\ref{TEXMIN}) is
attained, i.e., with $\lambda(t_Q(w_0)) = \tex(Q)$ is called a
\emph{trace generator}\index{trace generator} of $Q$. If even
$t_Q(w_0) = \pi^{\tex(Q)}$, we speak of a \emph{normalized} trace
generator\index{trace generator!normalized}, dependence on $\pi$
being understood. Trace generators always exist (as do normalized
 ones) and their traces generate the trace ideal of $C$.
Moreover, they satisfy the inequalities
\begin{align}
\label{ESTRAG} 0 \leq \lambda_Q(w_0) \leq 1 - \frac{1}{e_{Q/F}}.
\end{align}
Indeed, assuming $\lambda_Q(w_0) > 1 - \frac{1}{e_{Q/F}}$ implies
\[
\lambda_Q(\pi^{-1}w_0) = \lambda_Q(w_0) - 1 > -\frac{1}{e_{Q/F}},
\]
and since $\lambda_Q(\pi^{-1}w_0)$ belongs to
$\frac{1}{e_{Q/F}}\IZ$, we conclude $\pi^{-1}w_0 \in \mfo_Q $. But
now $t_Q(\pi^{-1}w_0) = \pi^{\tex(Q)-1} \notin \mfp^{\tex(Q)}$
leads to a contradiction.

\subsection{Tignol's invariant $\wid(C)$.}\index{Tignol's invariant $\wid$}
\label{ss.WID} Another invariant that fits into the present set-up
is due to Tignol \cite[pp.~9,17]{Tignol:wild} for separable field
extensions and central associative division algebras of degree $p$
over Henselian fields (not necessarily discrete) having residual
characteristic $p > 0$. As a straightforward adaptation of
Tignol's definition to the situation we are interested in, we
define %the \emph{width} of $C$ by
\begin{align}
\label{DEWI} \wid(Q) := \mbm\,\{\lambda\big(t_Q(x)\big) -
\lambda_Q(x) \mid x \in V_Q^\times \};
\end{align}
thanks to (\ref{ss.VADAL}.\ref{LAT}), it is a non-negative rational
number. Moreover, it is closely related to the trace exponent, as
the following proposition shows.

\subsection{Proposition.} \label{p.WITE} (a) \emph{$\wid(Q) =
\tex(Q)$ or $\wid(Q) = \tex(Q) - \frac{1}{2}$.} $\ssp$ \\
(b) \emph{There exists a trace generator $w_0$ of $Q$ with}
\begin{align}
\label{WIDTEX} \wid(Q) = \tex(Q) - \lambda_Q(w_0).
\end{align}
(c) \emph{If $e_{Q/F}\ = 1$, then $\wid(Q) = \tex(Q)$.}

\Proof Since the map
\[
\varphi\:V_Q^\times \longrightarrow \Gamma_Q =
\frac{1}{e_{Q/F}}\IZ, \quad x \longmapsto \varphi(x) :=
\lambda\big(t_Q(x)\big) - \lambda_Q(x),
\]
is homogeneous of degree zero, so $\varphi(\alpha x) = \varphi(x)$
for all $\alpha \in F^\times$, $x \in V_Q^\times$,
Lemma~\ref{l.COVA} implies $\wid(Q) = \mbm\,\{\varphi(x) \mid x
\in S\}$, where
\[
S := \{x \in V_Q^\times \mid 0 \leq \lambda_Q(x) \leq 1 -
\frac{1}{e_{Q/F}}\}.
\]
Accordingly, let $w_0 \in S$ satisfy
\begin{align}
\label{PHIWE} \varphi(w_0) = \wid(Q).
\end{align}
Given any trace generator $w_0^\prime$ of $Q$, the chain of
inequalities
\begin{align}
\label{TEXWID} 0 \leq \lambda_Q(w_0^\prime)
=\,\,&\lambda\big(t_Q(w_0^\prime)\big) - \varphi(w_0^\prime) \leq
\tex(Q) - \wid(Q) \leq \lambda\big(t_Q(w_0)\big) - \varphi(w_0) \\
=\,\,&\lambda_Q(w_0) \leq 1 - \frac{1}{e_{Q/F}} \notag
\end{align}
implies (a) and (c), while in (b) we may assume $e_{Q/F} = 2$ since
otherwise \eqref{WIDTEX} holds for $w_0^\prime$, i.e., for any trace
generator of $Q$. Since all quantities in \eqref{TEXWID} belong to
$\frac{1}{2}\IZ$, we either have $\lambda_Q(w_0^\prime) = \tex(Q) -
\wid(Q)$ or
\[
\tex(Q) - \wid(Q) = \lambda\big(t_Q(w_0)\big) - \varphi(w_0) =
\lambda_Q(w_0) = \frac{1}{2}.
\]
In the latter case, \eqref{PHIWE} shows that $w_0$ is a trace
generator of $Q$, and the proof of (b) is complete. \hfill
$\square$

\subsection{Regular trace generators.} \label{ss.RETREGEN} Trace
generators of $Q$ satisfying (\ref{p.WITE}.\ref{WIDTEX}) are
called \emph{regular}\index{trace generator!regular}. If $Q$ has
ramification index $1$, then every trace generator by
(\ref{ss.TRAG}.\ref{ESTRAG}) and Prop.~\ref{p.WITE}~(c) is regular
but, in general, this need not be so, cf. Cor.~\ref{c.TRGEUN} and
Thm.~\ref{t.WITY} below.

\subsection{Example} \label{e.BAFI} If $F$ has characteristic not
$2$, it is a composition division algebra over itself, with norm
and trace given by $n_F(\alpha) = \alpha^2$, $t_F(\alpha) =
2\alpha$ for all $\alpha \in F$. Hence (\ref{ss.HEFI}.\ref{EEF}),
(\ref{ss.TITE}.\ref{TIP}) and (\ref{ss.WID}.\ref{DEWI}) imply
\begin{align}
\label{TEEF} \tex(F) = \wid(F) = e_F,
\end{align}
and $w_0 = \frac{\pi^{e_F}}{2} \in \mfo^\times$ is the unique
normalized trace generator of $F$; it is obviously regular. From
\eqref{TEEF} and Props.~\ref{p.PROTE}~(c),~\ref{p.WITE} we now
conclude
\begin{align}
\label{ETE} 0 \leq \wid(Q) \leq \tex(Q) \leq e_F,
\end{align}
which for trivial reasons also holds in characteristic $2$.

\subsection{The norm exponent.} \label{ss.NE} We wish to measure
how far a given unit in the valuation ring of $F$ is removed from
being the norm of an element in $V_Q$ or, equivalently, in
$\mfo_Q^\times$. To this end, we observe that, given $\alpha \in
\mfo^\times$ and an integer $d \geq 0$, the following conditions
are equivalent.
\begin{itemize}
\item [(i)] $\alpha$ is a norm of $\mfo_Q^\times \bmod \mfp^d$,
i.e., there exists an element $v \in \mfo_Q^\times$ with
\[
\alpha - n_Q(v) \in \mfp^d.
\]

\item[(ii)] There exist elements $\beta \in \mfo$, $v \in
\mfo_Q^\times$ with
\[
\alpha = (1 - \pi^d\beta)n_Q(v).
\]
\end{itemize}
Thus the set $N_Q(\alpha)$ of non-negative integers $d$ satisfying
(i)/(ii) above contains $0$, and $d \in N_Q(\alpha)$ implies
$d^\prime \in N_Q(\alpha)$ for all integers $d^\prime$, $0 \leq
d^\prime \leq d$. We therefore put
\begin{align}
\label{SUN} \nex_Q(\alpha) := \mathrm{sup}\,N_Q(\alpha)
\end{align}
and call this the \emph{norm exponent}\index{norm exponent (nexp)}
of $\alpha$ relative to $Q$; it is either a non-negative integer or
$\infty$. Roughly speaking, the bigger the norm exponent becomes,
the closer $\alpha$ gets to being a norm of $Q$. \emph{More
precisely, $2^{-\nex_Q(\alpha)}$ is the (minimum) distance of
$\alpha$ from the subset $n_Q(\mfo_Q^\times) \subseteq F^\times$
relative to the metric induced by the absolute value $\xi \mapsto
\vert\xi\vert := 2^{-\lambda(\xi)}$.}

A number of useful elementary properties of the norm exponent are
collected in the following proposition, whose straightforward
proof is left to the reader.

\subsection{Proposition.} \label{p.PRONE} \emph{With the notations
and assumptions of} \ref{ss.NE}, \emph{let $\alpha, \alpha^\prime
\in \mfo^\times$.} $\ssp$ \\
(a) \emph{If $\alpha \in n_Q(\mfo_Q^\times)$, then $\nex_Q(\alpha)
= \infty$.} $\ssp$ \\
(b) \emph{If $d = \nex_Q(\alpha)$ is finite, then $\alpha$ can be
written in the form}
\[
\alpha = (1 - \pi^d\beta)n_Q(v), \quad \beta \in \mfo, \quad v \in
\mfo_Q^\times,
\]
\emph{and every such representation of $\alpha$ satisfies $\beta
\in \mfo^\times$.} $\ssp$ \\
(c) \emph{$\nex_Q(\alpha) = 0$ if and only if $\bar \alpha \notin
n_{\bar Q}(V_{\bar Q}^\times)$.} $\ssp$ \\
(d) \emph{$\alpha \equiv \alpha^\prime \bmod n_Q(V_Q^\times)$
implies
$\nex_Q(\alpha) = \nex_Q(\alpha^\prime)$.} $\ssp$ \\
(e) \emph{$\nex_P(\alpha) \leq \nex_Q(\alpha)$ for any round
non-singular pointed quadratic subspace $P \subseteq Q$.}
\begin{flushright}
\hfill $\square$
\end{flushright} \lz

\noindent By Prop.~\ref{p.PRONE}~(a), the elements of
$n_Q(\mfo_Q^\times)$ have infinite norm exponent. While the
converse is also true, we can do better than that by showing that
the norm exponents of elements in $\mfo^\times \setminus
n_Q(\mfo_Q^\times)$ are uniformly bounded from above. Indeed, we
have the following result.

\subsection{Local Norm Theorem.}\index{Local norm theorem} \label{t.HENOTH} \emph{Let
$\alpha \in \mfo^\times \setminus n_Q(\mfo_Q^\times) = \mfo^\times
\setminus n_Q(V_Q^\times)$. Then}
\[
\nex_Q(\alpha) \leq 2\wid(Q).
\]
\emph{More precisely, letting $w_0 \in \mfo_Q$ be a normalized
regular trace generator of $Q$, every $\beta \in \mfo$ admits a
$\gamma \in \mfo$ with}
\begin{align}
\label{NOEQ} 1 - \pi^{2\wid(Q)+1}\beta = n_Q\big(1_Q +
\pi^{\tex(Q)}\gamma w_0\big).
\end{align}

\begin{proof}  Arguing indirectly, and using roundness of $n_Q$, the first
part of the theorem follows from the second. To establish the second
part, we apply Prop.~\ref{p.WITE}~(a),(b) and obtain
\[
0 \leq d := \lambda\big(n_Q(w_0)\big) = 2\lambda_Q(w_0) =
2\big(\tex(Q) - \wid(Q)\big) \leq 1.
\]
Since $F$ is $2$-Henselian, the polynomial
\[
g := \pi^{1-d}n_Q(w_0)\bft^2 + \bft + \beta \in \mfo[\bft]
\subseteq F[\bft]
\]
by \ref{ss.HEFI}~(ii) is reducible, and the two roots
$\delta,\delta^\prime \in F$ of $\pi^{d-1}n_Q(w_0)^{-1}g$ satisfy
the relation $\delta\delta^\prime = \pi^{d-1}n_Q(w_0)^{-1}\beta
\in \mfp^{-1}$. Thus we may assume $\delta \in \mfo$ and have
$\pi^{1-d}n_Q(w_0)\delta^2 + \delta + \beta = 0$. Setting $\gamma
:= \pi^{1-d}\delta \in \mfo$, we therefore obtain
\begin{align*}
n_Q\big(1 + \pi^{\tex(Q)}\gamma w_0\big) =\,\,&1 +
\pi^{\tex(Q)+1-d}\delta t_Q(w_0) +
\pi^{2\tex(Q)}\pi^{2(1-d)}\delta^2n_Q(w_0) \\
=\,\,&1 + \pi^{2\tex(Q)-d+1}\big(\delta +
\pi^{1-d}n_Q(w_0)\delta^2\big) \\
=\,\,&1 - \pi^{2\wid(Q)+1}\beta.\qedhere
\end{align*}
\end{proof}

\noindent For a version of this result addressed to central
associative division algebras of degree $p = \mbc(\bar F) > 0$, see
Kato \cite[Prop.~2~(iii)]{Kato:gen1}.

\subsection{Corollary.} \label{c.NOCL} \emph{Let $\mu,\mu^\prime
\in \mfo^\times$. Then}
\[
\mu \equiv \mu^\prime \bmod\,n_Q(Q^\times) \Longleftrightarrow
\exists\,x \in \mfo_Q^\times\,:\,\mu \equiv \mu^\prime n_Q(x)
\bmod\,\mfp^{2\wid(Q)+1}.
\]
\begin{flushright}
\vspace{-8pt} $\square$
\end{flushright}

\subsection{The connection with the quadratic defect.}\index{quadratic defect}
\label{ss.QUADE } We assume $\mbc(F) \neq 2$ and return to the
composition division algebra $C = F$ of \ref{e.BAFI}. Comparing
\ref{ss.NE} with \cite[\S63A]{MR0152507}, we see that
$\mfp^{\nex_F(\alpha)}$ is the quadratic defect\index{quadratic
defect} of $\alpha \in \mfo^\times$. Moreover, $\tex(F) = \wid(F)
= e_F$, and $w_0 := \frac{\pi^{e_F}}{2} \in \mfo^\times$ is the
unique  normalized regular trace generator of $F$. In particular,
$\frac{4}{\pi^{2e_F}} \in \mfo^\times$, and the change of
variables $\beta = -\frac{4}{\pi^{2e_F}}\beta^\prime$, $\gamma =
\frac{4}{\pi^{2e_F}}\gamma^\prime$ converts
(\ref{t.HENOTH}.\ref{NOEQ}) into the relation
\[
1 + 4\pi\beta^\prime = (1 + 2\pi\gamma^\prime)^2.
\]
Hence the local norm theorem becomes the local square theorem of
\cite[63:1]{MR0152507} or \cite[VI.2.19]{Lam} in the special case
$Q = Q_F$. See also \cite[Prop.~4.1.2]{CT:kato} for an extension
of this result to residual characteristics other than $2$. \lz

\noindent The local norm theorem has numerous applications. One of
these can already be given in this section; a useful technical
lemma prepares the way.

\subsection{Lemma.} \label{l.NOCO} \emph{Suppose $Q$ is wild and
has ramification index $e_{Q/F} = 1$. For $d \in \IZ$, $0 \leq d
\leq 2\tex(Q) = 2\wid(Q)$ and $\beta \in \mfo^\times$, the
following conditions are equivalent.}
\begin{itemize}
\item [(i)] $1 - \pi^d\beta \in n_Q(\mfo_Q^\times)$.

\item [(ii)] \emph{$1 - \pi^d\beta \in \mfo^\times$ and there are
elements $m \in \IZ$, $w \in \mfo_Q^\times$ with}
\begin{align}
\label{BEWE} d = 2m, \quad \beta = -n_Q(w) + \pi^{-m}t_Q(w).
\end{align}
\end{itemize}

\Proof (i) $\Longrightarrow$ (ii). There exists an element $v \in
\mfo_Q^\times$ with $1 - \pi^d\beta = n_Q(v) \in \mfo^\times$. For
$d = 0$ we put $w = 1_Q - v \in \mfo_Q$ and obtain $1 - \beta =
n_Q(1_Q - w) = 1 - t_Q(w) + n_Q(w)$, hence $\beta = -n_Q(w) +
t_Q(w)$. But $t_Q(w) \in \mfp^{\tex(Q)} \subseteq \mfp$ since $Q$
is wild, forcing $w \in \mfo_Q^\times$. Thus (ii) holds with $m =
0$. We may therefore assume $d > 0$. Then $1 = \overline{n_Q(v)} =
n_{\bar Q}(\bar v)$, forcing $\bar v = 1$ since $n_{\bar Q}$ is
anisotropic and $n_{\bar Q}(1_{\bar Q} - \bar v) = 0$ by wildness
of $Q$. Combining with $e_{Q/F} = 1$, we find an integer $m > 0$
and a unit $w \in \mfo_Q^\times$ with $v = 1_Q - \pi^mw$.
Expanding the right-hand side of $1 - \pi^d\beta = n_Q(1_Q - \pi^m
w)$, we conclude
\begin{align}
\label{DEBET} \pi^d\beta = -\pi^{2m}n_Q(w) + \pi^mt_Q(w),
\end{align}
which in turn yields the estimate
\begin{align*}
2\,\tex(Q) \geq\,\,&d = \lambda(\pi^d\beta) =
\lambda\big(-\pi^{2m}n_Q(w) + \pi^mt_Q(w)\big) \geq \mbm\,\{2m,m +
\tex(Q)\}.
\end{align*}
This implies
\begin{align}
\label{EMTE} m \leq \tex(Q), \quad d \geq 2m,
\end{align}
and \eqref{DEBET} attains the form
\begin{align}
\label{DEBEG} \pi^d\beta = -\pi^{2m}n_Q(w)(1 - \gamma), \quad
\gamma := \pi^{-m}n_Q(w)^{-1}t_Q(w) \in \mfp^{\tex(Q)-m}.
\end{align}
If $m < \tex(Q)$, then \eqref{BEWE} follows from
\eqref{DEBET},\eqref{DEBEG}. On the other hand, if $m = \tex(Q)$,
then \eqref{EMTE} implies $d = 2\tex(Q) = 2m$, and \eqref{BEWE}
follows from \eqref{DEBET}.

(ii) $\Longrightarrow$ (i). Setting $v := 1_Q - \pi^mw \in \mfo_Q$
and applying \eqref{BEWE}, we conclude $n_Q(v) = 1 - \pi^d\beta
\in \mfo^\times$, hence $v \in \mfo_Q^\times$, and (i) holds.
\hfill $\square$

\Remark For $d > 0$, the condition $1 - \pi^d\beta \in
\mfo^\times$ in (ii) is of course automatic.

\subsection{Proposition.} \label{c.ODNE} \emph{Let $P$ be a
pointed quadratic space over $F$ that is non-singular, round and
anisotropic with $e_{P/F} = 1$, and let $d$ be an odd
integer with $0 \leq d < 2\tex(P)$.}  \\
(a) \emph{If $\beta \in \mfo$, then $\mu := 1 - \pi^d\beta \in
\mfo^\times$ and}
\[
d = \nex_P(\mu) \Longleftrightarrow \beta \in \mfo^\times.
\]
(b) \emph{If $\beta_i \in \mfo^\times$ and $\mu_i := 1 -
\pi^d\beta_i$ for $i = 1,2$, then}
\[
\dla\mu_1\dra \otimes P \cong \dla\mu_2\dra \otimes P
\Longrightarrow \overline{\beta_1} = \overline{\beta_2}.
\]

\Proof (a) $d = \nex_P(\mu)$ implies $\beta \in \mfo^\times$ by
Prop.~\ref{p.PRONE}~(b). Conversely, suppose $\beta \in
\mfo^\times$. Then $\mu \notin n_P(V_P^\times)$ by
Lemma~\ref{l.NOCO}. Thus the local norm theorem \ref{t.HENOTH}
implies $d \leq d^\prime := \nex_P(\mu) \leq 2\tex(P)$, and from
Prop.~\ref{p.PRONE}~(b) we obtain a representation $\mu =
\mu^\prime n_P(v^\prime)$, $\mu^\prime = 1 -
\pi^{d^\prime}\beta^\prime$ for some $\beta^\prime \in
\mfo^\times$, $v^\prime \in \mfo_P^\times$. Assuming $d <
d^\prime$ would imply that $\mu\mu^\prime = 1 - \pi^d\gamma$,
$\gamma := \beta + \pi^{d^\prime - d}\beta^\prime -
\pi^{d^\prime}\beta\beta^\prime \in \mfo^\times$, does not belong
to $n_P(V_P^\times)$, again by Lemma~\ref{l.NOCO}, in
contradiction to $\mu\mu^\prime = n_P(\mu^\prime v^\prime)$.

(b) Arguing indirectly, let us assume $\overline{\beta_1} \neq
\overline{\beta_2}$. Then $\mu_1\mu_2 = 1 - \pi^d\beta$, where
$\beta = \beta_1 + \beta_2 - \pi^d\beta_1\beta_2 \in \mfo$
satisfies $\bar \beta = \overline{\beta_1} - \overline{\beta_2}
\neq 0$, hence $\beta \in \mfo^\times$. By (a), the norm exponent
of $\mu_1\mu_2$ relative to $P$ is $d$, forcing $\mu_1\mu_2 \notin
n_P(V_P^\times)$ (this also follows from Lemma~\ref{l.NOCO}).
Hence $\dla\mu_1\dra \otimes P$ and $\dla\mu_2\dra \otimes P$ are
not isomorphic by Prop.~\ref{p.CDPOQUA}~(b). \hfill $\square$ \lz

\noindent We will see in Example~\ref{e.NEXP}~(b) below that the
converse of Prop.~\ref{c.ODNE}~(b) does not hold.

\subsection{Lemma.} \label{l.TECH} \emph{Let $\beta \in \mfo$, $w
\in \mfo_Q$ and suppose $d,m$ are non-negative integers. Determine
$\alpha \in \mfo$ by $t_Q(w) = \pi^{\tex(Q)}\alpha$. Then}
\begin{align*}
(1 - \pi^d\beta)n_Q(1_Q - \pi^m w) =\,\,&1 - \pi^d\beta +
\pi^{2m}n_Q(w) - \pi^{\tex(Q)+m}\alpha \,+ \\
\,\,&\pi^{\tex(Q)+d+m}\alpha\beta - \pi^{d+2m}\beta n_Q(w).
\end{align*}

\Proof Expand the left-hand side in the obvious way. \hfill
$\square$

%%%%%%%%%%%%%%%%%%%%%%%%%%%%%%%%%%%%%%%%%%%%%%%%%%
\section{Valuation data under enlargements.} \label{s.VADAEN} In
this section we will be concerned with the question of what happens
to the valuation data ramification index (\ref{ss.VADAL}~(b)),
pointed quadratic residue space (\ref{ss.VADAL}~(c)) and trace
exponent (\ref{ss.TITE}) when passing from a pointed quadratic space
$P$ to $\dla\mu\dra \otimes P$, $\mu \in F^\times$. We will answer
this question not in full generality but only under the additional
hypothesis that $P$ have ramification index $1$. This hypothesis
derives its justification from the fact that, if $F$ has
characteristic zero, every anisotropic pointed $(n + 1)$-Pfister
quadratic space over $F$ contains a pointed $n$-Pfister quadratic
subspace of ramification index $1$. We will prove this in
Prop.~\ref{unramsubsym} and Thm.~\ref{texp.depth}\eqref{texp.e} below by using
methods from algebraic $K$-theory. It would be interesting to know
whether the result in question also holds for $F$ having
characteristic $2$.

\subsection{The general set-up.} \label{ss.GEST} (a) We fix a
$2$-Henselian field $F$ and a pointed quadratic space $P$ over $F$
which is non-singular, round and anisotropic. We also assume
throughout that $P$ has ramification index $e_{P/F} = 1$, which
implies
\begin{align}
\label{WITE} \wid(P) = \tex(P)
\end{align}
by Prop.~\ref{p.WITE}~(c) and
\begin{align}
\label{GALA} \Gamma_P = \lambda_P(P^\times) = \IZ, \quad
\lambda\big(n_P(P^\times)\big) = 2\IZ
\end{align}
by \ref{ss.VADAL}~(b) and (\ref{ss.VADAL}.\ref{EXC}). $\ssp$ \\
(b) We are interested in pointed quadratic spaces $Q = \dla\mu\dra
\otimes P$, $\mu \in F^\times$, as in \ref{ss.ENLAR}; in particular,
we recall $V_Q = V_P + V_Pj$ as vector spaces over $F$. By
Prop.~\ref{p.CDPOQUA}~(b), we may and always will assume
$\lambda(\mu) \in \{0,1\}$, so $\mu$ is either a unit or a prime
element in $\mfo$. \lz

\noindent There are  two harmless cases which we treat first. One of
them arises when $\mu$ is a prime, the other when $\mu$ is a unit
and $P$ is tame.

\subsection{Proposition.} \label{p.CDPR} \emph{If $\mu$ is a prime
element in $\mfo$, then $Q:= \dla\mu\dra \otimes P$ is a
non-singular, round and anisotropic pointed quadratic space over
$F$ with}
\begin{align}
\label{UVEPR} \lambda_Q(u + vj) =\,\,&
\mbm\,\{\lambda_P(u),\lambda_P(v) + \frac{1}{2}\} && (u,v \in
V_P),
\\
\label{OCEPR} \mfo_Q = \mfo_P \oplus \mfo_Pj&, \quad \mfp_Q =
\mfp_P \oplus \mfo_Pj, \\
\label{VADAPR} e_{Q/F} = 2, \quad \bar Q =\,\,& \bar P, \quad
\tex(Q) = \tex(P).
\end{align}

\Proof From (\ref{ss.GEST}.\ref{GALA}) we conclude $\mu \notin
n_P(P^\times)$. Thus $Q$ is not only round and non-singular but
also anisotropic, so $\lambda_Q$ satisfying
(\ref{ss.VADAL}.\ref{EXC})$\--$(\ref{ss.VADAL}.\ref{LASC}) exists.
Since $\lambda_Q(vj) = \lambda_P(v) + \frac{1}{2} \neq
\lambda_P(u)$ for all $u,v \in V_P^\times$, by
(\ref{ss.ENLAR}.\ref{ENQU}),(\ref{ss.VADAL}.\ref{EXC}) and again
by (\ref{ss.GEST}.\ref{GALA}), we obtain \eqref{UVEPR}, hence
\eqref{OCEPR} and the first two relations of \eqref{VADAPR}, while
the last one is immediately implied by
(\ref{ss.ENLAR}.\ref{TEQU}). \hfill $\square$

\subsection{Proposition.} \label{p.LATAU} \emph{If $P$ as in}
\ref{ss.GEST}~(a) \emph{is tame and $\mu \in \mfo^\times$, then $Q
= \dla\mu\dra \otimes P$ is a non-singular and round pointed
quadratic space over $F$. Moreover the following conditions are
equivalent.}
\begin{itemize}
\item [(i)] \emph{$Q$ is anisotropic.}

\item [(ii)] \emph{$\mu \notin n_P(V_P^\times)$}.

\item [(iii)] \emph{$\bar \mu \notin n_{\bar P}(V_{\bar
P}^\times)$.}
\end{itemize}
\emph{In this case,}
\begin{align}
\label{LACTAU} \lambda_Q(u + vj)
=\,\,&\mbm\,\{\lambda_P(u),\lambda_P(v)\} &&(u,v \in P), \\
\label{OPECU} \mfo_Q = \mfo_P \oplus \mfo_P&j, \quad \mfp_Q =
\mfp_P \oplus \mfp_Pj,
\end{align}
\emph{and $Q$ is tame of ramification index $e_{Q/F} = 1$ with
$\bar Q \cong \dla\bar{\mu}\dra \otimes \bar P$.}

\Proof While the first statement is obvious, the equivalence of
(i),(ii),(iii) follows from Prop.~\ref{p.CDPOQUA}~(a) and
Cor.~\ref{c.NOCL} since $\wid(P) = \tex(P) = 0$ by
(\ref{ss.GEST}.\ref{WITE}) and tameness of $P$. If (i),(ii),(iii)
hold, then $\lambda_Q$ exists and it suffices to show that
\eqref{LACTAU} holds. Since $\lambda_Q(vj) = \lambda_P(v)$ by
(\ref{ss.ENLAR}.\ref{ENQU}),(\ref{ss.VADAL}.\ref{EXC}), we certainly
have $\lambda_Q(u + vj) \geq \mbm\,\{\lambda_P(u),\lambda_P(v)\}$.
To prove equality, Lemma~\ref{l.COVA} allows us to assume
$\lambda_P(u) = \lambda_P(v) = 0$. Here $\lambda_Q(u + vj) > 0$
would imply $n_{\bar P}(\bar u) = \bar \mu n_{\bar P}(\bar v)$,
hence $\bar \mu \in n_{\bar P}(\bar P^\times)$ since $n_{\bar P}$ is
round, and we obtain a contradiction to (iii). \hfill $\square$ \lz

\noindent The remaining cases where $P$ is wild and $\mu \in \mfo$
is a unit are much more troublesome.

\subsection{Some easy reductions.} \label{ss.EARE} For the rest of
this section, we assume that $P$ as given in \ref{ss.GEST}~(a) is
wild, so $\wid(P) = \tex(P) > 0$ by Prop.~\ref{p.PROTE}~(a). We
are interested in the pointed quadratic spaces $\dla\mu\dra
\otimes P$, $\mu \in \mfo^\times$, only when they are anisotropic.
By Prop.~\ref{p.CDPOQUA}~(a) and the local norm
theorem~\ref{t.HENOTH}, this is equivalent to $\mu$ having norm
exponent $\leq 2\tex(P)$, so by Prop.~\ref{p.PRONE}~(b) it will be
enough to consider units in $\mfo$ having the form $\mu = (1 -
\pi^d\beta)n_P(v)$, where $d \in \IZ$ satisfies $0 \leq d \leq
2\tex(P)$, $\beta \in \mfo^\times$ and $v \in \mfo_P^\times$. Here
Prop.~\ref{p.CDPOQUA}~(b) allows us to assume $v = 1_P$. We are
thus reduced to working with scalars $\mu$ that may be written as
\begin{align}
\label{REDSCA} \mu = 1 - \pi^d\beta, \quad d \in \IZ, \quad 0 \leq
d \leq2\tex(P), \quad \beta \in \mfo^\times.
\end{align}
Setting
\begin{align}
\label{EMTHE} 0 \leq m := \floor{\frac{d}{2}} \leq \tex(P), \Theta_d
:= \pi^{-m}(1_P + j) \in V_Q = V_P + V_Pj,
\end{align}
we define, inspired by the Cayley-Dickson construction of conic
algebras (cf. (\ref{ss.CDCO}.\ref{CDMU})),
\begin{align}
\label{VETHE} v\Theta_d := \pi^{-m}(v + vj) \quad(v \in V_P),
\quad V_P\Theta_d := \{v\Theta_d \mid v \in V_P\}
\end{align}
and obtain after a straightforward computation, involving
(\ref{ss.ENLAR}.\ref{ENQU}),(\ref{ss.ENLAR}.\ref{TEQU}),
\begin{align}
\label{DECTH} V_Q =\,\,&V_P \oplus V_P\Theta_d, \\
\label{NODECTH}n_Q(u + v\Theta_d) =\,\,&n_P(u) + \pi^{-m}n_P(u,v)
+ \pi^{d-2m}\beta n_P(v),  \\
\label{TRDECTH} t_Q(u + v\Theta_d) =\,\,&t_P(u) + \pi^{-m}t_P(v)
\end{align}
for all $v \in V_P$.

\subsection{Pointed quadratic residue spaces and inseparable
extensions.} \label{ss.POQUINS} (cf. \cite[Remark~10.4]{MR2427530})
For $P$ as in \ref{ss.GEST},\ref{ss.EARE}, we claim \emph{that
$V_{\bar P}$ carries a unique structure of a purely inseparable
extension field over $\bar F$ having exponent at most $1$ such that
$n_{\bar P}(u^\prime) = u^{\prime 2}$ for all $u^\prime \in V_{\bar
P}$.} \index{pointed quadratic space!wild!purely inseparable field
structure}To see this, it suffices to note that $n_{\bar P}$ is
round and anisotropic with $\partial n_{\bar P} = 0$, making
$n_{\bar P}(V_{\bar P})$ a subfield of $\bar F$, so an $\bar
F$-bilinear multiplication $V_{\bar P} \times V_{\bar P} \to V_{\bar
P}$, $(u^\prime,v^\prime) \mapsto u^\prime v^\prime$, gives a purely
inseparable extension field structure as indicated iff $n_{\bar
P}(u^\prime v^\prime) = n_{\bar P}(u^\prime)n_{\bar P}(v^\prime)$
for all $u^\prime,v^\prime \in V_{\bar P}$. The purely inseparable
extension field thus constructed will again be denoted by $V_{\bar
P}$ if there is no danger of confusion.

\subsection{Proposition.} \label{p.CDUNOD} \emph{Let $d$ be an}
odd \emph{integer with $0 \leq d \leq 2\tex(P)$ and $\beta \in
\mfo^\times$. Then}
\[
\mu := 1 - \pi^d\beta \in \mfo^\times
\]
\emph{and $Q := \dla\mu\dra \otimes P$ is a non-singular, round
and anisotropic pointed quadratic space over $F$. Moreover, $Q$ is
wild and}
\begin{align}
\label{CDPROD} \Pi := \Theta_d = \pi^{-\frac{d-1}{2}}(1_B + j) \in
V_Q
\end{align}
\emph{is a prime element of $\mfo_Q$ with}
\begin{align}
\label{DEPROD} n_Q(\Pi) = \pi\beta, &\quad t_Q(\Pi) =
2\pi^{-\frac{d-1}{2}},
\quad V_Q = V_P \oplus V_P\Pi, \\
\label{UVEPROD} \lambda_Q(u + v\Pi) =
&\,\,\mbm\,\{\lambda_P(u),\lambda_P(v)
+ \frac{1}{2}\} \quad (u,v \in V_P), \\
\label{OCEPROD} \mfo_Q = \mfo_P& \oplus \mfo_P\Pi, \quad \mfp_Q =
\mfp_P \oplus \mfo_P\Pi, \\
\label{TREXPROD} e_{Q/F} = 2, \quad \bar Q =\,\,&\bar P, \quad
\tex(Q) = \tex(P) - \frac{d-1}{2}.
\end{align}

\Proof By Prop.~\ref{c.ODNE}~(a), $\nex_P(\mu) = d$ is finite,
forcing $Q$ to be anisotropic. Applying
(\ref{ss.EARE}.\ref{DECTH})$\--$(\ref{ss.EARE}.\ref{TRDECTH}), we
obtain \eqref{DEPROD}; in particular, $\lambda_Q(\Pi) =
\frac{1}{2}$, so $\Pi$ is a prime element of $\mfo_Q$ and the first
formula of \eqref{TREXPROD} holds. We proceed to establish
\eqref{UVEPROD}. If $u,v \in V_P^\times$, then $\lambda_Q(u) =
\lambda_P(u)$ is an integer by (\ref{ss.GEST}.\ref{GALA}), while
$\lambda_Q(v\Pi) = \lambda_P(v) + \frac{1}{2}$ is not. This not only
proves \eqref{UVEPROD} but also \eqref{OCEPROD} and the second
formula of \eqref{TREXPROD}. The last one follows from the fact that
(\ref{ss.EARE}.\ref{TRDECTH}) establishes $\mfp^{\tex(P)} +
\mfp^{\tex(P)-\frac{d-1}{2}} = \mfp^{\tex(P)-\frac{d-1}{2}}$ as the
trace ideal of $Q$. \hfill $\square$

\subsection{Proposition.} \label{p.CDUNEV} \emph{Let $d$ be an}
even \emph{integer with $0 \leq d <2\tex(P)$ and suppose $\beta \in
\mfo$ satisfies the condition $\bar \beta \notin V_{\bar P}^2$} (cf.
\ref{ss.POQUINS}). \emph{Then}
\begin{align}
\label{UNEV} \mu = 1 - \pi^d\beta \in \mfo^\times
\end{align}
\emph{and $Q := \dla\mu\dra \otimes P$ is a non-singular, round
and anisotropic pointed quadratic space over $F$. Moreover, $Q$ is
wild and}
\begin{align}
\label{CDPREV} \Xi := \Theta_d = \pi^{-\frac{d}{2}}(1_P + j) \in V_Q
\end{align}
\emph{is a unit of $\mfo_Q$ with}
\begin{align}
\label{DEPREV}n_Q(\Xi) = \beta, \quad t_Q(&\Xi) =
2\pi^{-\frac{d}{2}}, \quad V_Q = V_P \oplus V_P\Xi, \\
\label{UVEPREV} \lambda_Q(u + v\Xi) =\,\,&
\mbm\,\{\lambda_P(u),\lambda_P(v)\} &&(u,v \in V_P), \\
\label{OCEPREV} \mfo_Q = \mfo_P \oplus \mfo_P\Xi, &\quad \mfp_Q =
\mfp_P \oplus \mfp_P\Xi, \\
\label{TREXPREV} e_{Q/F} = 1, \quad \bar Q = \dla\bar\beta\dra
\otimes \bar P, &\quad \tex(Q) = \tex(P) - \frac{d}{2}.
\end{align}

\Proof The assertion $\mu \in \mfo^\times$ is trivial for $d > 0$
but holds also for $d = 0$ since in this case $\bar \beta \notin
V_{\bar P}^2$ implies $\bar\mu = 1_{\bar F} - \bar\beta \notin
V_{\bar P}^2$ by wildness of $P$. Next we show that $Q$ is
anisotropic. Otherwise, $\mu \in n_P(V_P)$ by
Prop.~\ref{p.CDPOQUA}~(a), and Lemma~\ref{l.NOCO} yields an element
$w \in \mfo_P^\times$ with $\beta = -n_P(w) + \pi^{-m}t_P(w)$, $m =
\frac{d}{2}$, where the second summand belongs to $\mfp^{\tex(P) -
m} \subseteq \mfp$ by the hypothesis on $d$. Thus $\bar\beta =
n_{\bar P}(\bar w) = \bar w^2$, a contradiction, and we have proved
that $Q$ is indeed anisotropic.  Consulting
(\ref{ss.EARE}.\ref{DECTH}$\--$\ref{TRDECTH}) for $u = 0$, $v =
1_P$, we end up with \eqref{DEPREV}. Turning to \eqref{UVEPREV}, it
suffices to show, by Lemma~\ref{l.COVA}, that $u,v \in
\mfo_P^\times$ implies $u + v\Xi \in \mfo_Q^\times$. Otherwise,
observing (\ref{ss.EARE}.\ref{NODECTH}),
\[
\lambda_Q(u + v\Xi) = \frac{1}{2}\lambda\big(n_P(u) +
\pi^{-m}n_P(u,v) + \beta n_P(v)\big)
\]
were strictly positive, and since $\pi^{-m}n_P(u,v) \in \mfp$ by
Prop.~\ref{p.PROTE}~(c), we would again arrive at the
contradiction $\bar\beta \in V_{\bar P}^2$. Thus \eqref{UVEPREV}
holds, which directly implies \eqref{OCEPREV}, while
\eqref{TREXPREV} follows from \eqref{OCEPREV} and
(\ref{ss.EARE}.\ref{NODECTH},\ref{TRDECTH}). \hfill $\square$ \lz

\noindent In \ref{ss.EARE}, particularly
(\ref{ss.EARE}.\ref{REDSCA}), we are left with the case $d =
2\tex(P)$, which turns out to be the most delicate. In order to
get started, we require the following elementary but crucial
observations.

\subsection{Setting the stage for the case $d = 2\tex(P)$.}
\label{ss.SETEX} (a) Let $K/k$ be a purely inseparable field
extension of characteristic $2$, exponent at most $1$ and finite
degree. Consider a scalar $\alpha \in k$ and a unital linear form
$s\:K \to k$. We denote by by $Q_{K;\alpha,s}$ the pointed
quadratic space over $k$ with norm the Pfister quadratic form
$q_{K;\alpha,s}$ of \ref{ss.BILQUA} and with base point $1_K \in K
\subseteq K \oplus Kj$. Recall from Prop.~\ref{p.FOCOAL} that
$Q_{K;\alpha,s} = Q_{\mbC(K;\alpha,s)}$ is the pointed quadratic
space corresponding to the flexible conic algebra
$\mbC(K;\alpha,s)$ that arises form $K,\alpha,s$ by means of the
non-orthogonal Cayley-Dickson construction. $\ssp$ \\
(b) Put $m := \tex(P)$ and let $w_0$ be a normalized trace
generator of $P$. Then $w_0 \in \mfo_P^\times$ by
(\ref{ss.TRAG}.\ref{ESTRAG}) and the map $s_{w_0}\:V_P \to F$, $u
\mapsto \pi^{-m}n_P(u,w_0)$ is a unital linear form with
$s_{w_0}(\mfo_P) \subseteq \mfo$, $s_{w_0}(\mfp_P) \subseteq \mfp$
by Prop.~\ref{p.PROTE}~(c) and since $\mfp_P = \mfp\mfo_P$. Thus
we obtain an induced unital linear form
\begin{align}
\label{INDUN} \bar s_{w_0}\:V_{\bar P} \longrightarrow \bar F,
\quad \bar u \longmapsto \bar s_{w_0}(\bar u) =
\overline{\pi^{-m}n_P(u,w_0)},
\end{align}
and given any $\alpha^\prime \in \bar F$, the notational
conventions of (a) apply when $V_{\bar P}/\bar F$ is viewed via
\ref{ss.POQUINS} as a purely inseparable field extension of
exponent at most $1$.

\subsection{Theorem.} \label{t.CHAGORE} \emph{With the notations
of} \ref{ss.SETEX}, \emph{let}
\begin{align}
\label{MUTEX} \beta \in \mfo, \quad \mu := 1 - \pi^{2\tex(P)}\beta
\in \mfo^\times, \quad \beta_0 := n_P(w_0)\beta \in \mfo, \quad Q
:=\dla\mu\dra \otimes P.
\end{align}
\emph{Then the following conditions are equivalent.}
\begin{itemize}
\item [(i)] \emph{$Q$ is anisotropic and unramified.}

\item [(ii)] \emph{$Q$ is anisotropic and tame.}

\item [(iii)] \emph{$Q$ is anisotropic.}

\item [(iv)] $\nex_P(\mu) = 2\tex(P)$.
\end{itemize}
\emph{Moreover, if $P$ is a pointed Pfister quadratic space, then
these conditions are also equivalent to}
\begin{itemize}
\item [(v)] $\overline{\beta_0} \notin \Im(\wp_{V_{\bar P},\bar
s_{w_0}})$.
\end{itemize}

\Proof In this section, we will not be able to give the proof in
full but must restrict ourselves to showing the equivalence of
(i)$\--$(iv), relegating the rest to the next section. Since the
implications (i) $\Rightarrow$ (ii) $\Rightarrow$ (iii) are obvious,
it suffices to show (iv) $\Leftrightarrow$ (iii) $\Rightarrow$ (i).

(iv) $\Longleftrightarrow$ (iii). From \ref{ss.NE} we deduce
$\nex_P(\mu) \geq 2\tex(P)$. Hence (iv) holds iff $\nex_P(\mu)
\leq 2\tex(P)$ iff $Q$ is anisotropic by the local norm theorem
\ref{t.HENOTH} and Prop.~\ref{p.PRONE}~(a).

(iii) $\Longrightarrow$ (i). Setting $d := 2\tex(P)$, we apply
(\ref{ss.EARE}.\ref{NODECTH},\ref{TRDECTH}) for $u = 0$, $v = w_0$
and obtain $\Xi_0 := w_0\Theta_d \in \mfo_Q^\times$ and
$t_Q(\Xi_0) = 1$. Thus $Q$ is tame, and since $P$ is wild, we
conclude that $\bar P$ is a \emph{proper} pointed quadratic
subspace of $\bar Q$. Hence $e_{Q/F} = 1$, forcing $Q$ to be
unramified. \hfill $\square$ \lz

\noindent Given $\beta \in \mfo$, $d \in \IZ$ with $0 \leq d \leq
2\tex(P)$, it follows from \ref{ss.NE} that $\mu := 1 -
\pi^d\beta$ has norm exponent at least $d$, and if $d$ is odd,
Prop.~\ref{c.ODNE}~(a) yields a characterization in terms of
$\beta$ when equality holds. While a similar characterization for
$d = 2\tex(P)$ is presented in Thm.~\ref{t.CHAGORE}~(v) (though as
yet unproved), we are now able to provide one for $d$ even, $d <
2\tex(P)$.

\subsection{Corollary.} \label{c.CHARNEXOD} \emph{Let $d$ be an
even integer such that $0 \leq d < 2\tex(P)$.} $\ssp$ \\
(a) \emph{With $\beta \in \mfo$ and $\mu := 1 - \pi^d\beta$, the
following conditions are equivalent.}
\begin{itemize}
\item [(i)] \emph{$\mu \in \mfo^\times$ and $\nex_P(\mu) = d$.}

\item [(ii)] \emph{$\bar \beta \notin V_{\bar P}^2$.}
\end{itemize}
(b) \emph{If $d > 0$, $\beta_i \in \mfo^\times$ and $\mu_i := 1 -
\pi^d\beta_i \in \mfo^\times$ for $i = 1,2$, then}
\[
\dla\mu_1\dra \otimes P \cong \dla\mu_2\dra \otimes P
\Longrightarrow \overline{\beta_1} \equiv \overline{\beta_2} \bmod
\,V_{\bar P}^2.
\]

\Proof (a) (i) $\Longrightarrow$ (ii). We have $\beta \in
\mfo^\times$ by Prop.~\ref{p.PRONE}~(b). For $d = 0$ the assertion
follows from Prop.~\ref{p.PRONE}~(c). Now suppose $d \ne 0$. Arguing
indirectly, we assume that there exists an element $w \in
\mfo_P^\times$ with $\bar \beta = \bar w^2 = \overline{n_P(w)}$.
Then $n_P(w) = \beta + \pi\beta^\prime$ for some $\beta^\prime \in
\mfo$, and Lemma~\ref{l.TECH} with $m := \frac{d}{2}$ yields
\[
\mu = 1 - \pi^d\beta = (1 -
\pi^{d+1}\beta^{\prime\prime})n_P(v^{\prime\prime})^{-1}
\]
where, setting $r := \tex(P)$,
\begin{align*}
\beta^{\prime\prime} = -\beta^\prime + \pi^{r-m-1}\alpha -
\pi^{r+m-1}\alpha\beta + \pi^{2m-1}\beta n_P(w) \in \mfo, \quad
v^{\prime\prime} = 1_P - \pi^mw \in \mfo_P^\times
\end{align*}
since $0 < m < r$.  Now roundness of $P$ and the definition of the
norm exponent imply $\nex_P(\mu) \geq d + 1$, a contradiction.

(ii) $\Longrightarrow$ (i). By (\ref{p.CDUNEV}.\ref{UNEV}) we have
$\mu \in \mfo^\times$ and $Q$ is anisotropic by
Prop.~\ref{p.CDUNEV}, forcing $\mu \notin n_P(V_P^\times)$
(Prop.~\ref{p.CDPOQUA}~(a)) and $d \leq d^\prime := \nex_P(\mu) \leq
2\tex(P)$ (Thm.~\ref{t.HENOTH}). Furthermore, by
Prop.~\ref{p.PRONE}~(b), $\mu = \mu^\prime n_P(v^\prime)$,
$\mu^\prime = 1 - \pi^{d^\prime}\beta^\prime$ for some $\beta^\prime
\in \mfo^\times$, $v^\prime \in \mfo_P^\times$. In particular,
$\mu^\prime \in \mfo^\times$ and $Q \cong Q^\prime :=
\dla\mu^\prime\dra \otimes P$. This implies $e_{Q^\prime/F} =
e_{Q/F} = 1$ by (\ref{p.CDUNEV}.\ref{TREXPREV}), so $d^\prime <
2\tex(P)$ (Thm.~\ref{t.CHAGORE}) is even (Prop.~\ref{p.CDUNOD}), and
we are allowed to apply Prop.~\ref{p.CDUNEV} to $Q^\prime$ by the
implication (i) $\Rightarrow$ (ii) already established. Thus
(\ref{p.CDUNEV}.\ref{TREXPREV}) yields $d = d^\prime$.

(b) We put $Q_i := \dla\mu_i\dra \otimes P$ for $i = 1,2$. If
$V_{\bar P}^2$ contains $\overline{\beta_1}$ but not
$\overline{\beta_2}$, then $\mu_1,\mu_2$ have different norm
exponents by (a), so $Q_1,Q_2$ cannot be isomorphic
(Props.~\ref{p.PRONE}~(d),~\ref{p.CDPOQUA}~(b)). We may therefore
assume $\overline{\beta_1},\overline{\beta_2} \notin V_{\bar
P}^2$. As in the proof of Cor.~\ref{c.ODNE} we have $\mu :=
\mu_1\mu_2 = 1 - \pi^d\beta$, $\beta := \beta_1 + \beta_2 -
\pi^d\beta_1\beta_2 \in \mfo$. Assuming $\bar \beta =
\overline{\beta_1} - \overline{\beta_2} \notin V_{\bar P}^2$ would
force $\dla\mu\dra \otimes P$ to be anisotropic by
Prop.~\ref{p.CDUNEV}, hence $\mu_1,\mu_2$ to fall into distinct
norm classes relative to $P$. But then $Q_1,Q_2$ would not be
isomorphic, a contradiction. \hfill $\square$ \lz

%\noindent For a counter example to the converse of
%Cor.~\ref{c.CHARNEXOD}~(b), see the following example.

\subsection{Examples.} \label{e.NEXP} (a) Let $m$ be an integer
with $0 \leq m < \tex(P)$ and suppose we are given an element
$\gamma \in \mfo^\times$ such that $\bar \gamma \in V_{\bar P}^2$.
If $\mu := 1 - \pi^{2m}\gamma$ is a unit in $\mfo$ (automatic unless
$m = 0$), Cor.~\ref{c.CHARNEXOD}~(a) implies $\nex_P(\mu) > 2m$, so
it is a natural question to ask whether a more precise estimate for
the norm exponent of $\mu$ can be given. Unless specific properties
of $\gamma$ are taken into account, the answer is no. To see this,
let $d \in \IZ$ with $d > 2m$, $\beta \in \mfo$, and $w \in
\mfo_P^\times$ with $\bar w \neq 1_{\bar P}$. Then
Lemma~\ref{l.TECH} implies
\[
1 - \pi^{2m}\gamma \equiv 1 - \pi^d\beta \bmod\,n_P(P^\times),
\quad \bar \gamma = \bar w^2
\]
for
\begin{align*}
\gamma = -n_P(w) + \pi^{d-2m}\beta + \pi^{\tex(P) - m}\alpha -
\pi^{\tex(P)+d-m}\alpha\beta + \pi^d\beta n_P(w) \in \mfo^\times.
\end{align*}
Hence $\nex_P(1 - \pi^{2m}\gamma) = \infty$ for $d > 2\tex(P)$ by
the local norm theorem \ref{t.HENOTH}, and, by
Cors.~\ref{c.ODNE},\ref{c.CHARNEXOD}, $\beta$ may be so chosen that
$\nex_P(1 - \pi^{2m}\gamma)$ attains any finite pre-assigned value
$d$ with $2m < d < 2\tex(P)$ provided $V_{\bar P}^2 \neq
\bar F$. $\ssp$ \\
(b) Let $d \in \IZ$, $0 < d \leq 2\tex(P)$, $\beta,\beta^\prime
\in \mfo^\times$ and put $\mu := 1 - \pi^d\beta$, $\mu^\prime = 1
- \pi^d\beta^\prime \in \mfo^\times$. We wish to refute the
converse of Cors.~\ref{c.ODNE}~(b) and \ref{c.CHARNEXOD}~(b) by
showing that $\bar \beta = \bar \beta^\prime$ \emph{does not}
imply $\mu \equiv \mu^\prime \bmod\,n_P(V_P^\times)$. Indeed, if
$\beta \neq \beta^\prime$, then $\bar \beta = \bar \beta^\prime$
amounts to the same as $\mu^\prime = 1 - \pi^d\beta - \pi^q\gamma$
for some integer $q > d$ and some $\gamma \in \mfo^\times$,
allowing us to conclude
\[
\mu\mu^\prime = \mu^2(1 - \pi^q\mu^{-1}\gamma) \equiv 1 -
\pi^q\mu^{-1}\gamma \bmod\,n_P(V_P^\times).
\]
Therefore,
\begin{itemize}
\item $\mu \equiv \mu^\prime \bmod\,n_P(V_P^\times)$ for $q >
2\tex(P)$ (Thm.~\ref{t.HENOTH}),
\end{itemize}
but
\begin{itemize}
\item $\mu \not\equiv \mu^\prime \bmod\,n_P(V_P^\times)$ for $d <
q \leq 2\tex(P)$, provided $\bar \gamma \notin V_{\bar P}^2$ if $q <
2\tex(P)$ is even (Cors.~\ref{c.ODNE},\ref{c.CHARNEXOD}), and $\bar
w_0^2\bar \gamma \notin \Im(\wp_{V_{\bar P};\bar s_{w_0}})$ if $d =
2\tex(P)$ (Thm.~\ref{t.CHAGORE}).
\end{itemize}
\vspace{-7pt} \hfill $\square$ \lz

\noindent The preceding results on the behavior of the
ramification index, the pointed quadratic residue space and the
trace exponent under the passage from $P$ to $\dla\mu\dra \otimes
P$, $\mu \in F^\times$, can be stated in a particularly concise
way when addressed to pointed Pfister quadratic spaces by
combining them with the embedding property of
Prop.~\ref{p.EMPROQUA}. In order to do so, we introduce the
following terminology.

\subsection{Scalars of standard type.} \label{ss.CABARE}
We say a scalar $\mu \in F$ has \emph{standard
type}\index{standard type, scalar of} relative (or with respect)
to $P$ if it satisfies one of the following mutually exclusive
conditions.
\begin{itemize}
\item [(a)] $\mu$ is a prime element of $\mfo$ (possibly distinct
from $\pi$).

\item [(b)] $\mu = 1 - \pi^d\beta$ for some \emph{odd} integer $d$
with $0 \leq d < 2\tex(B)$ and some $\beta \in \mfo^\times$.

\item [(c)] $\mu = 1 - \pi^d\beta$ for some \emph{even} integer
$d$ with $0 \leq d < 2\,\tex(B)$ and some $\beta \in \mfo$ with
$\bar \beta \notin V_{\bar P}^2$.
\end{itemize}

\subsection{Theorem.} \label{t.ENBARE} \emph{Let $P$ be a pointed
$n$-Pfister quadratic space over $F$ that is anisotropic and wild
of ramification index $e_{P/F} = 1$. For $Q$ to be a wild
anisotropic pointed $(n+1)$-Pfister quadratic space over $F$ into
which $P$ embeds as a pointed quadratic subspace it is necessary
and sufficient that $Q$ be a pointed quadratic space isomorphic to
$\dla\mu\dra \otimes P$, for some scalar $\mu \in F$ of standard
type relative to $P$. In this case, precisely one of the following
implications holds.}
\begin{itemize}
\item [(a)] \emph{If $\mu$ is a prime element in $\mfo$, then}
\[
e_{Q/F} = 2, \quad \bar Q \cong \bar P, \quad \tex(Q) = \tex(P).
\]
\item [(b)] \emph{If $\mu = 1 - \pi^d\beta$ for some} odd
\emph{integer $d$ with $0 \leq d <2\tex(P)$ and some $\beta \in
\mfo^\times$, then}
\[
e_{Q/F} = 2, \quad \bar Q \cong \bar P, \quad \tex(Q) = \tex(P) -
\frac{d-1}{2}.
\]
\item [(c)] \emph{If $\mu = 1 - \pi^d\beta$ for some} even
\emph{integer $d$, $0 \leq d < 2\tex(P)$ and some $\beta \in \mfo$
with $\bar\beta \notin V_{\bar P}^2$, then}
\[
e_{Q/F} = 1, \quad \bar Q \cong \dla\bar\beta\dra \otimes \bar P,
\quad \tex(Q) = \tex(P) - \frac{d}{2}.
\]
\end{itemize}

\Proof By Props.~\ref{p.CDPR},\ref{p.CDUNOD},\ref{p.CDUNEV}, the
condition is sufficient, and (a)$\--$(c) hold. Conversely, suppose
$Q$ is a pointed anisotropic wild $(n+1)$-Pfister quadratic space
over $F$ containing $P$ as a pointed quadratic subspace. Up to
isomorphism, $Q = \dla\mu\dra \otimes P$ for some $\mu \in F^\times$
by the embedding property (Prop.~\ref{p.EMPROQUA}), where the
reduction of \ref{ss.GEST}~(b) allows us to assume that $\mu \in
\mfo^\times$ is a unit in $\mfo$. $Q$ being anisotropic implies $0
\leq d := \nex_P(\mu) \leq 2\tex(P)$ by the local norm
theorem~\ref{t.HENOTH} and without loss $\mu = 1 - \pi^d\beta$ for
some $\beta \in \mfo$. Since $Q$ is wild, we conclude $d < 2\tex(P)$
from the part of Thm.~\ref{t.CHAGORE} already established, and if
$d$ is odd, then $\beta \in \mfo^\times$ (Prop.~\ref{c.ODNE}~(a)),
while if $d$ is even, then $\bar\beta \notin V_{\bar P}^2$
(Cor.~\ref{c.CHARNEXOD}). In any event, $\mu$ has standard type
relative to $P$. \hfill $\square$

\subsection{Corollary.} \label{c.EQTE} \emph{With $P$ as in
Thm.~}\ref{t.ENBARE}, \emph{let $Q$ be a
pointed quadratic space over $F$.} $\ssp$ \\
(a) \emph{The following conditions are equivalent.}
\begin{itemize}
\item [(i)] \emph{$Q$ is an anisotropic pointed $(n+1)$-Pfister
quadratic space into which $P$ embeds as a pointed quadratic
subspace such that $e_{Q/F} = 2$ and $\tex(Q) = \tex(P)$.}

\item [(ii)] \emph{$Q \cong \dla\pi\beta\dra \otimes P$ or $Q
\cong \dla 1 - \pi\beta\dra \otimes P$ for some $\beta \in
\mfo^\times$.}
\end{itemize}
\emph{In this case $\bar Q \cong \bar P$.} $\ssp$ \\
(b) \emph{The following conditions are equivalent.}
\begin{itemize}
\item [(i)] \emph{$Q$ is an anisotropic pointed $(n+1)$-Pfister
quadratic space into which $P$ embeds as a pointed quadratic
subspace such such that $e_{Q/F} = 1$ and $\tex(Q) = \tex(P)$.}

\item [(ii)] \emph{$Q \cong \dla\mu\dra \otimes P$ for some $\mu
\in \mfo$ with $\bar \mu \notin V_{\bar P}^2$.}
\end{itemize}
\emph{In this case, $\bar Q \cong \dla\bar\mu\dra \otimes P$.}

\Proof In (a) and (b), condition (i) implies that $Q$ is wild
(Prop.~\ref{p.PROTE}~(a)). Hence the assertions follow immediately
from Thm.~\ref{t.ENBARE}. \hfill $\square$ \lz

\noindent There is an analogue of Thm.~\ref{t.ENBARE} dealing with
tame rather than wild enlargements of pointed $n$-Pfister quadratic
spaces. We omit the proof since it proceeds along the same lines as
the one of Thm.~\ref{t.ENBARE}, applying Thm.~\ref{t.CHAGORE} in
full rather than Props.~\ref{p.CDPR},\ref{p.CDUNOD},\ref{p.CDUNEV}.

\subsection{Theorem.} \label{t.ENGORE} \emph{Keeping the notations
of} \ref{ss.SETEX}~(b), \emph{let $P$ be a pointed $n$-Pfister
quadratic space over $F$ that is anisotropic and wild of
ramification index $e_{P/F} = 1$. For $Q$ to be a tame and
anisotropic $(n+1)$-Pfister pointed quadratic space over $F$ into
which $P$ embeds as a pointed quadratic subspace it is necessary
and sufficient that $Q$ be a pointed quadratic $F$-space and there
exist an element $\beta \in \mfo$ with}
\[
Q \cong \dla\mu\dra \otimes P, \quad \mu := 1 - \pi^{2\tex(B)}\beta,
\quad \overline{\beta_0} \notin \Im(\wp_{V_{\bar P},\bar s_{w_0}}),
\quad \beta_0 := n_P(w_0)\beta.
\]
\hfill $\square$

\section{$\lambda$-normed and $\lambda$-valued conic algebras.}
\label{s.LANOVA} In order to illuminate the intuitive background of
the present section, we recall the notion of an absolute-valued
algebra. Following Albert \cite{MR0020550} (see also  Palacios
\cite{MR2162414}), an absolute-valued algebra is a non-associative
real algebra $A$ equipped with a norm $x \mapsto \|x\|$ that permits
composition: $\|xy\| = \|x\|\,\|y\|$ for all $x,y \in A$. Since
absolute valued algebras obviously have no zero divisors, the
finite-dimensional ones are division algebras, hence exist only in
dimensions $1,2,4,8$ (Albert \cite{MR0020550} gave an ad-hoc proof
of this result, the Bott-Kervaire-Milnor thoerem not having been
known at the time). By contrast, natural analogues of
absolute-valued algebras over $2$-Henselian fields will be discussed
in the present section that exist in all dimensions $2^n$, $n =
0,1,2\dots$ .

Throughout we continue to work over a fixed $2$-Henselian field
$F$ as in \ref{ss.HEFI} and alert the reader to the terminological
conventions of \ref{ss.CONCON}. All vector spaces, algebras, etc.\
over $F$ are tacitly assumed to be finite-dimensional.

\subsection{The basic concepts.} \label{ss.BACO} A conic algebra
$C$ over $F$ is said to be \emph{$\lambda$-normed}\index{conic
algebra!$\lambda$-normed} if the following conditions hold.
\begin{itemize}
\item [(i)] $C$ is non-singular, round and anisotropic.

\item [(ii)] $\lambda_C$ is sub-multiplicative: $\lambda_C(xy)
\geq \lambda_C(x) + \lambda_C(y)$ for all $x,y \in C$.
\end{itemize}
We speak of a \emph{$\lambda$-valued} conic algebra\index{conic
algebra!$\lambda$-valued} if $C$ satisfies (i), and if instead of
(ii) the following stronger condition holds:
\begin{itemize}
\item [(iii)] $\lambda_C$ is multiplicative: $\lambda_C(xy) =
\lambda_C(x) + \lambda_C(y)$ for all $x,y \in C$.
\end{itemize}
The norm of a $\lambda$-valued conic algebra $C$ over $F$ will typically not permit composition (for example, if the dimension of
$C$ differs from $1,2,4,8$) but, remarkably, its failure to do so
is not detected by $\lambda$ since $\lambda_C$ being
multiplicative  by (\ref{ss.CONCON}.\ref{EXCO}) amounts to
$\lambda(n_C(xy)) = \lambda(n_C(x)n_C(y))$ for all $x,y \in C$.
This looks like a pretty far-fetched phenomenon but, in fact,
turns out to be quite common.

We start with a trivial but useful observation.

\subsection{Proposition.} \label{p.ELCOMCOM}
(a) \emph{$\lambda$-valued conic algebras over $F$ are division
algebras.}  \\
(b) \emph{A composition algebra over $F$ is a $\lambda$-valued
conic algebra if and only if it is a division algebra.} \hfill
$\square$

\begin{prop} \label{p.VARIID} Let $C$ be a
non-singular, round and anisotropic conic algebra over $F$. Then
$C$ is $\lambda$-normed if and only if $\mfo_C \subseteq C$ is an
$\mfo$-subalgebra and $\mfp_C \subseteq \mfo_C$ is an ideal with
$\mfp_C^2 \subseteq \mfp\mfo_C$. In this case, $\bar C :=
\mfo_C/\mfp_C$ is a conic algebra over $\bar F$ whose norm, trace
and conjugation are given by the formulas
\begin{align}
\label{VARINO} n_{\bar C}(\bar x) =\,\,&\overline{n_C(x)}, \\
\label{VARITR} t_{\bar C}(\bar x) =\,\,&\overline{t_C(x)}, \\
\label{VARICO} \left( \bar x\right)^\ast =\,\,& \overline{\left( x^* \right)}
\end{align}
for all $x \in \mfo_C$. Moreover, the norm of $\bar C$ is
round and anisotropic.
\end{prop}

\Proof By Lemma~\ref{l.COVA}, sub-multiplicativity of $\lambda_C$
amounts to
\begin{align*}
\lambda_C(xy) \geq \lambda_C(x) + \lambda_C(y) &&(x,y \in C, \quad
0 \leq \lambda_C(x),\lambda_C(y) \leq \frac{1}{2}).
\end{align*}
The first part of the proposition follows from this at once. The
second part is a restatement of Prop.~\ref{p.PORESQUA}. \hfill
$\square$

\Remark The conic algebra $\bar C = \mfo_C/\mfp_C$ is called the
\emph{residue algebra} of $C$\index{conic algebra!residue algebra}.
If $C$ is wild, we do not know whether this residue algebra
always agrees with the purely inseparable extension field of
$\bar F$ attached to $C$ via \ref{ss.POQUINS}, though it does if $C$
is a composition algebra \cite[Prop.~1]{MR51:635}.

\subsection{Corollary.} \label{c.LAVADIV} \emph{With the notations
of Prop.~}\ref{p.VARIID} \emph{suppose in addition that $C$ has
ramification index $e_{C/F} = 1$. Then}: $\ssp$ \\
(a) \emph{$C$ is $\lambda$-normed if and only if $\mfo_C \subseteq
C$ is an $\mfo$-subalgebra.} $\ssp$ \\
(b) \emph{$C$ is $\lambda$-valued if and only if $C$ is
$\lambda$-normed and $\bar C$ is a division algebra.}

\Proof (a) $\mfp_C = \mfp\mfo_C$.

(b) Consulting Lemma~\ref{l.COVA} again, $\lambda_C$ is
multiplicative iff $\lambda_C(xy) = 0$ for all $x,y \in
\mfo_C^\times$ iff $C$ is $\lambda$-normed and $\bar C$ is a
division algebra. \hfill $\square$

We now proceed to re-examine the main results of the preceding
section within the framework of $\lambda$-normed and
$\lambda$-valued conic algebras.

\subsection{Convention.} \label{ss.LACON} For the remainder of
this section, we fix a $\lambda$-normed conic algebra $B$ over $F$
having ramification index $e_{B/F} = 1$.

\subsection{Proposition.} \label{p.CONPR} \emph{If $\mu$ is a
prime element in $\mfo$, then $C := \mbC(B,\mu)$ is a
$\lambda$-normed conic algebra over $F$ with $\bar C = \bar B$ as
conic $\bar F$-algebras. Moreover, $C$ is $\lambda$-valued if and
only if $B$ is $\lambda$-valued.}

\Proof By Prop.~\ref{p.CDPR}, $C$ is non-singular, round and
anisotropic. Combining Prop.~\ref{p.VARIID} with
(\ref{ss.CDCO}.\ref{CDMU}), we conclude that $C$ is
$\lambda$-normed. It remains to show that if $B$ is
$\lambda$-valued, so is $C$. By Lemma~\ref{l.COVA}, we must show
$\lambda_C(x_1x_2) = \lambda_C(x_1) + \lambda_C(x_2)$ for all $x_i =
u_i + v_ij \in C$, $u_i,v_i \in B$, $0 \leq \lambda_C(x_i) \leq
\frac{1}{2}$, $i = 1,2$. There are four cases: (i) $\lambda_C(x_1) =
\lambda_C(x_2) = 0$, (ii) $\lambda_C(x_1) = 0$, $\lambda_C(x_2) =
\frac{1}{2}$, (iii) $\lambda_C(x_1) = \frac{1}{2}$, $\lambda_C(x_2)
= 0$, (iv) $\lambda_C(x_1) = \lambda_C(x_2) = \frac{1}{2}$. We only
treat (iv) and leave the other three cases to the reader. From
(\ref{p.CDPR}.\ref{UVEPR}) we deduce $u_i \in \mfp_B$, $v_i \in
\mfo_B^\times$ for $i = 1,2$ and (\ref{ss.CDCO}.\ref{CDMU}) yields
\[
x_1x_2 = u + vj, \quad u = u_1u_2 + \mu v_2^\ast v_1, \quad v =
v_2u_1 + v_1u_2^\ast,
\]
where $\lambda_B(u) = 1$, $\lambda_B(v) \geq 1$, hence
$\lambda_C(x_1x_2) = 1 = \lambda_C(x_1) + \lambda_C(x_2)$ \hfill
$\square$

\subsection{Example.} \label{e.LAPR} Specializing Prop.~\ref{p.CONPR}
to (iterated) Laurent series fields of characteristic not $2$, we
recover examples of $\lambda$-valued conic algebras that
originally go back to to Brown \cite[pp.~421-422]{MR0215891}.
In a slightly more general vein, let $k$ be any field, $L/k$ a
separable quadratic field extension and write
\begin{align*}
A = \mbC(L;\mu_1,\dots,\mu_{n-1}) &&(n \in \IZ,\;n \geq 1)
\end{align*}
for the $k$-algebra arising from $L$ and scalars
$\mu_1,\dots,\mu_{n-1} \in k^\times$ by means of the
Cayley-Dickson process as in \ref{ss.CDP}. Then $A$ is a flexible
conic algebra with norm an $n$-Pfister quadratic form. We now
assume that $A$ is a division algebra, forcing $A$ to be
non-singular, round and anisotropic. Consider the field $F =
k((\bft))$ of formal Laurent series in a variable $\bft$ with
coefficients in $k$, which is complete and therefore Henselian
under the standard discrete valuation $\lambda\:F \to \IZ_\infty$.
Setting
\[
B := A \otimes_k F = A((\bft)),
\]
we obtain a flexible conic division $F$-algebra whose norm is an
anisotropic $n$-Pfister quadratic form over $F$. Using
(\ref{ss.VADAL}.\ref{EXC}), a straightforward verification shows
\begin{align*}
\lambda_B\Big(\sum_{r \gg-\infty}^\infty\,a_r\bft^r\Big) = \mbm\,\{r
\in \IZ \mid a_r \neq 0\} &&(a_r \in A,\;r \in \IZ),
\end{align*}
which immediately implies that $B$ is an unramified $\lambda$-valued
conic algebra over $F$. By Prop.~\ref{p.CONPR} we thus find in $C :=
\mbC(B,\bft)$ a $\lambda$-valued conic algebra over $F$ having
dimension $2^{n+1}$. Starting from $A = L$ (i.e., from $n = 1$)
and continuing in this way, we obtain $\lambda$-valued conic
algebras over appropriate iterated Laurent series fields
in all dimensions $2^n$, $n = 0,1,2,\dots$.

\subsection{Example.} \label{e.LAPRESP} We specify
Example~\ref{e.LAPR} a bit further by setting $k = \IR$, $L =
\IC$, $n = 3$, $\mu_1 = \mu_2 = -1$. Then $A = \IO$, the real
algebra of Graves-Cayley octonions, and $B = \IO((\bft))$ is the
unique unramified octonion division algebra over $F =
\IR((\bft))$. Moreover, the $16$-dimensional conic division
algebra
\[
C = \mbC(B,\bft) = \mbC(F;-1,-1,-1,\bft)
\]
over $F$ contains $B^\prime := \mbC(F;-1,-1,\bft)$ as a ramified
octonion subalgebra. In particular, $B$ and $B^\prime$ are not
isomorphic, allowing us to conclude from \cite[Thm.~2]{MR0215891}
that the subalgebra $B^\prime \subseteq C$ does \emph{not} satisfy
the embedding property \ref{ss.EMPRO}.

\subsection{Proposition.} \label{p.CONTAU} \emph{If $B$ as in}
\ref{ss.LACON} \emph{is tame and $\mu \in \mfo^\times \setminus
n_B(B^\times)$, then $C := \mbC(B,\mu)$ is a $\lambda$-normed
conic algebra over $F$ with $\bar C = \mbC(\bar B,\bar\mu)$ as
conic $\bar F$-algebras. Moreover, $C$ is $\lambda$-valued if and
only if $\bar C$ is a division algebra.}

\Proof By Prop.~\ref{p.LATAU}, $C$ is a tame non-singular, round
and anisotropic conic algebra over $F$ with $e_{C/F} = 1$.
Moreover, $\mfo_B$ being a $\mfo$-subalgebra of $B$ by
Prop.~\ref{p.VARIID}, we conclude from (\ref{p.LATAU}.\ref{OPECU})
that $\mfo_C$ is an $\mfo$-subalgebra of $C$. Now everything
follows from Cor.~\ref{c.LAVADIV}. \hfill $\square$ \lz

\noindent Dealing with the case that $B$ as in \ref{ss.LACON} is
wild turns out to be more troublesome. We not only need a few
preparations but also have to add an extra hypothesis by requiring
that the conic algebras involved be flexible.

\subsection{Proposition.} \label{p.ELAC} \emph{Let $C$ be a
flexible $\lambda$-normed conic algebra over $F$. Then}
\begin{align}
\label{LACT} \lambda\big(t_C(x)\big) \geq\,\,&\wid(C) +
\lambda_C(x), \\
\label{LACCO} \lambda_C([x_1,x_2]) \geq\,\,&\wid(C) +
\lambda_C(x_1) + \lambda_C(x_2), \\
\label{LACAS} \lambda_C(x - x^\ast) \geq\,\,&\wid(C) +
\lambda_C(x)
\end{align}
\emph{for all $x,x_1,x_2 \in C^\times$.}

\Proof \eqref{LACT} follows immediately from
(\ref{ss.WID}.\ref{DEWI}). To establish \eqref{LACCO}, we combine
\eqref{LACT} with
(\ref{p.NOCOMT}.\ref{NOCOMT}),(\ref{ss.VADAL}.\ref{LAST}),(\ref{e.BAFI}.\ref{ETE})
and use the fact that $\lambda_C$ is sub-multiplicative. Finally,
applying \eqref{LACT} and (\ref{ss.IDCON}.\ref{CJ}), we obtain
\[
\lambda_C(x - x^\ast) = \lambda_C\big(2x - t_C(x)1_C\big) \geq
\mbm\,\{e_F + \lambda_C(x), \wid(C) + \lambda_C(x)\},
\]
and \eqref{LACAS} follows from (\ref{e.BAFI}.\ref{ETE}). \hfill
$\square$

\subsection{Lemma.} \label{l.PIXI} \emph{Suppose $B$ as in}
\ref{ss.LACON} \emph{is wild and $\mu \in \mfo^\times$ has the form}
\begin{align}
\label{MUDE} \mu = 1 - \pi^d\beta, \quad d \in \IZ, \quad 0 \leq d
\leq 2\tex(B), \quad \beta \in \mfo^\times.
\end{align}
\emph{Setting $C = \mbC(B,\mu)$ and}
\begin{align}
\label{EMTH} 0 \leq m := \floor{ \frac{d}{2}} < \tex(B), \quad
\Theta_d := \pi^{-m}(1_B + j) \in C
\end{align}
\emph{as in} (\ref{ss.EARE}.\ref{REDSCA}), \emph{the relations}
\begin{align}
%\label{ENTH} n_C(v\Theta) =\,\,&\pi^{d-2m}\beta n_B(v), \quad
%t_C(v\Theta) =
%\pi^{-m}t_B(v), \quad C = B \oplus B\Theta, \\
\label{UVTH} u(v\Theta_d) =\,\,&\pi^{-m}[u,v] + (vu)\Theta_d, \\
\label{VTHU} (v\Theta_d)u =\,\,&\pi^{-m}v(u - u^\ast) +
(vu^\ast)\Theta_d, \\
\label{VTHVTH} (v_1\Theta_d)(v_2\Theta_d)
=\,\,&\Big(\pi^{-2m}[v_1,v_2 - v_2^\ast] - \pi^{d-2m}\beta
v_2^\ast v_1\Big) + \pi^{-m}\Big(t_B(v_2)v_1 -
[v_1,v_2]\Big)\Theta_d
\end{align}
\emph{hold for all $u,v,v_1,v_2 \in B$. Moreover, if $B$ is
flexible, then $\mfO := \mfo_B \oplus \mfo_B\Theta_d$ is an
$\mfo$-subalgebra of $C$.}

\Proof A slightly involved but straightforward computation using
the transition formulas
\[
u + v\Theta_d = (u + \pi^{-m}v) + \pi^{-m}vj, \quad u + vj = (u -
v) + \pi^mv\Theta_d
\]
implies \eqref{UVTH}$\--$\eqref{VTHVTH}. Combining these with
Prop.~\ref{p.ELAC} leads to the final assertion of the lemma.
\hfill $\square$

\subsection{Theorem.} \label{t.CONUNOD} \emph{Suppose $B$ as in}
\ref{ss.LACON} \emph{is flexible and wild},
\[
\mu = 1 - \pi^d\beta, \quad d \in \IZ, \quad 0 \leq d <2\tex(B),
\quad \beta \in \mfo^\times,
\]
\emph{and $d$ is} odd. \emph{Then $C := \mbC(B,\mu)$ is a flexible
$\lambda$-normed conic algebra over $F$ with $\bar C = \bar B$ as
conic $\bar F$-algebras. Moreover, $C$ is $\lambda$-valued if and
only if $B$ is $\lambda$-valued.}

\Proof The proof is similar to, but a bit more complicated than,
the one of Prop.~\ref{p.CONPR}. By Prop.~\ref{p.CDUNOD}, $C$ is
non-singular, round and anisotropic. Combining the final statement
of Lemma~\ref{l.PIXI} with (\ref{p.CDUNOD}.\ref{OCEPROD}), we
conclude that $\mfo_C \subseteq C$ is an $\mfo$-subalgebra which,
thanks to (\ref{l.PIXI}.\ref{UVTH}$\--$\ref{VTHVTH}) and to $d$
being odd contains $\mfp_C$ as an ideal with $\mfp_C^2 \subseteq
\mfp\mfo_C$. Thus $C$ is $\lambda$-normed (Prop.~\ref{p.VARIID}),
and it remains to show that if $B$ is $\lambda$-valued, so is $C$.
Let $x_i = u_i + v_i\Pi \in C^\times$, $u_i,v_i \in B$, $i = 1,2$
and $x_1x_2 = u + v\Pi$, where
(\ref{l.PIXI}.\ref{UVTH}$\--$\ref{VTHVTH}) imply
\begin{align}
\label{U} u =\,\,&u_1u_2 + \pi^{-m}\big([u_1,v_2] + v_1(u_2 -
u_2^\ast)\big) + \pi^{-2m}[v_1,v_2 - v_2^\ast] -\pi\beta v_2^\ast
v_1, \\
\label{V} v =\,\,&v_2u_1 + v_1u_2^\ast + \pi^{-m}\big(t_B(v_2)v_1
- [v_1,v_2]\big).
\end{align}
We must show $\lambda_C(x_1x_2) = \lambda_C(x_1) + \lambda_C(x_2)$.
To this end, invoking Lemma~\ref{l.COVA}, we may assume $0 \leq
\lambda_C(x_i) \leq \frac{1}{2}$, $i = 1,2$. Since conjugation is an
algebra involution of $C$ leaving $\lambda_C$ invariant, there are
three cases: (i) $\lambda_C(x_1) = \lambda_C(x_2) = 0$, (ii)
$\lambda_C(x_1) = 0$, $\lambda_C(x_2) = \frac{1}{2}$, (iii)
$\lambda_C(x_1) = \lambda_C(x_2) = \frac{1}{2}$. Among these cases,
we treat only (iii) since the other ones can be treated analogously.
In (iii) we have $u_1,u_2 \in \mfp_B$, $v_1,v_2 \in \mfo_B^\times$,
observe Prop.~\ref{p.ELAC} and obtain $u \equiv -\pi\beta v_2^\ast
v_1 \bmod \mfp_B^2$ by \eqref{U}, hence $\lambda_B(u) = 1$, while
\eqref{V} yields $\lambda_B(v) \geq 1$. Therefore $\lambda_C(x_1x_2)
= 1 = \lambda_C(x_1) + \lambda_C(x_2)$. \hfill $\square$

\subsection{Theorem.} \label{t.CONUNEV} \emph{Suppose $B$ as in}
\ref{ss.LACON} \emph{is flexible and wild},
\begin{align}
\label{MUDEV} \mu = 1 - \pi^d\beta, \quad d \in \IZ, \quad 0 \leq d
<2\tex(B), \quad \beta \in \mfo, \quad \bar\beta \notin V_{\bar B}^2
\end{align}
\emph{and $d$ is} even. \emph{Then $C:= \mbC(B,\mu)$ is a wild
$\lambda$-normed conic algebra over $F$ with $\bar C \cong
\mbC(\bar B,\bar\beta)$ as conic $\bar F$-algebras. Moreover, $C$
is $\lambda$-valued if and only if $\bar C$ is a division
algebra.}

\Proof By Prop.~\ref{p.CDUNEV}, $C$ is a non-singular, round and
anisotropic conic $F$-algebra. Moreover, $C$ is wild of
ramification index $1$. The final statement of Lemma~\ref{l.PIXI}
combined with (\ref{p.CDUNEV}.\ref{OCEPREV}) shows that $\mfo_C
\subseteq C$ is an $\mfo$-subalgebra, forcing $C$ to be
$\lambda$-normed (Cor.~\ref{c.LAVADIV}~(a)). Moreover, writing
$\mbC(\bar B,\bar\beta) = \bar B \oplus \bar Bj^\prime$,
$j^{\prime 2} = \bar\beta1_{\bar B}$ as in \ref{ss.CDCO},
consulting (\ref{l.PIXI}.\ref{UVTH}$\--$\ref{VTHVTH}) and
observing $d < \tex(B) = \wid(B)$, Prop.~\ref{p.ELAC} shows that
(in the notations of Prop.~\ref{p.CDUNEV}) the assignment $u +
v\Xi \mapsto \bar u + \bar vj^\prime$ determines an isomorphism
$\bar C \overset{\sim} \rightarrow \mbC(\bar B,\bar\beta)$ of
$\bar F$-algebras. The final statement of the theorem follows
immediately from Cor.~\ref{c.LAVADIV}~(b). \hfill $\square$

\subsection{Corollary.} \label{c.LAPINS} \emph{With the notations
and assumptions of Thm.} \ref{t.CONUNEV}, \emph{suppose in
addition that $\bar B/\bar F$ is a purely inseparable field
extension of exponent at most $1$. Then $C$ is a $\lambda$-valued
conic algebra over $F$ with $\bar C \cong \bar B(\sqrt{\bar
\beta})$.}

\Proof It suffices to note that the last condition of
(\ref{t.CONUNEV}.\ref{MUDEV}) makes $\mbC(\bar B, \bar\beta) =
\bar B(\sqrt{\bar\beta})$ (cf. \ref{ss.CDP}, Case~2) a division
algebra. \hfill $\square$

\Remark The additional hypothesis in Cor.~\ref{c.LAPINS} is
fulfilled if, e.g., $B$ is a composition algebra (Remark to
\ref{p.VARIID}).

\subsection{Examples.} \label{e.LAU} In Brown's
examples of conic division algebras (cf.~\ref{e.LAPR}), one
basically keeps building up ramified $\lambda$-valued conic algebras
over iterated Laurent series fields of characteristic not $2$. By
contrast, we will now be able to construct wild $\lambda$-valued
conic algebras of ramification index $1$ over appropriate Henselian
fields of characteristic zero. Let $k$ be any field of
characteristic $2$ and write $K$ for the field of rational functions
in an infinite number of variables over $k$. Then $[K:K^2] =
\infty$. Pick an infinite chain
\[
K = K_0 \subset K_1 \subset \cdots \subset K_{n-1} \subset K_n
\subset \cdots
\]
of purely inseparable field extensions of $K$ having exponent at
most $1$ and $[K_n:K] = 2^n$ for all integers $n \geq 0$. Following
Teichm\"uller \cite{Te37}, there is an essentially unique complete
field $F$ under a discrete valuation $\lambda\:F \to \IZ_\infty$
such that $F$ has characteristic zero, residue field $\bar F = K$
and absolute ramification index $e_F = 1$. For $n \geq 1$ choose
$\beta_n \in \mfo^\times$ such that $K_n =
K_{n-1}(\sqrt{\bar\beta_n})$, put $\mu_n = 1 - \beta_n$, observe
(\ref{e.BAFI}.\ref{TEEF}) and apply Cor.~\ref{c.LAPINS} successively
for $d = 0$ and $n = 1,2,3,\dots$ to conclude that the
Cayley-Dickson process leads to a wild $\lambda$-valued conic
algebra
\[
C_n := \mbC(F;\mu_1,\dots,\mu_n)
\]
over $F$ having dimension $2^n$ and ramification index $1$ such
that $\bar C_n \cong K_n$.

\subsection{Corollary.} \label{c.PFICON} \emph{Let $Q$ be a
pointed Pfister quadratic space over $F$ that is anisotropic and
wild of ramification index $e_{Q/F} = 1$. Then there exists a
flexible $\lambda$-valued conic algebra $C$ over $F$ such that
$Q_C \cong Q$ and $\bar C/\bar F$ is a purely inseparable field
extension of exponent at most $1$.}

\Proof Arguing by induction, we let $Q$ be a pointed
$(n+1)$-Pfister quadratic space and pick a pointed $n$-Pfister
quadratic subspace $P \subseteq Q$. Clearly, $P$ is wild with
$e_{P/F} = 1$. By Theorem~\ref{t.ENBARE}, some scalar $\mu$ of
standard type relative to $P$ satisfies $Q = \dla\mu\dra \otimes
P$ up to isomorphism, and since the implications (a),(b) of that
theorem do not hold for $\mu$, implication (c) does. On the other
hand, the induction hypothesis leads to a flexible
$\lambda$-valued conic algebra $B$ over $F$ with $Q_B \cong P$
such that $\bar B/\bar F$ is a purely inseparable field extension
of exponent at most $1$. By Cor.~\ref{c.LAPINS}, $C :=
\mbC(B,\mu)$ is a flexible $\lambda$-valued conic algebra over $F$
with $Q_C = \dla\mu\dra \otimes P = Q$ and $\bar C \cong \bar
B(\sqrt{\bar\beta})$. \hfill $\square$ \lz

\noindent We still haven't closed the gap in our proof of
Thm.~\ref{t.CHAGORE} but will now be able to do so by appealing to
the connection with conic algebras. In view of
Cors.~\ref{c.LAPINS},\ref{c.PFICON}, the missing equivalence of
(v) and (i)$\--$(iv) in Thm.~\ref{t.CHAGORE} will be a consequence
of the following result.

\subsection{Theorem.} \label{t.CONCHA} \emph{Suppose $B$ as in}
\ref{ss.LACON} \emph{is flexible, wild and $\lambda$-normed having
$\bar B/\bar F$ as a purely inseparable field extension of exponent
at most $1$ and, with the notations of} \ref{ss.SETEX}~(b),
\emph{let}
\begin{align}
\label{BETEX} \beta \in \mfo,\quad \mu := 1 - \pi^{2\tex(B)}\beta
\in \mfo^\times, \quad \beta_0 := n_B(w_0)\beta.
\end{align}
\emph{Then $C := \mbC(B,\mu)$ is anisotropic if and only if
$\overline{\beta_0} \notin \Im(\wp_{\bar B,\bar s_{w_0}})$. In
this case, setting}
\begin{align}
\label{XINUL} \Xi_0 := \pi^{-\tex(B)}(w_0 + w_0j) \in C,
\end{align}
\emph{we obtain the relations}
\begin{align}
\label{NTXINUL} t_C(\Xi_0) =\,\,&1, \quad n_C(\Xi_0) = \beta_0,
\quad C = B \oplus B\Xi_0, \\
\label{VAXINUL} \lambda_C(u\,+\,&v\Xi_0) =
\mbm\,\{\lambda_B(u),\lambda_B(v)\} && (u,v \in B), \\
\label{VARIXINUL} \mfo_C = \mfo_B&\oplus\mfo_B\Xi_0, \quad \mfp_C
= \mfp_B \oplus \mfp_B\Xi_0
\end{align}
\emph{and $C$ is $\lambda$-normed with $Q_{\bar C} \cong Q_{\bar
B;\bar\beta_0,\bar s_{w_0}}$.}

\Proof As in the proof of the implication (iii) $\Rightarrow$ (i) in
Thm.~\ref{t.CHAGORE}, we put $d := 2\tex(B)$, $\Xi := \Theta_d$ and
have $\Xi_0 = w_0\Xi$. Thus the first two relations of
\eqref{NTXINUL} follow from
(\ref{ss.EARE}.\ref{NODECTH},\ref{TRDECTH}), while the last one is a
straightforward consequence of the fact that $B$ is $\lambda$-valued
by Cor.~\ref{c.LAVADIV}~(b), hence a division algebra
(Prop.~\ref{p.ELCOMCOM}~(a)). Setting $m := \tex(B)$, we let $v \in
B$ and compute, using
(\ref{ss.CDCO}.\ref{CDMU},\ref{CDNO}),(\ref{p.NOCOMT}.\ref{NOCO})
and \eqref{BETEX},
\begin{align*}
n_C(v\Xi_0) =\,\,&\pi^{-2m}n_C\big(v(w_0 + w_0j)\big) =
\pi^{-2m}n_C\big(vw_0 + (w_0v)j\big) \\
=\,\,&\pi^{-2m}\big(n_B(vw_0) - \mu n_B(w_0v)\big) =
\pi^{-2m}n_B(vw_0)(1 - \mu) = n_B(vw_0)\beta.
\end{align*}
But since $\bar B/\bar F$ is a purely inseparable field extension
of exponent at most $1$, we have $n_B(vw_0) \equiv n_B(v)n_B(w_0)
\bmod \mfp$ and conclude
\begin{align}
\label{MULXI} n_C(v\Xi_0) \equiv n_B(v)\beta_0 \bmod \mfp.
\end{align}
Suppose first that $C$ is anisotropic. Then $\bar C$ is an
anisotropic conic algebra over $\bar F$ containing $\bar B$ as a
subalgebra and $l := \overline{\Xi_0} \in \bar C$ as a distinguished
element with $t_{\bar C}(l) = 1$, $n_{\bar C}(l) =
\overline{\beta_0}$, $n_{\bar C}(\bar u,l) = \bar s_{w_0}(\bar u)$
for all $u \in \mfo_B$. Moreover, \eqref{MULXI} implies $n_{\bar
C}(\bar vl) = \overline{\beta_0}\bar v^2$ for all $v \in \mfo_B$.
Writing $\bar B \oplus \bar Bj^\prime$ for the vector space
underlying the pointed quadratic space $Q_{\bar B;\bar\beta_0,\bar
s_{w_0}}$, the assignment $\bar u + \bar vj^\prime \mapsto \bar u +
\bar vl$ therefore and by (\ref{ss.FLECON}.\ref{NAS}) gives an
embedding $\varphi$ from $Q_{\bar B;\bar\beta_0,\bar s_{w_0}}$ to
$Q_{\bar C}$ of pointed quadratic spaces. Comparing dimensions
(observe $e_{C/F} = 1$ by (i) of Thm.~\ref{t.CHAGORE}),
$\varphi\:Q_{\bar B;\bar\beta_0,\bar s_{w_0}} \overset{\sim} \to
Q_{\bar C}$ is, in fact, an isomorphism, and since $Q_{\bar C}$ is
anisotropic, so is $Q_{\bar B;\bar\beta_0,\bar s_{w_0}}$. Now
$\overline{\beta_0} \notin \Im(\wp_{\bar B,\bar s_{w_0}})$ follows
from Cor.~\ref{c.ISAS}~(a). Moreover, we claim that \eqref{VAXINUL}
holds (which immediately implies \eqref{VARIXINUL}). As usual, we
may assume $u,v \in \mfo_B^\times$, which yields $\lambda_C(u +
v\Xi_0) \geq 0$, and if this were strictly positive, we would end up
with $\varphi(\bar u + \bar vj^\prime) = \bar u + \bar vl = 0$,
forcing $\bar u = \bar v = 0$, a contradiction. Suppose next
$\overline{\beta_0} \notin \Im(\wp_{\bar B,\bar s_{w_0}})$ and
consider the full $\mfo$-lattice
\begin{align}
\label{FULLA} \mfO := \mfo_B \oplus \mfo_B\Xi_0 \subseteq C,
\end{align}
on which $n_C$ takes integral values. More precisely,
\eqref{MULXI} and (\ref{ss.FLECON}.\ref{NAS}) imply
\begin{align}
\label{NORO} n_C(u + v\Xi_0) \equiv n_B(u) + \pi^{-m}n_B(v^\ast
u,\Xi_0) + n_B(v)\beta_0 \bmod \mfp
\end{align}
for all $u,v \in \mfo_B$. Reducing $\bmod\,\mfp$, we obtain a
pointed quadratic space
\[
\bar\mfO := \mfO \otimes_\mfo \bar F = \bar B \oplus \bar B
l^\prime, \quad l^\prime := \Xi_0 \otimes_\mfo 1_{\bar F},
\]
over $\bar F$, and the assignment $\bar u \oplus \bar vj^\prime
\mapsto \bar u + \bar vl^\prime$ by \eqref{FULLA},\eqref{NORO} gives
an isomorphism from $Q_{\bar B;\bar\beta_0,\bar s_{w_0}}$ onto $\bar
\mfO$. The former being anisotropic by Cor.~\ref{c.ISAS}~(a), so is
the latter. But then $C$ must be anisotropic as well
since every non-zero element $x \in C$ satisfies $\pi^mx \in \mfO
\setminus\mfp\mfO$ for some integer $m$.

It remains to show that $C$ is $\lambda$-normed provided it is
anisotropic. By Cor.~\ref{c.LAVADIV}~(a), it suffices to show that
$\mfo_C \subseteq C$ is an $\mfo$-subalgebra. In order to do so, we
use Lemma~\ref{l.PIXI} to derive the following formulas by a
straightforward  computation, for all $u,v,v_1,v_2 \in B$.
\begin{align}
\label{VEXINUL} v\Xi_0 =\,\,&\pi^{-m}[v,w_0] + (w_0v)\Xi, \\
\label{VEXIB} v\Xi =\,\,&\pi^{-m}[w_0,L_{w_0}^{-1}v] +
(L_{w_0}^{-1}v)\Xi_0, \\
\label{UVEXINUL} u(v\Xi_0) =\,\,&\pi^{-m}\big(u[v,w_0] +
[u,w_0v]\big) + \big((w_0v)u\big)\Xi, \\
\label{VEXINULU} (v\Xi_0)u =\,\,&\pi^{-m}\big([v,w_0]u + (w_0v)(u
- u^\ast)\big) + \big((w_0v)u^\ast\big)\Xi \\
\label{VEXINULVEXINUL} (v_1\Xi_0)(v_2\Xi_0)
=\,\,&\pi^{-2m}\Big([v_1,w_0][v_2,w_0] + [[v_1,w_0],w_0v_2] + \\
\,\,&[w_0v_1,w_0v_2 - (w_0v_2)^\ast] + (w_0v_1)\big([v_2,w_0] -
[v_2,w_0]^\ast\big)\Big) - \notag \\
\,\,&\beta(w_0v_2^\ast)(w_0v_1) + \pi^{-m}\Big((w_0v_2)[v_1,w_0] +
(w_0v_1)[v_2,w_0]^\ast + \notag \\
\,\,&t_B(w_0v_2)w_0v_1 - [w_0v_1,w_0v_2]\Big)\Xi. \notag
\end{align}
Since $B$ is $\lambda$-valued, \eqref{VEXIB} implies $\mfo_B\Xi
\subseteq \mfo_C$, and then
\eqref{UVEXINUL}$\--$\eqref{VEXINULVEXINUL} combine with
Prop.~\ref{p.ELAC} to establish $\mfo_C$ as an $\mfo$-subalgebra of
$C$. \hfill $\square$

\subsection{Corollary.} \label{c.CONTAM} \emph{Suppose in
Thm.}~\ref{t.CONCHA} \emph{that $B$ is an associative composition
division algebra and $\overline{\beta_0} \notin \Im(\wp_{\bar
B;\bar\beta_0,\bar s_{w_0}})$. Then $C$ is an unramified
composition division algebra over $F$ with $\bar C \cong \mbC(\bar
B;\overline{\beta_0},\bar s_{w_0})$ as a non-orthogonal
Cayley-Dickson construction.}

\Proof Composition algebras are classified by their norms. \hfill
$\square$

\Remark If $C$ as in Thm.~\ref{t.CONCHA} is anisotropic, its
pointed quadratic residue space is described explicitly by the
theorem. But $C$ is also a $\lambda$-normed conic algebra, making
$\bar C$ canonically a conic algebra (over $\bar F$) in its own
right. It would be interesting to obtain an equally explicit
description of that algebra. Cor.~\ref{c.CONTAM} provides one if
$C$ is a composition algebra but it is not at all clear whether
this description prevails in the general case, nor whether $C$ is
always a $\lambda$-\emph{valued} conic algebra.

%%%%%%%%%%%%%%%%%%%%%%%%%%%%%%%%%%%%%%%%%%%%%%%%%%
\section{Applications to composition algebras.} \label{s.APCOM}
There are obvious and less obvious applications of the preceding
results to composition algebras. Working over a fixed
$2$-Henselian field $F$ of arbitrary characteristic as before (cf.
\ref{ss.HEFI}), the obvious ones may be described as follows.

\subsection{Translations.} \label{ss.TRANS} Since composition
division algebras over $F$ are classified by their norms and are
$\lambda$-valued conic $F$-algebras by Prop.~\ref{p.ELCOMCOM}~(b),
the results of Sections $8,9$ translate immediately into this more
special setting, where the ones in Section $9$ in particular yield
explicit descriptions of how the valuation data ramification
index, residue algebra and trace exponent behave under the
Cayley-Dickson construction. Rather than carrying out these
translations in full detail, suffice it to point out that all one
has to do is replace
\begin{itemize}
\item the pointed quadratic space $P$ of
\ref{c.ODNE},\ref{ss.GEST},\ref{ss.EARE} by an associative
composition division algebra $B$ over $F$ having ramification
index $e_{B/F} = 1$,

\item the pointed quadratic space $Q$ by a composition algebra
$C$, and the condition of $Q$ being anisotropic by the one of $C$
being a division algebra,

\item the passage from $P$ to $\dla\mu\dra \otimes P$, $\mu \in
F^\times$, by the Cayley-Dickson construction $\mbC(B,\mu)$.
\end{itemize}
The less obvious applications of our results to composition
algebras are all related, in one way or another, to the following
innocuous observation, which we have not been able to extend to
pointed Pfister quadratic spaces.

\subsection{Proposition.} \label{p.LIFT} \emph{Let $C$ be a
composition division algebra over $F$ and $B^\prime \subseteq \bar
C$ a unital subalgebra. Assume $\mbc(F) \neq 2$ or $\dim_{\bar
F}(B^\prime) > 1$. Then there exists a composition subalgebra $B
\subseteq C$ having ramification index $e _{B/F} = 1$ and
satisfying $\bar B = B^\prime$.}

\Proof One adapts the proof of \cite[Lemma~3]{MR51:635} to the
present more general set-up; for completeness, we include the
details. Since $B^\prime$ is either a composition division algebra
over $\bar F$ or a purely inseparable field extension of
characteristic $2$ and exponent at most $1$, it has dimension
$2^m$, $m \in \IZ$, $0 \leq m \leq 3$. Moreover $B^\prime$ is
generated by $m$ elements $\overline{x_1},\dots,\overline{x_m}$,
for some $x_1,\dots,x_m \in \mfo_C$, where we may assume $m \geq
1$ since $m = 0$ implies $\mbc(F) \neq 2$ by hypothesis and $B :=
F$ does the job. The elements of non-zero trace in $C$ form a
Zariski open and dense subset, which therefore is open and dense
in the valuation topology as well, so we may assume $t_C(x_1) \neq
0$. Then $B$, the unital subalgebra of $C$ generated by
$x_1,\dots,x_m$, is a composition division algebra with $\mbd_F(B)
\leq 2^m$, $B^\prime \subseteq \bar B$. Now Prop.~\ref{p.EFEN}
implies $\bar B = B^\prime$ and $e_{B/F} = 1$. \hfill $\square$
\lz

\noindent The preceding result can be refined in various ways. For
example, given a composition division algebra $C$ over $F$, we will
exhibit (chains of) proper composition subalgebras of $C$ having
ramification index $1$ and the same trace exponent as $C$. From this
we derive normal forms for octonion and quaternion algebras over $F$
and show that quantities subject to a few obvious constraints are
the valuation data of an appropriate composition division algebra.
We begin by listing a few properties of wild separable quadratic
field extensions which should be well known but seem to lack a
convenient reference. We therefore include the details.

\subsection{Proposition.} \label{p.SEQUA} \emph{Let $L$ be an
$F$-algebra and suppose $\bar F$ has characteristic $2$. Then the
following conditions are equivalent.}
\begin{itemize}
\item [(i)] \emph{$L/F$ is a wild separable quadratic field
extension and $e_{L/F} = 1$.}

\item[(ii)] \emph{There are a positive integer $r$ and elements
$\alpha \in \mfo^\times$, $\beta \in \mfo$ such that $\bar \beta
\notin \bar F^2$ and}
\[
L \cong F[\bft]/(\bft^2 - \pi^r \alpha\bft + \beta).
\]
\end{itemize}
\emph{If these conditions hold, $\bar L = \bar F(\sqrt{\bar
\beta})$. Moreover, setting}
\begin{align}
\label{GEN} \vartheta := \bft \bmod(\bft^2 - \pi^r \alpha\bft +
\beta) \in L
\end{align}
\emph{in} (ii), \emph{the following relations hold.}
\begin{align}
\label{VALL} \lambda_L(\gamma + \delta\vartheta)
=\,\,&\mathrm{min}\,\{\lambda(\gamma),\lambda(\delta)\}
&&(\gamma,\delta \in F), \\
\label{OL} \mfo_L = \mfo 1_L \oplus \mfo \vartheta, &\quad
\mfp_L = \mfp 1_L \oplus \mfp\vartheta, \\
\label{TRL} \tex(L) =\,\,&\mathrm{min}\,\{e_F,r\}.
\end{align}

\Proof (i) $\Longrightarrow$ (ii). By (i) there exists an element
$\beta \in \mfo$ such that $\bar \beta \notin \bar F^2$ and $\bar
L = \bar F(\sqrt{\bar \beta})$. Pick an element $\vartheta \in
\mfo_L$ satisfying $\bar \vartheta = \sqrt{\bar \beta}$. Then
$t_L(\vartheta) \in \mfp$, and replacing $\beta$ by
$n_L(\vartheta)$ if necessary, we may assume $n_L(\vartheta) =
\beta$. Here $\bar \beta \notin \bar F^2$ forces $L = F(\vartheta)
= F[\vartheta]$. We claim there is no harm in assuming
$t_L(\vartheta) \neq 0$. Indeed, for $\mathrm{char}(F) = 2$, this
is automatic while, if $\mbc(F) = 0$, we may replace $\vartheta$
by $\vartheta^\prime := 1_L + \vartheta$. Thus, without loss,
$t_L(\vartheta) \neq 0$. But this yields a unit $\alpha \in
\mfo^\times$ such that $t_L(\vartheta) = \pi^r \alpha$, $r :=
\lambda(t_L(\vartheta)) \in \IZ$, $r > 0$, and $L$ has the form
described in (ii).

(ii) $\Longrightarrow$ (i). By the hypotheses on $\beta$, the
monic polynomial $f := \bft^2 - \pi^r\alpha\bft + \beta \in
\mfo[\bft] \subseteq F[\bft]$ is irreducible over $F$, and $L/F$
is a separable quadratic field extension. Define $\vartheta$ as in
\eqref{GEN}. We have $n_L(\vartheta) = \beta \in \mfo^\times$,
hence $\vartheta \in \mfo_L^\times$, and ${\bar \vartheta}^2 =
\bar \beta$, which implies $\bar L = \bar F(\sqrt{\bar \beta})$,
so $L$ is wild and has ramification index $e_{L/F} = 1$. Hence (i)
holds.

A standard argument now yields \eqref{VALL}, which immediately
implies the remaining assertions of the proposition. \hfill
$\square$

%\Remark If $F$ has characteristic $0$, Prop.~\ref{p.SEQUA} is just
%a special case of Prop.~\ref{p.CDUNEV} (combined wit the
%translation formalism of \ref{ss.TRANS}) for $B = F$.

\subsection{Corollary.} \label{c.TREL} \emph{For a separable
quadratic field extension $L/F$ to be wild and to have
ramification index $1$ it is necessary and sufficient that
$\tex(L) > 0$ and there exist a trace generator $u$ of $L$ with
$\overline{n_L(u)} \notin \bar F^2$. In this case, $L = k[u]$ and
$u$ may be so chosen as to satisfy the additional relation}
\begin{align}
\label{TREC} \lambda_L(u - u^\ast) = \tex(L).
\end{align}

\Proof \emph{Necessity and the final statement.} If $L$ is wild and
$e_{L/F} = 1$, we obtain $\tex(L) > 0$ and first deal with the case
$e_F > \tex(L)$. Then we find a positive integer $r$ and elements
$\alpha \in \mfo^\times$, $\beta \in \mfo$ as in
Prop.~\ref{p.SEQUA}~(ii) and conclude $\tex(L) = r$ from
(\ref{p.SEQUA}.\ref{TRL}). Moreover, $u :=\vartheta$ as defined in
(\ref{p.SEQUA}.\ref{GEN}) is a trace generator of $L$ satisfying
$\overline{n_L(u)} \notin \bar F^2$ and $L = F[u]$. Finally,
$\lambda_L(u - u^\ast) = \lambda_L(2u - \pi^r\alpha 1_L)$, where
$\lambda_L(2u) = e_F > r = \lambda_L(\pi^r\alpha 1_L)$, which
implies \eqref{TREC} as well. By (\ref{e.BAFI}.\ref{ETE}), we are
left with the case $e_F = \tex(L) < \infty$. Then $F$ has
characteristic $0$, allowing us to apply Cor.~\ref{c.EQTE}~(b) with
$B = P = F$: there exists a scalar $\mu \in \mfo$ such that $\bar
\mu \notin \bar F^2$ and $L = F(\sqrt{\mu})$, so some $y \in
\mfo_L^\times$ has
\begin{align}
\label{TREZE} L = F[y], \quad t_L(y) = 0, \quad \overline{n_L(y)}
\notin \bar F^2.
\end{align}
Hence $u := 1_L + y \in \mfo_L$ is a trace generator of $L$
satisfying $\overline{n_L(u)} \notin \bar F^2$. Moreover, since
$y^\ast = -y$ by \eqref{TREZE}, $\lambda_L(u - u^\ast) =
\lambda_L(2y) = e_F = \tex(L)$, and the proof is complete. $\ssp$ \\
\emph{Sufficiency.} We have $L = k[u]$, $t_L(u) = \pi^r\alpha$, $r
:= \tex(L)$, $\alpha \in \mfo^\times$, and condition (ii) of
Prop.~\ref{p.SEQUA} holds with $\beta = n_L(u)$. \hfill $\square$

%\subsection{Remark.} \label{r.TRELS} If $L/F$ is an unramified
%(separable) quadratic field extension, so $e_{L/F} = 1$ and $F$ is
%tame, then (\ref{c.TREL}.\ref{TREC}) can still be solved in
%$\mfo_L^\times$. Indeed, the conjugation of $\bar L$ as a conic
%$\bar F$-algebra being non-trivial, some $u \in \mfo_L^\times$ has
%$\bar u^\ast \neq \bar u$, i.e., $\lambda_L(u - u^\ast) = 0 =
%\tex(L)$.

\subsection{Proposition.} \label{p.SEQUAT} \emph{Let $L$ be an
$F$-algebra and suppose $\bar F$ has characteristic $2$. Then the
following conditions are equivalent.}
\begin{itemize}
\item [(i)] \emph{$L/F$ is a separable quadratic field extension
of ramification index $e_{L/F} = 2$.}

\item [(ii)] \emph{There are a positive integer $r$ and elements
$\alpha,\beta \in \mfo^\times$ such that}
\[
L \cong F[\bft]/(\bft^2 - \pi^r\alpha\bft + \pi \beta).
\]
\end{itemize}
\emph{In this case, for any prime element $\Pi \in \mfo_L$}
(\emph{e.g., for}
\begin{align}
\label{PRAM} \Pi := \bft \bmod (\bft^2 - \pi^r\alpha\bft +
\pi\beta))
\end{align}
\emph{the following relations hold.}
\begin{align}
\label{VALR} \lambda_L(\gamma + \delta\Pi)
=\,\,&\mbm\,\{\lambda(\gamma),\lambda(\delta) +
\frac{1}{2}\} &&(\gamma,\delta \in F), \\
\label{OLR} \mfo_L = \mfo \oplus \mfo\Pi&, \quad
\mfp_L = \mfp \oplus \mfo\Pi, \\
\label{TRLR} \tex(L) =\,\,&\mbm\,\{e_F,r\}.
\end{align}

\Proof Everything is standard once it has been shown in (i)
$\Rightarrow$ (ii) that $\mfo_L$ contains prime elements $\Pi$
with $t_L(\Pi) \neq 0$. But this follows from the fact that the
set of elements in $L$ with non-zero trace, by separability being
open and dense in the Zariski topology, is open and dense in the
valuation topology as well. \hfill $\square$ \lz

\subsection{Theorem.} \label{t.SLIFT} \emph{Let $C$ be a
composition division algebra of dimension $2^n$, $n = 2,3$, over
$F$. Then there exists a separable quadratic subfield $L \subseteq
C$ having ramification index $1$ and the same trace exponent as
$C$}: $e_{L/F} = 1$, $\tex(L) = \tex(C)$.

\Proof Setting $r := \tex(C)$, we proceed in four steps. $\ssp$ \\
$1^0$. Let us first consider the case $r = 0$. Then $\bar C$ is a
composition division algebra of dimension at least $2$ over $\bar
F$ and hence contains a separable quadratic subfield $L^\prime
\subseteq \bar C$, which by Prop.~\ref{p.LIFT} or
\cite[Lemma~3]{MR51:635}, lifts to a separable quadratic subfield
$L \subseteq C$ with $e_{L/F} = 1$, $\tex(L) = 0 = r$. We may
therefore assume from now on that $r > 0$, so $C$ is wild. $\ssp$ \\
$2^0$. Next we deal with the case $e_{C/F} = 1$. Pick a trace
generator $w_0$ of $C$, which belongs to $\mfo_C^\times$ by
(\ref{ss.TRAG}.\ref{ESTRAG}). If $w_0 \notin F1_C$, then $L
:=F[w_0]$ is a separable quadratic subfield of $C$ satisfying $1
\leq e_{L/F} \leq e_{C/F} = 1$, hence $e_{L/F} = 1$. From $\mfp^r
= \pi^r\mfo = t_L(\mfo w_0) \subseteq t_L(\mfo_L) =
\mfp^{\tex(L)}$ we conclude $\tex(L) \leq r$, which implies
$\tex(L) = r$ by Prop.~\ref{p.PROTE}~(b). On the other hand, if
$w_0 \in F1_C$, then $r = e_F$ and \emph{any} separable quadratic
subfield $L \subseteq C$ satisfies $e_{L/F} = 1$ as well as $r
\leq \tex(L) \leq e_F = r$ by Prop.~\ref{p.PROTE}~(b) and
(\ref{e.BAFI}.\ref{ETE}). $\ssp$ \\
$3^0$. We are left with the case $e_{C/F} = 2$. Then
Prop.~\ref{p.LIFT} yields a compositions subalgebra $B \subseteq
C$ with $e_{B/F} = 1$ and $\mathrm{dim}_F(B) = 2^{n-1}$. If
$\tex(B) = r$, we are done for $n = 2$ and my apply $2^0$ to $B$
for $n = 3$ to arrive at the desired conclusion. By
Prop.~\ref{p.PROTE}~(b), we may therefore assume $\tex(B) > r$.
But then Thm.~\ref{t.ENBARE} justifies the assumption $C =
\mbC(B,\mu)$, $\mu = 1 - \pi^d\beta$, for some \emph{odd} integer
$d$, $0 \leq d  < 2\tex(B)$, and some unit $\beta \in
\mfo^\times$. In other words, we are in the
situation of Prop.~\ref{p.CDUNOD}. $\ssp$ \\
$4^0.$ Applying $2^0$ to $B$ if $B$ is a quaternion algebra, we
find a separable quadratic subfield $L \subseteq B$ satisfying
$e_{L/F} = 1$, $\tex(L) = \tex(B)$. From Cor.~\ref{c.TREL} we
therefore conclude that there is an element $u \in \mfo_L^\times$
such that $t_L(u) = \pi^{\tex(B)}$, $\overline{n_L(u)} \notin \bar
F^2$. Now consider the element $w := u + u\Pi$, which is a unit of
$\mfo_C$ by (\ref{p.CDUNOD}.\ref{OCEPROD}). Moreover, setting $m
:= \frac{d-1}{2}$ and observing Prop.~\ref{p.CDUNOD},
\begin{align*}
t_C(w) =\,\,&\pi^r\varepsilon, \quad \varepsilon := 1 +
\pi^{\tex(B)-r} \in \mfo^\times, \\
n_C(w) =\,\,&n_L(u) + 2\pi^{-m}n_L(u) + \pi \beta n_L(u),
\end{align*}
hence $\overline{n_C(w)} = \overline{n_L(u)} \notin \bar F^2$, so
$w$ by Prop.~\ref{p.SEQUA} generates a separable quadratic
subfield $L^\prime \subseteq C$ with $e_{L^\prime/F} = 1$,
$\tex(L^\prime) = r$. \hfill $\square$

\subsection{Corollary.} \label{c.TRGEUN} \emph{Every composition
division algebra $C$ of dimension $2^n$, $n = 2,3$, over $F$
contains a trace generator which is a unit in $\mfo_C$.}

\Proof Picking $L \subseteq C$ as in Thm.~\ref{t.SLIFT}, every
trace generator of $L$ is one of $C$ and by
(\ref{ss.TRAG}.\ref{ESTRAG}) belongs to $\mfo_L^\times \subseteq
\mfo_C^\times$. \hfill $\square$

\subsection{Corollary.} \label{c.TEXPCONJ} \emph{Let $C$ be a
composition division algebra over $F$ with $f_{C/F} > 1$. Then}
\[
\tex(C) = \mbm\,\{\lambda_C(u - u^\ast) \mid u \in \mfo_C\} =
\mbm\,\{\lambda_C(u - u^\ast) \mid u \in \mfo_C^\times\}.
\]

\Proof For $u \in \mfo_C$ we apply
(\ref{ss.IDCON}.\ref{CJ}),(\ref{ss.HEFI}.\ref{EEF}),(\ref{e.BAFI}.\ref{ETE})
and obtain
\[
\lambda_C(u - u^\ast) = \lambda_C\big(2u - t_C(u)1_C\big) \geq
\mbm\,\{e_F + \lambda_C(u), \lambda\big(t_C(u)\big)\} \geq
\tex(C).
\]
Hence it suffices to show that there is an element $u \in
\mfo_C^\times$ with $\lambda_C(u - u^\ast) = \tex(C)$. If $C$ is
tame, the conjugation of $\bar C$ cannot be the identity, so
$\lambda_C(u - u^\ast) = 0 = \tex(C)$ for some $u \in
\mfo_C^\times$. We may therefore assume that $C$ is wild. Combining
the hypothesis $f_{C/F} > 1$ with Thm.~\ref{t.SLIFT}, we find a
separable quadratic subfield $L \subseteq C$ with $e_{L/F} = 1$ and
$\tex(L) = \tex(C)$. Now Cor.~\ref{c.TREL} yields and element $u \in
\mfo_L^\times \subseteq \mfo_C^\times$ such that $\lambda_C(u -
u^\ast) = \lambda_L(u - u^\ast) = \tex(L) = \tex(C)$. \hfill
$\square$ \lz

\noindent  We will see in Example~\ref{e.SEPQUA} below that the
hypothesis $f_{C/F} > 1$ in the preceding corollary cannot be
avoided.

\subsection{Corollary.} \label{c.COSUPRE} \emph{Let $C$ be a
composition division algebra of dimension $2^n$, $n = 2,3$, over
$F$. Then there exists a composition subalgebra $B \subseteq C$
with}
\[
\mathrm{dim}_F(B) = 2^{n-1}, \quad e_{B/F} = 1, \quad \tex(B) =
\tex(C).
\]

\Proof For $n = 2$, this is just Thm.~\ref{t.SLIFT}, so we may
assume $n = 3$, i.e., that $C$ is an octonion algebra. By
\cite[Lemma~3]{MR51:635}, we may also assume $\tex(C) > 0$ and,
applying Thm.~\ref{t.SLIFT} again, we find a separable quadratic
subfield $L \subseteq C$ with $e_{L/F} = 1$, $\tex(L) = \tex(C)$.
Hence it suffices to prove the following lemma.

\subsection{Lemma.} \label{l.COSUPRE} \emph{Let $C$ be a wild
octonion division algebra over $F$ and $L \subseteq C$ a separable
quadratic subfield such that $e_{L/F} = 1$, $\tex(L) = \tex(C)$.
Then there exists a quaternion subalgebra $L \subseteq B \subseteq
C$ with $e_{B/F} = 1$, $\tex(B) = \tex(C)$.}

\Proof Since $C$ is wild, its residue algebra $\bar C$ is a purely
inseparable field extension of exponent $1$ and degree at least $4$
over $\bar F$ containing $\bar L$ as a quadratic subfield. Pick an
element $y \in \mfo_C^\times$ satisfying $\bar y \in \bar C
\setminus \bar L$. Then $L$ and $y$ generate a composition
subalgebra $B \subseteq C$ of dimension at most $4$ whose residue
algebra contains $\bar L(\bar y)$. Hence $B \subseteq C$ is a
quaternion subalgebra containing $L$ and having $e_{B/F} = 1$,
$\tex(C) = \tex(L) \geq \tex(B) \geq \tex(C)$ by
Prop.~\ref{p.PROTE}~(b). \hfill $\square$ \lz

\subsection{Normal Form Theorem: octonion algebras.} \label{t.NOFOO}\index{normal form theorem!octonion algebras}
\emph{For $0 < r \leq e_F$ and an $F$-algebra $C$, the following
conditions are equivalent.}
\begin{itemize}
\item [(i)] \emph{$C$ is an octonion division algebra over $F$
having trace exponent $\tex(C) = r$.}

\item [(ii)] \emph{There exist}
\begin{itemize}
\item [(I)] \emph{a separable quadratic field extension $L/F$ such
that $e_{L/F} = 1$ and $\tex(L) = r$} (\emph{in particular, $L$ is
wild}),

\item [(II)] \emph{a scalar $\alpha \in \mfo$ such that $\bar
\alpha \notin \bar L^2$,}

\item [(III)] \emph{either a scalar $\beta \in \mfo$ such that
$\bar \beta \notin \bar L(\sqrt{\bar \alpha})^2$ and}
\begin{align}
%\label{ONR}
C \cong \mbC(L;\alpha,\beta) \notag
\end{align}
\emph{or a scalar $\beta \in \mfo^\times$ such that}
\begin{align}
%\label{OR}
C \cong \mbC(L;\alpha,1 - \pi \beta) \quad \text{\emph{or}} \quad
C \cong \mbC(L;\alpha,\pi \beta). \notag
\end{align}
\end{itemize}
\end{itemize}
\emph{If this is so, then $e_{C/F} = 1$, $\bar C \cong \bar
L(\sqrt{\bar \alpha},\sqrt{\bar \beta})$ in the first alternative
of} (III), \emph{while $e_{C/F} = 2$, $\bar C \cong \bar
L(\sqrt{\bar \alpha})$ otherwise.}

\Proof (i) $\Longrightarrow$ (ii). Thm.~\ref{t.SLIFT} yields a
separable quadratic subfield $L \subseteq C$ satisfying (I), and
applying Lemma~\ref{l.COSUPRE}, we find a quaternion subalgebra $L
\subseteq B \subseteq C$ with $e_{B/F} = 1$, $\tex(B) = \tex(C) =
r = \tex(L)$. Now Cor.~\ref{c.EQTE}~(b) yields a quantity $\alpha$
satisfying (II) and $B \cong \mbC(L,\alpha)$, $\bar B \cong \bar
L(\sqrt{\bar \alpha})$, while (a) and (b) of the same corollary
yield a quantity $\beta$ satisfying (III) as well as the remaining
assertions of the theorem.

(ii) $\Longrightarrow$ (i). This follows immediately from
Prop.~\ref{p.CDPR} and Cor.~\ref{c.EQTE}. \hfill $\square$ \lz

\noindent Basically the same arguments, inserting the obvious
simplifications at the appropriate places, leads to the following
quaternionic version of the preceding result.

\subsection{Normal Form Theorem: quaternion algebras.}\index{normal form theorem!quaternion algebras}
\label{t.NOFOQ} \emph{For $0 < r \leq e_F$ and an $F$-algebra $B$,
the following conditions are equivalent.}
\begin{itemize}
\item [(i)] \emph{$B$ is a quaternion division algebra over $F$
having trace exponent $\tex(B) = r$.}

\item [(ii)] \emph{There exist}
\begin{itemize}
\item [(I)] \emph{a separable quadratic field extension $L/F$ such
that $e_{L/F} = 1$ and $\tex(L) = r$} (\emph{in particular, $L$ is
wild}),

\item [(II)]  \emph{either a scalar $\alpha \in \mfo$ such that
$\bar \alpha \notin \bar L^2$ and}
\begin{align}
%\label{QNR}
B \cong \mbC(L,\alpha), \notag
\end{align}
\emph{or a scalar $\alpha \in \mfo^\times$ such that}
\begin{align}
%\label{QR}
B \cong \mbC(L;1 - \pi \alpha) \quad \text{\emph{or}} \quad B
\cong \mbC(L,\pi \alpha). \notag
\end{align}
\end{itemize}
\end{itemize}
\emph{If this is so, then $e_{B/F} = 1$, $\bar B \cong \bar
L(\sqrt{\bar \alpha})$ in the first alternative of} (II),
\emph{while $e_{B/F} = 2$, $\bar B \cong \bar L$ otherwise.}
\hfill $\square$ \lz

\noindent On the basis of the preceding results, it can now be
shown that the three valuation data we are interested in can be
pre-assigned in advance pretty much arbitrarily once the obvious
constraints are taken into account:

\subsection{Corollary.} \label{c.PREASINV} \emph{Let $e,n,r$ be
integers with $n \geq 0$ and $A$ an $\bar F$-algebra of dimension
$2^n$. There exists a composition division algebra $C$ over $F$
satisfying $\bar C \cong A$, $e_{C/F} = e$, $\tex(C) = r$ if and
only if the following conditions are fulfilled.}
\begin{itemize}
\item [(a)] \emph{$e \in \{1,2\}$, $0 \leq n \leq 4 - e$, $0 \leq
r \leq e_F$, $r = e_F$} (\emph{for $n = 0$, $e = 1$}).

\item [(b)] \emph{$A$ is a composition division algebra for $r =
0$, and $A/\bar F$ is a purely inseparable field extension of
characteristic $2$ and exponent at most $1$ for $r > 0$.}
\end{itemize}

\Proof If $C$ is a composition division algebra of dimension
$2^m$, $0 \leq m \leq 3$, over $F$ having the prescribed valuation
data, then $e = e_{C/F} \in \{1,2\}$ by
(\ref{ss.VADAL}.\ref{RAIN}), $0 \leq r = \tex(C) \leq e_F$ by
(\ref{e.BAFI}.\ref{ETE}) and $n = m$ or $n = m - 1$ according as
$e_{C/F}$ is $1$ or $2$ by Prop.~\ref{p.EFEN}; moreover, $n = 0$
and $e = 1$ imply $C \cong F$. Summing up and observing
(\ref{e.BAFI}.\ref{TEEF}), we obtain (a), while (b) follows from
Prop.~\ref{p.PROTE}~(d) combined with the remark to
Prop.~\ref{p.VARIID}. Conversely, suppose (a) and (b) are
fulfilled. If $r = 0$, the existence of a composition division
algebra over $F$ with the desired properties follows from from
\cite[Thms.~1,2]{MR51:635}. We may therefore assume $r > 0$, which
by (b) implies that $A/\bar F$ is a purely inseparable field
extension of characteristic $2$, exponent at most $1$ and degree
$2^n$. We first consider the case $n = 0$. If $e = 1$, then $C =
F$ has the prescribed valuation data since $r = e_F$ by (a). If $e
= 2$, Prop.~\ref{p.SEQUAT} yields a separable quadratic field
extension $L/k$ having $e_{L/F} = 2$, $\tex(L) = r$. We may
therefore assume $n > 0$. By Prop.~\ref{p.SEQUA}, there exists a
separable quadratic field extension $L/F$ such that $e_{L/F} = 1$,
$\tex(L) = r$ and $\bar L \subseteq A$. Hence $A = \bar L$, or $A
= \bar L(\sqrt{\bar \alpha})$ for some $\alpha \in \mfo$, $\bar
\alpha \notin \bar L^2$, or $A = \bar L(\sqrt{\bar
\alpha},\sqrt{\bar \beta}\,)$ for some $\alpha,\beta \in \mfo$,
$\bar \alpha \notin \bar L^2$, $\bar \beta \notin \bar
L(\sqrt{\bar \alpha})^2$ according as $n = 1,2,3$. In each case,
Prop.~\ref{p.SEQUA} and Thms.~\ref{t.NOFOO},~\ref{t.NOFOQ} yield a
composition division algebra $C$ over $F$ having the prescribed
valuation data. \hfill $\square$ \lz

\noindent Not surprisingly, however, the above valuation data are
far away from \emph{classifying} composition division algebras
over $2$-Henselian fields. This fact is underscored by the
following example.

\subsection{Example.} \label{e.PAIRS} Suppose $\bar F$ has
characteristic $2$ and let $\bar F \subseteq K^\prime \subseteq
L^\prime$ be a chain of purely inseparable field extensions having
exponent at most $1$ over $\bar F$ such that $[L^\prime:\bar F] =
2^n$, $[K^\prime:\bar F] = 2^{n-1}$, $n = 2,3$. Furthermore,
suppose $r,s \in \IZ$ satisfy the relations $0 < s < r \leq e_F$.
Then Cor.~\ref{c.PREASINV} yields a composition division algebra
$B$ over $F$ with
\begin{align}
\label{BEEF} e_{B/F} = 1, \quad \bar B \cong K^\prime, \quad
\tex(B) = r.
\end{align}
We claim \emph{there are an infinite number of mutually
non-isomorphic composition division algebras $C$ over $F$
containing $B$ as a subalgebra and having $e_{C/F} = 1$, $\bar C
\cong L^\prime$, $\tex(C) = s$.} This is in stark contrast to the
fact that unramified composition division algebras up to
isomorphism are uniquely determined by their residue algebras
\cite[Thm.~1]{MR51:635}.

To prove our claim, we fix an element $\beta \in \mfo$ with
$L^\prime = K^\prime(\sqrt{\bar \beta})$ and have $\bar \beta
\notin \bar B^2$ by \eqref{BEEF}, which implies
$\overline{\alpha\beta} \notin \bar B^2$ for any $\alpha \in
\mfo^\times$ with $\bar \alpha \in \bar B^2$. Setting
\[
d := 2(r - s) = 2 \left(\tex(B) - s\right) > 0, \quad \mu_\alpha := 1
- \pi^d\alpha\beta \in \mfo, \quad C_\alpha :=\mbC(B,\mu_\alpha),
\]
we conclude from Prop.~\ref{p.CDUNEV} that $\mu_\alpha \in
\mfo^\times$ and $C_\alpha$ is a composition division algebra over
$F$ satisfying
\[
e_{C_\alpha/F} = 1, \quad \overline{C_\alpha} \cong L^\prime,
\quad \tex(C_\alpha) = s.
\]
It therefore suffices to show, for any additional element
$\alpha^\prime \in \mfo^\times$ with $\overline{\alpha^\prime} \in
\bar B^2$, that $C_\alpha \cong C_{\alpha^\prime}$ implies $\bar
\alpha = \overline{\alpha^\prime}$. To see this, we write $\bar
\alpha = \delta^2$, $\overline{\alpha^\prime} = \delta^{\prime 2}$
with $\delta,\delta^\prime \in \bar B$ and invoke
Cor.~\ref{c.CHARNEXOD}~(b) to derive the following chain of
implications.
\begin{align*}
C_\alpha \cong C_{\alpha^\prime} &\Longrightarrow
\overline{\alpha\beta} \equiv \overline{\alpha^\prime\beta}
\bmod\,\bar B^2 \\
&\Longrightarrow \exists\,\gamma \in \bar B: \delta^2\bar \beta =
\delta^{\prime 2}\bar \beta + \gamma^2 \\
&\Longrightarrow \exists\,\gamma \in \bar B: (\delta -
\delta^\prime)\sqrt{\bar \beta} + \gamma = 0 \\
&\Longrightarrow \delta = \delta^\prime \Longrightarrow\, \bar
\alpha = \overline{\alpha^\prime}.
\end{align*}
\hfill $\square$ \lz

%\noindent In the remainder of this section, we aim at a refinement
%of Cor.~\ref{c.PREASINV}. However, there are natural obstructions
%to the validity of this refinement in full generality. In order to
%describe these obstructions, we introduce the following concepts.

\section{Types of composition algebras and heights.} \label{s.PRAS}
The trace exponent has been our principal tool so far to detect
wildness in pointed quadratic spaces and related objects. Other
tools of this kind may be obtained by consulting the literature.
It is the purpose of the present section to recast these tools in
the setting of composition algebras and to compare them with the
trace exponent. We begin by describing what will turn out later
(see Cor.~\ref{c.DYCTYP} below) to be a dichotomy of composition
division algebras over a $2$-Henselian field $F$.

\subsection{Types of composition algebras.} \label{ss.TYCO} A
composition division algebra $C$ over $F$ is said to be of
\emph{unitary} (resp.~of \emph{primary})
\emph{type}\index{unitary type}
\index{primary type} if there exist an
associative composition division algebra $B$ over $F$ with
$e_{B/F} = 1$ and a scalar $\mu$ which is a unit (resp.~a prime
element) in $\mfo$ such that $C \cong \mbC(B,\mu)$. We record a
few easy but useful
observations. $\ssp$ \\
(a) By Prop.~\ref{p.LIFT}, $C$ is of unitary or of primary type,
provided $f_{C/F} > 1$ or $\mbc(F) \neq 2$. \vspace{-8pt}$\ssp$ \\
(b) For $C$ to be of primary type it is necessary by
Prop.~\ref{p.CDPR} that $C$ have ramification index $e_{C/F} = 2$.
$\ssp$ \\
(c) If $C$ is tame, then $C$ is of unitary (resp. of primary) type
if and only if $C$ is unramified (resp. ramified) \cite{MR51:635}.
$\ssp$ \\
(d) Suppose $F$ has characteristic $0$ and $\bar F$ has
characteristic $2$. A quadratic field extension of $F$ may or may
not be wild. But if it is, it is of primary (resp.~of unitary)
type if and only if it is tamely ramified (resp.\ wildly ramified
or wildly unramified) in the sense of \cite[p.~60]{MR868606}.

\subsection{Remark.} \label{ss.TRELS} At this stage we cannot rule
out the possibility that a composition division algebra $C$ over
$F$ is both of unitary and of primary type: conceivably, there
could exist composition division algebras $B,B^\prime$ over $F$ of
ramification index $1$, and a unit $\mu$ as well as a prime
element $\mu^\prime$ in $\mfo$ such that $\mbC(B,\mu) \cong C
\cong \mbC(B^\prime,\mu^\prime)$. For showing that this scenario
is actually impossible, the following improvement of
Prop.~\ref{p.WITE} will play a crucial role.

\subsection{Theorem.} \label{t.WITY} \emph{Let $C$ be a
composition division algebra over $F$.} $\ssp$ \\
(a) \emph{If $C$ is of primary type, then $\wid(C) = \tex(C)$.}
$\ssp$ \\
(b) \emph{Consider the following conditions on $C$.}
\begin{itemize}
\item [(i)] \emph{$\wid(C) = \tex(C) - \frac{1}{2}$.}

\item [(ii)] \emph{There are trace generators of $C$ belonging to
$\mfp_C$.}

\item [(iii)] \emph{$C$ is wild of unitary type and ramification
index $2$.}
\end{itemize}
\emph{Then the implications}
\begin{align*}
\text{(i)} \Longleftrightarrow \text{(ii)} \Longleftarrow
\text{(iii)}
\end{align*}
\emph{hold. Moreover, if $f_{C/F} > 1$ or $\mbc(F) \neq 2$, then
all three conditions are equivalent.}

\Proof We begin with the first part of (b).

(i) $\Longrightarrow$ (ii). Let $w_0$ be a regular trace generator
of $C$. Then (i) and (\ref{p.WITE}.\ref{WIDTEX}) show
$\lambda_C(w_0) > 0$, hence $w_0 \in \mfp_C$.

(ii) $\Longrightarrow$ (i). Let $w_0 \in \mfp_C$ be a trace
generator of $C$. Then $\wid(C) \leq \lambda(t_C(w_0)) -
\lambda_C(w_0) < \tex(C)$, and Prop.~\ref{p.WITE}~(a) gives (ii).

(iii) $\Longrightarrow$ (ii). If $C$ is wild of unitary type and
ramification index $2$, we have $C = \mbC(B,\mu)$, $B$ an
associative composition division algebra over $F$ with $e_{B/F} =
1$, $\mu \in \mfo^\times$, and Thm.~\ref{t.ENBARE} shows that we
are in the situation of Prop.~\ref{p.CDUNOD} with $Q = C$, $P = B$
(cf. the translation formalism \ref{ss.TRANS}). Picking a trace
generator $u$ of $B$, we apply
(\ref{p.CDUNOD}.\ref{CDPROD}),(\ref{p.CDUNOD}.\ref{TREXPROD}) to
obtain
\begin{align*}
\lambda\big(t_C(u\Pi)\big) =\,\,&
\lambda\Big(t_C\big(\pi^{-\frac{d-1}{2}}(u + uj)\big)\Big) =
\lambda\big(t_B(u)\big) - \frac{d-1}{2} \\
=\,\,&\tex(B) - \frac{d-1}{2} = \tex(C),
\end{align*}
so $u\Pi \in \mfp_C$ is a trace generator of $C$, showing (ii).
Before completing the proof of (b), we turn to

(a) By hypothesis, we are in the situation of Prop.~\ref{p.CDPR}
with $Q = C$, $P = B$ as before, and by
(\ref{p.CDPR}.\ref{OCEPR}), an element of $\mfp_C$ has the form $y
= u + vj$, $u \in \mfp_B$, $v \in \mfo_B$. Hence $\pi^{-1}u \in
\mfo_B$ since $e_{B/F} = 1$, and from (\ref{p.CDPR}.\ref{VADAPR})
we conclude
\[
\lambda\big(t_C(y)\big) = \lambda\big(t_B(u)\big) =
\lambda\big(t_B(\pi^{-1}u)\big) + 1 \geq \tex(B) + 1 = \tex(C) + 1
\tex(C).
\]
Thus $\mfp_C$ does not contain trace generators of $C$, violating
(ii), hence (i), in (b). Now Prop.~\ref{p.WITE}~(a) implies (a).

It remains to show (i) $\Rightarrow$ (iii) in (b) under the
assumption $f_{C/F} > 1$ or $\mbc(F) \neq 2$. Then $e_{C/F} = 2$ by
Prop.~\ref{p.WITE}~(c), and $\tex(C) = \wid(C) + \frac{1}{2}
> 0$ by (i) forces $C$ to be wild. Moreover, by
\ref{ss.TYCO}~(a), it is of unitary or of primary type, the latter
case being excluded by (i) and (a). \hfill $\square$

\subsection{Corollary.} \label{c.DYCTYP} \emph{A composition
division algebra $C$ over $F$ cannot be both of unitary and of
primary type.}

\Proof If $C$ is tame, the assertion follows immediately from
\ref{ss.TYCO}~(c). If $C$ is wild of ramification index $1$, it
cannot be of primary type by \ref{ss.TYCO}~(b). Hence we are left
with the case that $C$ is wild of ramification index $2$. Then, by
Thm.~\ref{t.WITY}, $\wid(C) = \tex(C)$ is an integer if $C$ is of
primary type, while $\wid(C) = \tex(C) - \frac{1}{2}$ is not if
$C$ is of unitary type. \hfill $\square$

\subsection{Corollary.} \label{c.PRIRA} \emph{Let $C$ be a
composition division algebra of primary type over $F$ and suppose
$B \subseteq C$ is a composition subalgebra with}
\[
e _{B/F} = 1, \quad \mbd_k(B) = \frac{1}{2}\mbd_k(C).
\]
\emph{Then $\tex(B) = \tex(C)$.}

\Proof Since $C$ is not of unitary type by Cor.~\ref{c.DYCTYP}, we
obtain $C \cong \mbC(B,\mu)$ for some prime element $\mu \in
\mfo$. Hence $\tex(B) = \tex(C)$ by (\ref{p.CDPR}.\ref{VADAPR}).
\hfill $\square$ \lz

\noindent Both Cor.~\ref{c.TEXPCONJ} and Thm.~\ref{t.WITY}~(b)
fail in the exceptional cases stated therein. This is the upshot
of the following example.

\subsection{Example.} \label{e.SEPQUA} Let $L/F$ be a separable
quadratic field extension of ramification index $2$ and suppose
$\mbc(\bar F) = 2$. Then we are in the situation of of
Prop.~\ref{p.SEQUAT}, and if $r \leq e_F$, then $\tex(L) = r$ by
(\ref{p.SEQUAT}.\ref{TRLR}). Moreover, (\ref{p.SEQUAT}.\ref{PRAM})
gives $\lambda(t_L(\Pi)) = \lambda(\pi^r\alpha) = r = \tex(L)$, so
$\Pi \in \mfp_L$ is a trace generator of $L$; in particular, for
$\mbc(F) = 2$, parts (i),(ii) of Thm.~\ref{t.WITY}~(b) hold but
(iii) doesn't. On the other hand, if $r > e_F$, then $F$ has
characteristic zero and
\[
L = F\big(\sqrt{\pi\gamma}\big) = \mbC(F,\pi\gamma), \quad \gamma
:= -\beta + \frac{\pi^{2r-1}}{4}\alpha^2 \in \mfo^\times,
\]
is of primary type,  forcing $\wid(L) = \tex(L)$ by
Thm.~\ref{t.WITY}~(a).

Let $u = \gamma + \delta\Pi \in \mfo_L$, $\gamma,\delta \in \mfo$
(cf. (\ref{p.SEQUAT}.\ref{OLR})). Applying
(\ref{p.SEQUAT}.\ref{VALR}),
\begin{align*}
\lambda_L(u - u^\ast) =\,\,&\lambda_L(-\pi^r\alpha \delta +
2\delta\Pi) =
\mbm\,\{r + \lambda(\delta),e_F + \lambda(\delta) + \frac{1}{2}\} \\
\geq\,\,&\mbm\,\{r,e_F + \frac{1}{2}\},
\end{align*}
and this minimum is attained for, e.g., $u = \Pi$. Thus,
by(\ref{p.SEQUAT}.\ref{TRLR}),
\[
\mbm\,\{\lambda_L(u - u^\ast) \mid u \in \mfo_L\} =
\begin{cases}
\tex(L) &\text{for $r \leq e_F$,} \\
\tex(L) + \frac{1}{2} &\text{for $r > e_F$},
\end{cases}
\]
so in the latter case, the conclusion of Cor.~\ref{c.TEXPCONJ}
does not hold.

We remark in closing that Cor.~\ref{c.TEXPCONJ} also fails if
$L/F$ is \emph{tame} of ramification index $2$ since this implies
$\tex(L) = 0$ while $^\ast$ induces the identity on $\bar L$, so
$\lambda_L(u - u^\ast) > 0$ for all $u \in \mfo_L$.

\subsection{Remark.} \label{r.COSLO} There is an alternate way of
proving Cor.~\ref{c.DYCTYP}, by working with quadratic forms. Let
$C$ be a composition division algebra over $F$ and suppose $C$ is
of unitary and of primary type. Then there are composition
division algebras $B,B^\prime$ over $F$ of ramification index $1$
and a unit $\mu$ as well as a prime element $\mu^\prime$ in $\mfo$
such that
\[
\mbC(B,\mu) \cong C \cong \mbC(B^\prime,\mu^\prime).
\]
Then Prop.~\ref{p.PFISUB} yields a scalar $\gamma \in F^\times$
satisfying
\begin{align}
\label{COMSLO} \mbC(B,\mu) \cong \mbC(B,\gamma), \quad
\mbC(B^\prime,\mu^\prime) \cong \mbC(B^\prime,\gamma)
\end{align}
Since $e_{B/F} = e_{B^\prime/F} = 1$, we conclude from
(\ref{ss.GEST}.\ref{GALA}) that $\lambda(n_B(B^\times)) =
\lambda(n_{B^\prime}(B^{\prime\times})) = 2\IZ$. Hence
Prop.~\ref{p.CDPOQUA}~(b) and the first relation of \eqref{COMSLO}
show that $\lambda(\gamma)$ is even, while
Prop.~\ref{p.CDPOQUA}~(b) and the second relation of
\eqref{COMSLO} show that $\lambda(\gamma)$ is odd, a
contradiction. \lz

\noindent With the aim of generalizing Thm.~\ref{t.SLIFT},
Cor.~\ref{c.COSUPRE} and Lemma~\ref{l.COSUPRE}, we next turn to the
problem of finding (chains of) subalgebras having ramification index
$1$ and pre-assigned trace exponents. Once the obvious constraints
are taken into account (provided, e.g., by Prop.~\ref{p.PROTE}~(b)
and Cor.~\ref{c.PRIRA}), we will show that chains of such
subalgebras always exist.

\subsection{Theorem.}
\label{t.SESFIPATE} \emph{Suppose $\bar F$ has characteristic $2$
and let $C$ be a composition division algebra of dimension $2^n$,
$n = 2,3$, over $F$ that is not a quaternion division algebra of
primary type. Given $r \in \IZ$, $\tex(C) \leq r \leq e_F$, there
exists a separable quadratic subfield $L \subseteq C$ with
$e_{L/F} = 1$, $\tex(L) = r$.}

\Proof The case $r = \tex(C)$ having been settled by
Thm.~\ref{t.SLIFT}, we may assume $r > \tex(C)$. Applying
Cor.~\ref{c.COSUPRE} and then Thm.~\ref{t.SLIFT}, we find a
composition subalgebra $B \subseteq C$ and a separable quadratic
subfield $K \subseteq B$ with
\begin{align}
%\label{MINFIL}
\mbd_F(B) = 2^{n-1}, \quad e_{K/F} = e_{B/F} = 1, \quad \tex(K) =
\tex(B) = \tex(C). \notag
\end{align}
Suppose for the time being that the case $n = 2$ has been solved
and let $n = 3$. The quaternion algebra $B$, having ramification
index $1$, cannot be of primary type, allowing us to apply the
case $n = 2$ to $B$ in place of $C$ and leading us to the desired
conclusion.

We are thus reduced to the case $n = 2$, which we will assume for
the rest  of the proof. Then $K = B$ and $C$ has dimension $4$,
hence is not of primary type by hypothesis. We are therefore lead
to a unit $\mu \in \mfo^\times$ with
\begin{align}
\label{CDKMU} C = \mbC(K,\mu) = K \oplus Kj.
\end{align}
Let us first assume $\tex(K) = \tex(C) = 0$, so $K$ is tame. Since
$C$ is a division algebra, we conclude $\bar \mu \notin n_{\bar
K}(\bar K^\times)$ from Prop.~\ref{p.LATAU}. In particular, we have
$\bar \mu \notin \bar F^2$. The hypothesis $\tex(K) = 0$ yields an
element $v \in \mfo_K$ having $t_K(v) = 1$. Observing \eqref{CDKMU},
we put
\[
w := \pi^rv + j \in C
\]
and obtain $t_C(w) = \pi^r$, $n_C(w) = \pi^{2r}n_K(v) - \mu$,
$\overline{n_C(w)} = \bar \mu \notin \bar F^2$. In particular, $w
\in \mfo_C^\times \setminus F1_C$ and we conclude from
Prop.~\ref{p.SEQUA} that $L = F[w] \subseteq C$ is a separable
quadratic subfield of the desired kind.

We are left with the case $\tex(K) = \tex(C) > 0$. Then $K$ is
wild, and Cor.~\ref{c.TREL} yields an element $u \in
\mfo_K^\times$ with
\begin{align}
\label{KATREL} t_K(u) = \pi^{\tex(K)}, \quad \overline{n_K(u)}
\notin \bar F^2.
\end{align}
By \eqref{CDKMU} and Cor.~\ref{c.EQTE}, we may assume
\begin{align}
\label{MUALT} \bar \mu \notin \bar F^2 \quad \text{or} \quad \mu =
1 - \pi \beta \quad \text{for some $\beta \in \mfo^\times$}
\end{align}
according as $C$ has ramification index $1$ or $2$. We now put
\begin{align*}
w :=
\begin{cases}
\pi^{r-\tex(K)}u + j \in C &\text{if $e_{C/F} = 1$}, \\
\pi^{r-\tex(K)}u + uj \in C &\text{if $e_{C/F} = 2$.}
\end{cases}
\end{align*}
Then $t_C(w) = \pi^r$ by \eqref{KATREL}, and for $e_{C/F} = 1$ we
obtain $\overline{n_C(w)} = \bar\mu \notin \bar F^2$ by
\eqref{MUALT}. Similarly, for $e_{C/F} = 2$, we obtain
$\overline{n_C(w)} = \overline{n_K(u)} \notin \bar F^2$ by
\eqref{KATREL},\eqref{MUALT}. In either case, $w \in \mfo_C^\times
\setminus k1_C$ generates a separable quadratic subfield $L :=
F[w] \subseteq C$ of ramification index $e_{L/F} = 1$ and trace
exponent $\tex(L) = r$ (Prop.~\ref{p.SEQUA}). \hfill $\square$

\subsection{Corollary.} \label{c.FILPREAS} \emph{Suppose
$\bar F$ has characteristic $2$ and let $C$ be an octonion
division algebra over $F$ that is not of primary type. Given $r,s
\in \IZ$, $\tex(C) \leq r \leq s \leq e_F$, there exists a
filtration $L \subseteq B \subseteq C$ consisting of a quaternion
subalgebra $B \subseteq C$ and a separable quadratic subfield $L
\subseteq B$ such that $e_{L/F} = e_{B/F} = 1$ and $\tex(B) = r$,
$\tex(L) = s$.}

\Proof It suffices to construct a quaternion subalgebra $B
\subseteq C$ with $e_{B/F} = 1$, $\tex(B) = r$ because
Thm.~\ref{t.SESFIPATE} applies to such a $B$ and also yields an
$L$ with the desired properties.

To construct $B$, we first invoke Thm.~\ref{t.SLIFT} and
Cor.~\ref{c.COSUPRE} to find a quaternion subalgebra $B_1
\subseteq C$ and a separable quadratic subfield $L_1 \subseteq
B_1$ satisfying
\begin{align}
\label{MINFILONE}e_{L_1/F} = e_{B_1/F} = 1, \quad \tex(L_1) =
\tex(B_1) = \tex(C).
\end{align}
By hypothesis, \eqref{MINFILONE} and Cor.~\ref{c.EQTE} yield units
$\mu_1,\mu_2 \in \mfo^\times$ with
\begin{align}
\label{CDMUONE} B_1 =\,\,&\mbC(L_1,\mu_1) = L_1 \oplus L_1j_1,
\quad \overline{\mu_1} \notin \bar F^2, \\
C =\,\,&\mbC(B_1,\mu_2) = B_1 \oplus B_1j_2, \notag
\end{align}
where
\begin{align*}
\text{$\overline{\mu_2} \notin \bar F^2$ for $e_{C/F} = 1$ and
$\mu_2 = 1 - \pi \beta_2$, $\beta_2 \in \mfo^\times$ for $e_{C/F}
= 2$.}
\end{align*}
By the same token,
\[
B_2 := \mbC(L_1,\mu_2) = L_1 \oplus L_1j_2 \subseteq C
\]
is a quaternion subalgebra not of primary type by
Cor.~\ref{c.DYCTYP} with $\tex(B_2) = \tex(L_1) = \tex(C) \leq r
\leq e_F$, $e_{B_2/F} = e_{C/F}$. Hence Thm.~\ref{t.SESFIPATE}
leads us to a separable quadratic subfield $L \subseteq B_2$ with
$e_{L/F} = 1$, $\tex(L) = r$ and \eqref{CDMUONE} combines with
Cor.~\ref{c.EQTE}~(b) to show that the quaternion subalgebra
\[
B = \mbC(L,\mu_1) = L \oplus Lj_1 \subseteq C
\]
satisfies $e_{B/F} = 1$, $\tex(B) = \tex(L) = r$. \hfill $\square$

\subsection{Heights.} \label{ss.HGT} Over a Henselian field
having residual characteristic $p > 0$, the height as an important
invariant of a central associative division algebra of degree $p$
over $F$ has been considered by Saltman \cite[pp.~
1757, 1765-6]{MR589083} (who uses the term ``level"\index{Saltman's level $\hcom$ ($= \wid$)}), Kato \cite[\S~1]{Kato:gen1} (who calls it the ``ramification number")\index{Kato's ramification number $\hcom$}, and Tignol
\cite[3.2]{Tignol:wild}. It is, in particular, Tignol's approach
that suggests two immediate translations to the setting of
composition algebras.

Let $C$ be a composition division algebra over our $2$-Henselian
field $F$. We use the maps $\hcom\:C^\times \times C^\times \to
\IQ_\infty$, $\hass:C^\times \times C^\times \times C^\times \to
\IQ_\infty$ given by
\begin{align}
\label{HCOM} \hcom(x,y) :=\,\,&\lambda_C([x,y]) - \lambda_C(x) -
\lambda_C(y) \geq 0, \\
\label{HASS} \hass(x,y,z) :=\,\,&\lambda_C([x,y,z]) - \lambda_C(x)
- \lambda_C(y) - \lambda_C(z) \geq 0
\end{align}
for all $x,y,z \in C^\times$ to define
\begin{align}
\label{HGCOM} \hgcom(C) :=\,\,&\mathrm{inf}\,\{\hcom(x,y) \mid x,y
\in
C^\times\}, \\
\label{HGASS} \hgass(C) :=\,\,&\mathrm{inf}\,\{\hass(x,y,z) \mid
x,y,z \in C^\times\}
\end{align}
and to call these numbers the \emph{commutative
height}\index{height ($\hcom$ and $\hass$)} and the \emph{associative
height} of $C$, respectively. If $C$ has
dimension at most $2$, then $\hgcom(C) = \hgass(C) = \infty$. On
the other hand, if $C$ is a quaternion algebra, then $\hgcom(C) <
\infty = \hgass(C)$ and $\hgcom(C)$ agrees with what Tignol calls
its height, while if $C$ is an octonion algebra, its commutative
and its associative height are both finite.

A general theorem of Tignol \cite[3.12]{Tignol:wild} implies
$\hgcom(C) = \wid(C)$ for any quaternion division algebra $C$ over
$F$. This special observation is part of a much more
general picture that will be summarized in the following theorem,
whose proof in the quaternionic case is independent of
\cite{Tignol:wild} and, in fact, works uniformly in the octonionic
case as well.

\subsection{Theorem.} \label{t.HEWI} \emph{If $C$ is a quaternion
division algebra over $F$, then}
\[
\hgcom(C) = \wid(C).
\]
\emph{If $C$ is an octonion division algebra over $F$, then}
\[
\hgcom(C) = \hgass(C) = \wid(C).
\]

\begin{proof} Let $C$ be a composition division algebra of dimension
$2^n$, $n = 2,3$, over $F$. We must show
\begin{align}
\label{HECWI} \hgcom(C) =\,\,&\wid(C), \\
\label{HECAWI} \hgcom(C) =\,\,&\hgass(C) = \wid(C) &&(\text{for $n
= 3$}).
\end{align}
To do so, we combine Thm.~\ref{t.NOAS} with
(\ref{p.ELAC}.\ref{LACT}) to obtain
%\begin{align}
%\label{NOCOT} n_C([x_1,x_2]) =\,\,&4n_C(x_1)n_C(x_2) -
%t_C(x_1)^2n_C(x_2) - t_C(x_2)^2n_C(x_1) + \\
%\,\,&t_C(x_1x_2)t_C(x_1x_2^\ast), \notag \\
%\label{NOAST} n_C([x_1,x_2,x_3]) =\,\,&4n_C(x_1)n_C(x_2)n_C(x_3) -
%\sum\,t_C(x_i)^2n_C(x_j)n_C(x_l) + \\
%\,\,&\sum\,t_C(x_ix_j)t_C(x_ix_j^\ast)n_C(x_l) -
%t_C(x_1x_2)t_C(x_2x_3)t_C(x_3x_1) + \notag\\
%\,\,&t_C(x_1x_2x_3)t_C(x_2x_1x_3) \notag
%\end{align}
%for all $x_1,x_2,x_3 \in C$. Since
%\[
%\lambda\big(t_C(x)\big) \geq \wid(C) + \lambda_C(x)
%\]
%for all $x \in C$ by (\ref{ss.WID}.\ref{DEWI}), a straightforward
%application of
%(\ref{ss.MOSTLY}.\ref{VACE}),(\ref{e.BAFI}.\ref{ETE}),\eqref{NOCOT},\eqref{NOAST}
%yields
\begin{align*}
%\lambda_C([x_1,x_2]) \geq\,\,&\wid(C) + \lambda_C(x_1) +
%\lambda_C(x_2), \\
\lambda_C([x_1,x_2,x_3]) \geq\,\,&\wid(C) + \lambda_C(x_1) +
\lambda_C(x_2) + \lambda_C(x_3)
\end{align*}
for all $x_1,x_2,x_3 \in C$. Combining this and
(\ref{p.ELAC}.\ref{LACCO}) with
(\ref{ss.HGT}.\ref{HCOM}$\--$\ref{HGASS}), we conclude
\begin{align}
\label{HGTGEQ} \hgcom(C) \geq \wid(C), \quad \hgass(C) \geq
\wid(C).
\end{align}
To complete the proof of \eqref{HECWI},\eqref{HECAWI}, it therefore
suffices to show that
\begin{align}
\label{HELWI} \hcom(x_1,x_2) =\,\,&\wid(C) &&\text{(for some $x_1,x_2 \in C$)}, \\
\label{HELAWI} \hass(x_1,x_2,x_3) =\,\,&\wid(C) &&(\text{for $n =
3$ and some $x_1,x_2,x_3 \in C$}).
\end{align}
For this purpose, we require two additional formulas: suppose $C =
\mbC(B,\mu) = B \oplus Bj$ is a Cayley-Dickson construction as in
\ref{ss.CDCO}, for some associative composition algebra $B$ over
$F$ and some scalar $\mu \in F^\times$. Then a straightforward
application of (\ref{ss.CDCO}.\ref{CDMU}) yields
\begin{align}
\label{CDCOM} [u,j] =\,\,&(u - u^\ast)j &&(u \in B), \\
\label{CDASS} [u_1,u_2,j] =\,\,&[u_1,u_2]j &&(u_1,u_2 \in B).
\end{align}
In order to establish \eqref{HELWI},\eqref{HELAWI}, we distinguish
the following cases. $\ssp$ \\
\emph{Case~$1$.} $e_{C/F} = 1$. \\
Then $\wid(C) = \tex(C)$ by Prop.~\ref{p.WITE}~(c). If $C$ is
tame, the $\bar C$ is a composition division algebra of dimension
$2^n$ over $\bar F$, so there are $x_1,x_2,x_3 \in \mfo_C^\times$
with $[x_1,x_2] \in \mfo_C^\times$, and even $[x_1,x_2,x_3] \in
\mfo_C^\times$ for $n = 3$. By
(\ref{ss.HGT}.\ref{HCOM}),(\ref{ss.HGT}.\ref{HASS}), this implies
$\hcom(x_1,x_2) = 0 = \tex(C)$, and even $\hass(x_1,x_2,x_3) = 0 =
\tex(C)$ for $n = 3$, proving \eqref{HELWI},\eqref{HELAWI} in the
tame case.

If $C$ is wild, it must be of unitary type since $e_{C/F} = 1$, and
Thm.~\ref{t.ENBARE} combined with Prop.~\ref{p.LIFT} implies $C=
\mbC(B,\mu)$, $B$ an associative composition division algebra over
$F$ with $e_{B/F} = 1$, $\mu = 1 - \pi^d\beta$, $d \in \IZ$ even, $0
\leq d <2\tex(B)$, $\beta \in \mfo$, $\bar\beta \notin \bar B^2$. In
particular, taking into account \ref{ss.TRANS}, we are in the
situation of Prop.~\ref{p.CDUNEV}. Applying Cor.~\ref{c.TEXPCONJ},
we find an element $u \in \mfo_B^\times$ such that $\lambda_B(u -
u^\ast) = \tex(B)$. Hence
(\ref{p.CDUNEV}.\ref{CDPREV}),(\ref{p.CDUNEV}.\ref{DEPREV}),(\ref{p.CDUNEV}.\ref{TREXPREV})
and \eqref{CDCOM} yield
\begin{align*}
\hcom(u,\Xi) =\,\,&\lambda_C([u,\Xi]) - \lambda_C(u) -
\lambda_C(\Xi) = \lambda_C(\pi^{-\frac{d}{2}}[u,j]) \\
=\,\,&\lambda_B(u - u^\ast) - \frac{d}{2} = \tex(B) - \frac{d}{2}
= \tex(C).
\end{align*}
Thus \eqref{HELWI} holds for $x_1 = u$, $x_2 = \Xi$, and we have
established \eqref{HECWI} in Case~$1$.

If $n = 3$ ($C$ still assumed to be wild), Case~$1$ applies to $B$,
so \eqref{HELWI} yields elements $u_1,u_2 \in B$ such that
$\hcom(u_1,u_2) = \tex(B)$. Hence \eqref{CDASS} and
(\ref{p.CDUNEV}.\ref{TREXPREV}) imply
\begin{align*}
\hass([u_1,u_2,\Xi])
=\,\,&\lambda_C(\pi^{-\frac{d}{2}}[u_1,u_2,j]) -\lambda_B(u_1) -
\lambda_B(u_2) \\
=\,\,&\hcom(u_1,u_2) - \frac{d}{2} = \tex(B) - \frac{d}{2} =
\tex(C),
\end{align*}
giving \eqref{HELAWI} for $x_1 = u_1$, $x_2 = u_2$, $x_3 = \Xi$,
and settling Case~$1$ completely. $\ssp$ \\
\emph{Case~$2$.} $e_{C/F} = 2$. \\
If $C$ is of primary type, then $\wid(C) = \tex(C)$ by
Thm.~\ref{t.WITY}~(a), and $C = \mbC(B,\mu)$ with $B$ as before and
$\mu$ a prime element in $\mfo$. Applying Thm.~\ref{t.ENBARE}~(a)
and picking $u \in \mfo_B^\times$ with $\tex(C) = \tex(B) =
\lambda_B(u - u^\ast)$, we obtain
\begin{align*}
\hcom(u,j) =\,\,&\lambda_C([u,j]) - \lambda_B(u) - \lambda_C(j) \\
=\,\,&\lambda_C\big((u - u^\ast)j\big) - \lambda_C(j) =
\lambda_B(u - u^\ast) = \tex(C),
\end{align*}
and \eqref{HELWI} holds for $x_1 = u$, $x_2 = j$. Moreover, for $n
= 3$, Case~$1$ applies to $B$ and yields $u_1,u_2 \in B^\times$
having $\hcom(u_1,u_2) = \tex(B) = \tex(C)$, allowing us to
compute
\begin{align*}
\hass(u_1,u_2,j) =\,\,&\lambda_C([u_1,u_2]j) - \lambda_B(u_1) -
\lambda_B(u_2) - \lambda_C(j) \\
=\,\,&\hcom(u_1,u_2) = \tex(C)
\end{align*}
and completing the proof of \eqref{HELAWI} for $x_1 = u_1$, $x_2 =
u_2$, $x_3 = j$.

Now suppose $C$ is of unitary type. Then $C$ is wild since
$e_{C/F} = 2$, and Thm.~\ref{t.WITY}~(b) shows $\wid(C) = \tex(C)
- \frac{1}{2}$. This time, Thm.~\ref{t.ENBARE} implies $C =
\mbC(B,\mu)$ with $B$ as before, $\mu = 1 - \pi^d\beta$, $d \in
\IZ$ odd, $0 \leq d <2\tex(B)$, $\beta \in \mfo^\times$, so in
view of \ref{ss.TRANS} we are in the situation of
Prop.~\ref{p.CDUNOD}. Picking again an element $u \in
\mfo_B^\times$ with $\tex(B) = \lambda_B(u - u^\ast)$, we obtain
\begin{align*}
\hcom(u,\Pi) =\,\,&\lambda_C([u,\Pi]) - \lambda_B(u) -
\lambda_C(\Pi) =\lambda_C(\pi^{-\frac{d-1}{2}}[u,j]) - \frac{1}{2}
\\
=\,\,&\lambda_C\big((u - u^\ast)j\big) - \frac{d-1}{2} -
\frac{1}{2} = \lambda_B(u - u^\ast) - \frac{d-1}{2} - \frac{1}{2}
\\
=\,\,&\tex(B) - \frac{d-1}{2} - \frac{1}{2} = \tex(C) -
\frac{1}{2} = \wid(C),
\end{align*}
and \eqref{HELWI} holds for $x_1 = u$, $x_2 = \Pi$. If, in
addition, $n = 3$, then Case~$1$ applies to $B$, yielding $u_1,u_2
\in B^\times$ with $\hcom(u_1,u_2) = \tex(B)$. Hence
\begin{align*}
\hass(u_1,u_2,\Pi) =\,\,&\lambda_C([u_1,u_2,\Pi]) - \lambda_B(u_1)
- \lambda_B(u_2) - \lambda_C(\Pi) \\
=\,\,&\lambda_C(\pi^{-\frac{d-1}{2}}[u_1,u_2,j]) - \lambda_B(u_1)
- \lambda_B(u_2) - \frac{1}{2} \\
=\,\,&\lambda_C([u_1,u_2]j) - \frac{d-1}{2} - \lambda_B(u_1) -
\lambda_B(u_2) - \frac{1}{2} \\
=\,\,&\hcom(u_1,u_2) - \frac{d-1}{2} - \frac{1}{2} = \tex(B) -
\frac{d-1}{2} - \frac{1}{2}  \\
=\,\,&\tex(C) - \frac{1}{2} = \wid(C).\qedhere
\end{align*}
\end{proof}

%%%%%%%%%%%%%%%%%%%%%%%%%%%%%%%%%%%%%%%%%%%%%%%%%%
\part{Connections with $K$-theory} \label{p.CONKAT}

\section{Introduction to Part \ref{p.CONKAT}} \label{s.CONKAT}

\subsection{} \label{p3.intro}
The goal of this part of the paper is to translate the results of
Part \ref{hensel.part} into the language of Kato's filtration on
Milnor $K$-theory mod 2.  Because it costs little extra, we will
also give a dictionary relating traditional valuation-theoretic
terms on associative division algebras of prime degree $p$ with
Milnor $K$-theory mod $p$.

\emph{For the remainder of the paper}, we fix a field $F$ of
characteristic zero that has a Henselian discrete valuation $\val$
and residue field $\Fb$ of prime characteristic $p$.  We assume
that $F$ contains a primitive $p$-th root of unity $\zeta$ and set
\[
m := p \cdot \frac{\val(p)}{p - 1}.
\]
This is an integer divisible by $p$ because $\val(p)/(p - 1) =
\val(\zeta - 1)$, see, e.g., \cite[4.1.2(i)]{CT:kato}.

Recall that the Milnor $K$-ring of $F$, denoted $K^M_*(F)$, is the
tensor algebra (over $\Z$) of the abelian group $\Fx$ modulo the
``Steinberg relation" $a \ot (1 -  a) = 0$ for $a \in F$, $a
\neq 0,1$.

One writes $\{ a_1, \ldots, a_q \}$ for the image of $a_1 \ot \cdots
\ot a_q$ in $K^M_q(F)$.  We put $k_q(F)$ for $K^M_q(F)/p$, and we
abuse notation by writing $\{ a_1, \ldots, a_q \}$ also for the
image of that element in $k_q(F)$; such a class in $k_q(F)$ is
called a \emph{symbol}. See, e.g., \cite{GilleSz} for basic
properties.  Kato, Bloch, and Gabber proved that $k_q(F)$ is
isomorphic to $H^q(F, \mu_p^{\ot q})$ via the ``Galois symbol",
which sends $\{ a_1, \ldots, a_q \} \mapsto (a_1) \cdot (a_2) \cdots
(a_q)$; this identifies nonzero symbols in $k_q(F)$ with nonzero
symbols in $H^q(F, \mu_p^{\ot q})$.  We are mainly interested in the
following cases:
\begin{enumerate}
\refstepcounter{equation}\renewcommand{\theenumi}{\theequation}
\item \underline{$q=1$}: $k_1(F)$ and $H^1(F, \mu_p)$ are
naturally identified with $\Fx / F^{\times p}$.  A nonzero element
$xF^{\times p}$ defines a degree $p$ extension $F(\chi)$ such that
$\chi^p = x$.\label{p3.field}\refstepcounter{equation}

\item
\underline{$q=2$}: $H^2(F, \mu_p^{\ot 2})$ is identified (via
$\zeta$) with $H^2(F, \mu_p)$, i.e., the $p$-torsion in the
Brauer group of $F$.  We fix the identification with the Brauer
group so that the nonzero symbol $\{ x, y \}$ in $k_2(F)$
is sent to the associative central division $F$-algebra of
dimension $p^2$ generated by elements $\chi, \psi$ satisfying\label{p3.div}\refstepcounter{equation}
\[
\chi^p = x, \quad \psi^p = y, \quad \text{and} \quad \chi \psi =
\zeta \psi \chi.
\]

\item \underline{$p = 2$}: In this case, there is a bijection between
symbols in $k_q(F)$ and $q$-Pfister (quadratic) forms given by
sending $\{ a_1, \ldots, a_q \}$ to $\pform{a_1, \ldots, a_q}$.  For
$q \le 3$, we can of course further identify anisotropic $q$-Pfister
forms with composition algebras of dimension $2^q$.\label{p3.comp}
\end{enumerate}

For a nonzero symbol $\gamma \in k_q(F)$ from cases \eqref{p3.field}
or \eqref{p3.div}, write $D$ for the corresponding division
$F$-algebra.  As the valuation $\val$ is Henselian, it extends to a
discrete valuation $\val_D$ on $D$ via the usual formula $\val_D(x)
:= \val(N_{F(x)/F}(x))/[F(x):F]$, cf.~(\ref{ss.CONCON}.\ref{EXCO}).
The definition of residue division algebra $\Db$, ramification index
$e_{D/F}$, etc., is the same as for quaternion algebras, and the
fundamental relation of Prop.~\ref{p.EFEN} holds, see
\cite[p.~393]{W:survey} for references.

\subsection{}
Below, we will recall the filtration\index{filtration on $K$-theory}
on $k_q(F)$ and define invariants $e_\gamma$\index{$e_\gamma$} (= 1
or $p$) and $\depth(\gamma)$\index{depth} for a symbol $\gamma \in
k_q(F)$ in terms of $K$-theory.  We will prove $K$-theoretic
analogues of the Local Norm Theorem \ref{t.HENOTH} (\S\ref{s.LNTK}),
the Normal Form Theorems \ref{t.NOFOO} and \ref{t.NOFOQ}
(\S\ref{gather.sec}), and Theorem \ref{t.CHAGORE}
(\S\ref{K.CHAGORE.sec}).  These proofs use the background results on
quadratic forms over Henselian fields from \S\ref{s.COALHE} but not
the deeper results from the rest of Part \ref{hensel.part}.

In the final sections of the paper
(\S\S\ref{dict.sec}--\ref{pf.texp.depth}) we give a dictionary
between properties of symbols in $k_q(F)$ in the cases $q =1$, $q
=2$, or $p = 2$.  The proofs in the case $p = 2$ rely heavily on the
full strength of the results in Part \ref{hensel.part}.

%%%%%%%%%%%%%%%%%%%%%%%%%%%%%%%%%%%%%%%%%%%%%%%%
\section{The filtration on $K$-theory} \label{filtration.sec}

We now recall the Kato filtration on Milnor
$K$-theory over $F$, together with the isomorphisms of the graded
components with various modules of differential forms, etc., over
$\Fb$.

\subsection{Filtration and depth} \label{filtration}
Write $\mfo$ for the valuation ring on $F$ and $\mfp$ for its maximal
ideal.  One can filter $\mfo$ as
\[
\mfo \setminus \{ 0 \} = U_0 \supseteq U_1 \supseteq U_2 \supseteq \cdots
\]
where $U_i := 1 + \mfp^i$ for $i \ge 1$.

For $q \ge 1$, setting $U^i k_q(F)$ to be additively generated by
elements $\{ u \} \cdot k_{q-1}(F)$ for $u \in U_i$ defines a
filtration
\[
k_q(F) = U^0 k_q(F) \supseteq U^1 k_q(F) \supseteq \cdots
\]
For $i > m$, $U_i$ consists of $p$-th powers
\cite[4.1.2(ii)]{CT:kato}, so $U^i k_q(F)$ is zero.

The \emph{depth} of $\gamma \in k_q(F)$ is the supremum of $\{ i
\mid \gamma \in U^i k_q(F) \}$.  The only element of depth $> m$
is zero, which has depth $\infty$.  The filtration is compatible
with the product in the sense that $U^i k_r(F) \cdot U^j k_s(F)
\subseteq U^{i+j} k_{r+s}(F)$ by \cite[4.1.1b]{CT:kato}.  Said
differently, for elements $\alpha \in k_r(F)$ and $\beta \in
k_s(F)$, we have
\begin{equation} \label{filt.1}
\depth(\alpha) + \depth(\beta) \le \depth(\alpha \cdot \beta);
\end{equation}
this inequality can be strict, see Example \ref{ram.eg}.

\subsection{Kato isomorphisms}\index{Kato isomorphisms}
For nonzero $\gamma \in k_q(F)$ with $q \ge 1$ of depth $d$, we
consider the (nonzero) image of $\gamma$ in $\gr^d k_q(F) := U^d
k_q(F) / U^{d+1} k_q(F)$; this is the \emph{initial form} of
$\gamma$.  The results of Kato, et al, include specific
isomorphisms:
\[
\gr^d k_q(F) \cong \begin{cases}
k_q(\Fb) \oplus k_{q-1}(\Fb) &\text{if $d = 0$;} \\
\Omega^{q-1} & \text{if $0 < d < m$ and $p$ does not divide $d$;}\\
\frac{\Omega^{q-1}}{Z^{q-1}} \oplus \frac{\Omega^{q-2}}{Z^{q-2}} &
\text{if $0 < d < m$ and $p$ divides $d$;} \\
H^1(\Fb, \nu(q-1)) \oplus H^1(\Fb, \nu(q-2))&\text{if $d = m$.}
\end{cases}
   \]
Here $\Omega^1$ denotes the $\Fb$-vector space of derivations $\Fb
\rightarrow \Fb$, $\Omega^q := \wedge^q \Omega^1$ for $q \geq
1$, $\Omega^0 = \Fb$ and $\Omega^{-1} = \{0\}$. The subspace
$Z^q$ is the kernel of the differential $\Omega^q \rightarrow
\Omega^{q+1}$, i.e., $Z^q$ is the subspace of exact forms. The
groups $\nu(q)$ are defined in terms of the Cartier operator
\cite[pp.~4,5]{CT:kato}; they are chosen so that $H^1(\Fb,
\nu(q-1))$ in characteristic $p$ plays the role of the Galois
cohomology group $H^q(K, \mu_p^{\otimes (q -1)})$ for $K$ of
characteristic $\ne p$.

%%\begin{rmk}\label{killed.by.p}
%%The targets of the Kato isomorphisms are abelian groups that are
%killed by $p$.  For $\Omega^q$ this is clear because it is a
%vector space over $\Fb$, and similarly for its subspace $Z^q$. The
%group $\nu(q)$ is a subgroup of $\Omega^q$ and so is also
%$p$-torsion, hence the same holds for the Galois cohomology group
%$H^1(\Fb, \nu(q))$.  This shows that multiplication by an integer
%not divisible by $p$ is invertible on these groups.
%%\end{rmk}

We refer to these isomorphisms as the ``Kato isomorphisms".  Fix a uniformizer $\pi$ for
$\val$ and write $a_i$ and $b$ for elements of $\Ox$.  The
isomorphisms are:
\[
\begin{array}{|c|r@{\mapsto}l|} \hline
d & \multicolumn{2}{c|}{\text{map}}\\ \hline\hline
d = 0&\{ a_1, a_2, \ldots, a_q \} & (\{ \ba_1, \ldots, \ba_q \}, 0) \\
& \{ \pi, a_1, \ldots, a_{q-1} \} & (0, \{ \ba_1, \ldots, \ba_{q-1} \}) \\
\hline \text{$p$ does not divide $d$}& \{ 1 + b \pi^d, a_1, \ldots,
a_{q-1} \}& \bb \dlog{\ba_1} \wedge \cdots \wedge \dlog{\ba_{q-1}}
\\ \hline
\text{$p$ divides $d$}&\{ 1 + b \pi^d, a_1, \ldots, a_{q-1} \} & (\bb
\dlog{\ba_1} \wedge \cdots \wedge \dlog{\ba_{q-1}}, 0) \\
\text{and $d \ne 0, m$}&\{ \pi, 1 + b \pi^d, a_1, \ldots, a_{q-2} \}
& (0, \bb \dlog{\ba_1} \wedge \cdots \wedge \dlog{\ba_{q-2}}) \\
\hline
d = m&\{ 1 + b (\zeta - 1)^p, a_1, \ldots, a_{q-1} \} &
(\bb \dlog{\ba_1} \wedge \cdots \wedge \dlog{\ba_{q-1}}, 0) \\
& \{ \pi, 1 + b (\zeta - 1)^p, a_1, \ldots, a_{q-2} \} & (0, \bb
\dlog{\ba_1} \wedge \cdots \wedge \dlog{\ba_{q-2}})  \\ \hline
\end{array}
\]
The description of $\gr^0 k_q(F)$ is a result of Bass-Tate that
holds without restriction on the characteristic of $\Fb$.  For
depth $m$, of course $(\zeta - 1)^p$ has value $m$, and in case $p
= 2$ the expression $(\zeta - 1)^p$ is 4.

\subsection{} \label{Kdiff}
We mention for later reference a useful fact about $k_q(\Fb)$ for $q
\ge 1$.  As with any field, there is a group homomorphism
$K^M_q(\Fb) \rightarrow \Omega^q$ defined by $\{ x_1, \ldots, x_q \}
\mapsto \dlog{x_1} \wedge \cdots \wedge \dlog{x_q}$ (one checks the
Steinberg relation).  But $\Fb$ has characteristic $p$, so this
homomorphism induces a homomorphism $\psi \!: k_q(\Fb) \rightarrow
\Omega^q$.  Moreover, $\psi$ is injective by \cite[2.1]{BlochKato}
or \cite[9.7.1]{GilleSz}.  In summary, we have: \emph{for $x_1,
\ldots, x_q \in \Fb^\times$, the following are equivalent:
\begin{enumerate}
\item \label{Kdiff.symb} The symbol $\{ x_1, \ldots, x_q \}$ is zero in $k_q(\Fb)$.
\item \label{Kdiff.der} $\der{x_1} \wedge \cdots \wedge \der{x_q}$ is zero in $\Omega^q$.
\item \label{Kdiff.pfree} The elements $\sqrt[p]{x_1}, \ldots, \sqrt[p]{x_q}$ are not $p$-free over $\Fb$.
\end{enumerate}}
\noindent The equivalence of \eqref{Kdiff.der} and
\eqref{Kdiff.pfree} is \cite[\S{V.13.2}, Th.~1]{Bou:alg2}.  For
further statements along these lines, see e.g.\ \cite[8.1]{Hoff:p}.

Given $y,x_1,\dots,x_q \in \Fb^\times$ with $\der x_1 \wedge
\cdots \wedge \der x_q \neq 0$ in $\Omega^q$, the
preceding equivalence implies that $y\,\der x_1 \wedge \cdots \wedge
\der x_q = 0$ in $\Omega^q/Z^q$ if and only if $y$ is a $p$-th
power in $\Fb(\sqrt[p]{x_1},\dots,\sqrt[p]{x_q})$.

%%%%%%%%%%%%%%%%%%%%%%%%%%%%%%%%%
\section{The Local Norm Theorem \ref{t.HENOTH} revisited} \label{s.LNTK}

We now prove an analogue of the Local Norm Theorem \ref{t.HENOTH}.
We continue the notation of \ref{p3.intro}, and focus on a nonzero
symbol $\gamma \in k_q(F)$ where $q = 1$, $q = 2$, or $p = 2$.  The
symbol $\gamma$ corresponds to a Galois field extension of $F$ of
degree $p$, a (central) associative division algebra of
dimension $p^2$ over $F$, or a $q$-Pfister quadratic form over $F$.
Write $V$ for the underlying vector space, which has dimension
$p^q$.  In each case, there is a canonical choice of homogeneous
polynomial $f \!: V \rightarrow F$ of degree $p$ and representing 1:
the norm, the reduced norm, or the quadratic form itself.  Further,
the valuation $\val$ extends to a valuation $\val_V$ on $V$ via the
formula $\val_V(v) := \val(f(v))/p$, and we write $\mfo_V^\times$
for the set of $v \in V$ with value zero.

For a given $a \in \mfo^\times$, we ask: How close is $f$ to
representing $a$?  Imitating the definition in \ref{ss.NE}, we
put\index{norm exponent (nexp)}
\[
\nexp_\gamma(a) := \sup \left\{ \val(a - f(v)) \mid v \in \mfo_V^\times \right\}.
\]
It is obviously equivalent to define $\nexp_\gamma(a)$ to be the
supremum of all $d \ge 0$
such that there exist $\beta \in \mfo$ and $v \in \mfo^\times_V$
such that $a = (1 - \pi^d \beta)\, f(v)$.

\begin{eg} \label{Ndepth}
As $f$ represents 1,
\[
\depth \{ a \} \le \nexp_\gamma(a) \quad (a \in \mfo^\times).
\]
Moreover, we have equality in case $p = 2$ and $q = 0$ (so
necessarily $f$ is the quadratic form $\qform{1}$ and $\gamma = 1
\in \Zm2 = k_0(F)$).
\end{eg}

\begin{LNT}\index{Local norm theorem!in terms of $K$-theory} \label{LNTK}
Suppose $p = 2$ or $1 \leq q \leq 2$.  For $a \in
\mfo^\times$ and a symbol $\gamma \in k_q(F)$, we have:
\begin{equation} \label{LNTK.1}
\nexp_\gamma(a) + \depth \gamma \le \depth(\{ a \} \cdot \gamma).
\end{equation}
Furthermore, the following are equivalent:
\begin{enumerate}
\item \label{LNTK.0} $\{ a \} \cdot \gamma = 0$.
\item \label{LNTK.rep} $a \in f(V)$.
\item \label{LNTK.big} $\nexp_\gamma(a) > m - \depth \gamma$.
\item \label{LNTK.infty} $\nexp_\gamma(a) = \infty$.
\end{enumerate}
\end{LNT}

The number $m$ was defined in \ref{p3.intro} to be $p \,\val(p) / (p
- 1)$. Note that the inequality (\ref{LNTK}.\ref{LNTK.1}) apparently
strengthens (\ref{filtration}.\ref{filt.1}) by Example \ref{Ndepth}.

\begin{proof}
\eqref{LNTK.0} and \eqref{LNTK.rep} are known to be equivalent:
for $q = 1$, it is \cite[4.7.5]{GilleSz} and for $p = 2$ it
is Prop.~\ref{p.CDPOQUA}.  For $q = 2$, the implication
\eqref{LNTK.rep} $\Rightarrow$ \eqref{LNTK.0} is elementary and the
converse is due to Merkurjev-Suslin \cite[12.2]{MS:Kcoh}.

Assume \eqref{LNTK.rep}, i.e., $f(v) = a$ for some $v \in V$.  Then
$\val_V(v) = \val(a)/p = 0$, so $v \in \mfo_V^\times$ and
\eqref{LNTK.infty} is obvious, hence also \eqref{LNTK.big}.

We now prove equation \eqref{LNTK.1}.  Suppose $a = (1 - \pi^d
\beta)\, f(v)$ for some $d \ge 0$, $\beta \in \mfo$ and
$v \in \mfo^\times_V$. Then
\[
\{ a \} \cdot \gamma =   \{ f(v) \} \cdot \gamma + \{ 1 - \pi^d \beta \} \cdot \gamma.
\]
But the first term on the right side is zero by the equivalence
of (i),(ii) already established. Hence
\[
\depth( \{ a \} \cdot \gamma ) = \depth (\{ 1 - \pi^d \beta \} \cdot \gamma) \ge d + \depth \gamma,
\]
and we have proved \eqref{LNTK.1}.

Finally, suppose \eqref{LNTK.big}.  By \eqref{LNTK.1}, the symbol
$\{ a \} \cdot \gamma$ has depth $> m$, hence the symbol
is zero, proving \eqref{LNTK.0}.
\end{proof}

In order to compare this result with the Local Norm Theorem
\ref{t.HENOTH}, we must relate $\wid$ with $\depth$; for the
purposes of this discussion, let us focus on the case $p = 2$ and
put $\wid(\gamma) := m - (\depth\gamma)/2$.  (This agrees with the
definition of $\wid$ for composition algebras  in case additionally
$q = 2$ or 3 by Cor.~\ref{texp.cor}(ii) below.)
Translating Equation \ref{LNTK}.\ref{LNTK.1} into this notation
gives:
\[
\nexp_\gamma(a) \le 2 \wid(\gamma) - 2 \wid(\{ a \} \cdot \gamma).
\]
That is, Theorem \ref{LNTK} sharpens Theorem \ref{t.HENOTH}.

\begin{rmk}
If every finite extension of $F$ has dimension a power of $p$ (``$F$
is $p$-special") for some prime $p$, then Theorem \ref{LNTK} holds for that prime $p$ and all $q \ge 1$ if
one adjusts slightly the statement of \eqref{LNTK.rep}. The adjusted
statement should be in terms of a norm variety for $\gamma$ as is
obvious from \cite[Prop.~2.4]{SusJ}; we leave the details to the
reader.  The paper \cite{SusJ} also provides the proof of the
equivalence of \eqref{LNTK.0} and the adjusted form of \eqref{LNTK.rep}.
\end{rmk}

%%%%%%%%%%%%%%%%%%%%%%%%%%%%%%%%%%%%%%%%%%%
\section{Gathering the depth} \label{gather.sec}

We maintain the notation of \ref{p3.intro}.
Recall
that $U^d k_q(F)$ is generated as an abelian group by $U_d \cdot
k_{q-1}(F)$.  In this section, we prove that a symbol in $U^d
k_q(F)$ can be written as $\{ u \} \cdot \alpha$ where $u \in U_d$
and $\alpha$ is a symbol in $k_{q-1}(F)$ (and not just as a sum of such things).  More precisely, we have:

\begin{gathlem} \label{gathlem} \index{Gathering Lemma}
For every nonzero symbol $\gamma \in k_q(F)$ with $q \ge 2$, there
is a $u \in \mfo^\times$ and a symbol $\alpha \in
k_{q-1}(F)$ such that $\gamma = \{ u \} \cdot \alpha$, $\depth \{ u
\} = \depth \gamma$, and $\depth \alpha = 0$.  The symbol $\alpha$
may be chosen to be $\{ a_2, \ldots, a_q \}$ with $0 \le \val(a_2) <
p$ and $\val(a_i) = 0$ for $3 \le i \le q$.  If $\depth \gamma$ is
not divisible by $p$, we may further arrange that $\val(a_2)$
has any pre-assigned value $j = 0, \dots, p - 1$ as desired.
\end{gathlem}

One should compare this lemma in the case $p = 2$ and $q = 2, 3$
with the Normal Form Theorems \ref{t.NOFOO} and \ref{t.NOFOQ}.
Heuristically speaking, here we ``gather the depth in the first
slot" $u$.  The Normal Form Theorems (in view of
Th.~\ref{texp.depth} below and with $u$ replaced by $L$) do the
same, except when $C$ is of unitary type, where they take $\depth L
= \depth C - 1$.

We first amass some preliminary results; the proof of the Gathering
Lemma will come at the end of the section. We use only background
material on $K$-theory including the material summarized in
\S\ref{filtration.sec}; we don't use anything else from this paper.

\subsection{} \label{BK.O}
We write $\ord(\pi^i)$ for an unspecified element (possibly zero) of
$\mfo$ divisible by $\pi^i$.  We have the trivial but useful
observation:
\[
u + \ord(\pi^j) = u(1 + \ord(\pi^j)) \qquad \text{for $u \in \mfo^\times$ and $j \ge 0$}.
\]
Indeed, for $b \in \mfo$, we have:
$u + b\pi^j = u\left(1 + (u^{-1}b)\pi^j\right)$.

\begin{eg}[cf.~\protect{\cite[p.~122]{BlochKato}}] \label{BK.eg}
Suppose that  $y \in \mfo$ has $\by = \bc^p$ for some $c \in \mfo$
and fix $0 \le s < \val(p) / (p - 1)$ such that $1 + y \pi^{ps}
\in \mfo^\times$.  Then
\[
1 + y \pi^{ps}  = 1 + c^p \pi^{ps} + \ord(\pi^{ps+1})
= (1 + c \pi^s)^p + \ord(\pi^{ps+1})
\]
where the second equality is because
\begin{equation} \label{BK.1}
\val\left(\pi^{is} \textstyle\binom{p}{i} \right) \ge is + \val(p) >
is + (p - 1) s \ge ps \quad \text{for $1 \le i < p$.}
\end{equation}
As $1 + y \pi^{ps}$ is a unit, so is $(1 + c\pi^s)^p$, and
\ref{BK.O} gives that $1 + y \pi^{ps}  = (1 + c\pi^s)^p\, (1 +
\ord(\pi^{ps+1}))$, hence $\{ 1 + y \pi^{ps} \} = \{ 1 +
\ord(\pi^{ps+1}) \}$ in $k_1(F)$. Noting that the hypotheses on $y$
obviously depend only on $\by$ (or applying \ref{BK.O} once more),
we find:
\[
\{ 1 + y\pi^{ps} + O(\pi^{ps+1}) \} = \{ 1 + O(\pi^{ps+1}) \}
\quad \text{in $k_1(F)$.}
\]
\end{eg}

\begin{lem} \label{pumping}
Let $a, b \in \mfo$, $i \geq 0$, $j \geq 1$ with $1 + a\pi^i
\in \mfo^\times$. Then in $k_q(F)$ we have:
\[
\{ 1 + a \pi^i, 1 + b \pi^j \} = \left\{ 1 + c
\pi^{i+j}, d \pi^{i(p - 1)} \right\}
\]
for some nonzero $c, d\in \mfo$.  Further, if $a, b \in
\mfo^\times$, then also $c, d \in \mfo^\times$.
\end{lem}

\begin{proof}
The computations in the proof of \cite[4.1.1b]{CT:kato} or
\cite[p.~122]{BlochKato} yield:
\[
\{ 1 + a \pi^i, 1 + b \pi^j \} = -\left\{1 + \frac{ab}{1 + a \pi^i}
\pi^{i+j}, -a\pi^i (1 + b \pi^j) \right\} \quad \in k_2(F).
\]
As $-1
= p - 1$ in $k_0(F)$, we have:
\[
\{ 1 + a \pi^i, 1 + b \pi^j \} = \left\{ 1 + \frac{ab}{1 + a\pi^i}
\pi^{i+j}, \left(-a(1+b\pi^j)\right)^{p-1} \pi^{i(p - 1)} \right\}.  \qedhere
\]
\end{proof}

\begin{lem} \label{gath.div}
Suppose that $\{ 1 + b\pi^{ps}, a_2, \ldots, a_q \}$ satisfies
$\der{\bb} \wedge \der{\ba_2} \wedge \der{\ba_3} \wedge \cdots
\wedge \der{\ba_q} = 0$ in $\Omega^q$, with $b, 1+ b\pi^{ps}, a_2,
\ldots, a_q \in \mfo^\times$ and $0 \le s <
\val(p)/(p-1)$. Then $\{ 1 + b\pi^{ps}, a_2, \ldots, a_q \}$ is
equal to $\{ u', a_2', \ldots, a_q' \}$ for some $u', a'_i \in
\mfo^\times$ where $\depth \{ u' \} > ps$.
%Let $\{ a_2, \ldots, a_q \}$ be a nonzero symbol in
%$k_{q-1}(F)$ such that $a_i \in \mfo^\times$, and let $u = 1 +
%b\pi^{ps} \in \mfo^\times$ for some $b \in \mfo^\times$ and $ps < m$.  If
%$\der{\bb} \wedge \der{\ba_2} \wedge \der{\ba_3} \wedge \cdots
%\wedge \der{\ba_q} = 0$ in $\Omega^q$, then $\{ u, a_2, \ldots, a_q
%\}$ is equal to $\{ u', a_2', \ldots, a_q' \}$ where $\depth \{ u'
%\} > ps$ and $a'_i \in \mfo^\times$.
\end{lem}

\begin{proof}
We may assume that $\der{\bb}$ is not zero---hence that $\bb$ is not
a $p$-th power---by Example \ref{BK.eg}.  This settles the $q = 1$
case.

Suppose $q \ge 2$ and $\{ \ba_2, \ldots, \ba_q \}$ is zero in
$k_{q-1}(\Fb)$, so $\der{\ba_2} \wedge \cdots \wedge \der{\ba_q} = 0$ by
\ref{Kdiff}.  We apply the $q-1$ case of the lemma with $u = a_q$
and $s = 0$ (so $u = a_q = 1 + c$ with $c = a_q - 1$, hence
$\der{\bc} = \der{\ba_q}$) to see that $\{ a_2, \ldots, a_q
\} = \{ a'_2, \ldots, a'_q\}$ where $\depth \{ a'_q \}$ is positive.
Then Lemma \ref{pumping} gives the claim.

\newcommand{\tu}[1]{{[#1]}}
So we may assume that $\{ \ba_2, \ldots, \ba_q \}$ is not zero.
Write $\tu{i}$ for a $(q-1)$-tuple $(i_2, \ldots, i_q)$ with $0 \le
i_j < p$ and put $\ba^{\tu{i}}$ for $\ba_2^{i_2} \ba_3^{i_3} \cdots
\ba_q^{i_q} \in \Fb$.  By \ref{Kdiff} there are $c_{\tu{i}} \in
\mfo$ such that
\[
\bb = \left( \sum\nolimits_\tu{i} \bc_\tu{i} \sqrt[p]{\ba^{\tu{i}}}
\right)^p = \sum\nolimits_\tu{i} \bc^p_\tu{i} \ba^\tu{i}.
\]
If it happens that $\bc_\tu{i} = 0$ for all nonzero $\tu{i}$, then
$\bb$ is a $p$-th power in $\Fb$ and we are done.  We now show,
roughly speaking, that we can make $\bc_\tu{i}$ zero for all nonzero
$\tu{i}$.

More precisely, fix a nonzero $\tu{i}$ with $\bc_\tu{i}$ nonzero;
choose a specific $j_0$ such that $i_{j_0} \ne 0$. Take $E/F$ to be
the extension obtained by adjoining a $p$-th root $\alpha$ of
$a_{j_0}^{i_{j_0}} \prod_{j \ne j_0} (-a_j)^{i_j}$; obviously
$\alpha$ is integral.  Take $v := c_\tu{i} \alpha$.  As $\bc_\tu{i}$
is not zero, the residue of $v$ does not belong to $\Fb$, hence
$v\pi^s$ is not in $F$ and so has minimal polynomial $x^p -
(c_\tu{i} \alpha \pi^{s})^p$ in $F[x]$.  By degree count, this is
also the characteristic polynomial $\chpoly_{v\pi^s}(x)$ of $v\pi^s$
as an element of the $F$-algebra $E$.  It follows that
\begin{equation} \label{brussel.1}
N_{E/F}(1 - v\pi^s) = \chpoly_{v\pi^s}(1) = 1 + (-1)^p N_{E/F}(v) \pi^{ps}.
\end{equation}

Next observe that $E$ kills $\gamma := \{ a_2, \ldots, a_q \}$:
we renumber the $a_j$'s so that $j_0 = 2$.  In $k_{q-1}(E)$,
we have:
\[
i_2 \gamma = \{ a_2^{i_2}, a_3, \ldots, a_q \} =
\sum_{j=3}^q (p-i_j) \{ -a_j, a_3, \ldots, a_q \} = 0,
\]
where the middle equality is because $\alpha$ is in $E$.
As $i_2$ is not divisible by $p$, we deduce that $\gamma$ is zero
in $k_{q-1}(E)$, as required.

Now the projection formula \cite[7.2.7]{GilleSz} gives:
\[
0 = N_{E/F}\left(  \{ 1 - v \pi^s \} \cdot \gamma \right)  =
\{ N_{E/F}(1 - v \pi^s) \} \cdot \gamma
\]
in $k_q(F)$.  Combining this with \eqref{brussel.1}, we find:
\[
\{ 1 + b \pi^{ps} \} \cdot \gamma = \{ 1 + (b + (-1)^p N_{E/F}(v))
\pi^{ps} + O(\pi^{ps+1}) \} \cdot \gamma.
\]
But the residue of $N_{E/F}(v)$ is $\bv^p = \bc_\tu{i}^p
\ba^\tu{i}$. In this way, we have replaced $b$ with a new one that
has fewer nonzero coefficients $\bc_\tu{i}$.  Repeating this process
completes the proof.
\end{proof}

\begin{proof}[Proof of the Gathering Lemma \ref{gathlem}]
First, consider a symbol $\{ x, y \} \in k_2(F)$.
Suppose that neither $\lambda(x)$ nor $\lambda(y)$ are
divisible by $p$. Then there is some $s$ such that $\val(x) \equiv s
\val(y) \pmod{p}$ and $\{ x, y \} = \{ x (-y)^{-s}, y \}$ because
$\{ -y, y \}$ is zero in $k_2(F)$.  As $\val(x (-y)^{-s}) = \val(x)
- s \val(y)$, we may assume that $x$ has value 0.

Second, we may shuffle the entries in a symbol $\{ x_1, \ldots, x_q
\}$ by a permutation $\s$.  We have $\{ x_1, \ldots, x_{q-1}, x_q \}
= \{ x_{\s(1)}, \ldots, x_{\s(q-1)}, x_{\s(q)}^{\sgn \s} \}$.

Combining the two preceding paragraphs shows that we may write
$\gamma = \{ u \} \cdot \alpha$ for
$\alpha = \{ a_2, \ldots, a_q \}$ with $u, a_i \in \mfo^\times$ for $3 \le i \le q$.
Amongst all such ways of writing $\gamma$,
fix one
with $\depth \{ u \}$ maximal.  For sake of contradiction, suppose
that $d := \depth \{ u \} < \depth \gamma \le m$.
Put $r = 2$ if $\val(a_2) = 0$ and $r = 3$ if $0< \val(a_2)
< p$.

We now inspect the Kato isomorphism at depth $d$.
By hypothesis, $\gamma$ is zero in $\gr^d k_q(F)$, so has zero
image.  If $d = 0$, then $\{ \bu, \ba_r, \ldots, \ba_q \}$ is zero
in $k_*(\Fb)$, hence $\der{(\bu - 1_{\Fb})} \wedge \der{\ba_r} \wedge \cdots
\wedge \der{\ba_q}$ is zero by \ref{Kdiff}.    Lemma \ref{gath.div}
gives a contradiction. The case where $d = ps$ for some $0 < s =
d/p < (\depth\gamma)/p < \val(p) / (p-1)$ is similar.

Finally we suppose $d$ is not divisible by $p$. We claim that
$\val(a_2)$ may be freely chosen.  Indeed, for $j \in \Z$,
$\val(a_2) - j \equiv sd \pmod{p}$ for some $s \ge 0$. Writing $u =
1 + b \pi^d$ for $b \in \mfo^\times$, we have:
\[
\{ 1+ b\pi^d, a_2 \} = \{ 1 + b\pi^d, a_2 \} - s \{ 1 + b\pi^d,
-b\pi^d \} = \{ 1 + b\pi^d, a_2 (-b)^{-s} \pi^{-sd} \},
\]
where $\val(a_2 (-b)^{-s} \pi^{-sd}) \equiv j \pmod{p}$, proving the
claim. The hypothesis that $\{ u \} \cdot \alpha$ has depth greater
than $d$ implies that $\der{\ba_2} \wedge \cdots \wedge \der{\ba_q}$
is zero, and applying Lemma \ref{gath.div} to $\alpha$ with $s = 0$
shows that we may assume that one of the $a_i$ has residue 1; Lemma
\ref{pumping} gives a contradiction.
\end{proof}

\begin{rmk} \label{preprmk}
We can quickly deduce a useful restatement of the Gathering Lemma.
Fix some $r \ge 1$ and a permutation $\s$ of $\{ 2, \ldots, q \}$
and put
\[
\beta := \{ u, a_{\s(2)}, \ldots, a_{\s(r)} \} \eand \delta := \{
a_{\s(r+1)}, a_{\s(r+2)}, \ldots, a_{\s(q)} \}.
\]
Because of the identity $\{ x, y \} = -\{y, x \}$ in $K^M_2(F)$, we
find:
\[
\gamma = \pm \beta \cdot \delta, \quad \depth \beta = \depth \gamma,
\eand \depth \delta = 0.
\]
Indeed, for the second and third equalities, $\ge$ is obvious and
$\le$ follows from the first equality and equation
(\ref{filtration}.\ref{filt.1}). We will apply this in the case $p =
2$, so the sign in the first equation will be irrelevant.
\end{rmk}

%\begin{rmk}
%The proofs in the literature that the Kato isomorphisms at depth $d$ from \ref{Kato.iso} are isomorphisms is done by proving that the inverse function $\overline{\rho_d}$ is an isomorphism at depth $d$---see \cite[p.~125]{BlochKato} for a proof that $\overline{\rho_d}$ is injective for $d < m$.  That result does not have any restriction on $q$ and $p$, but one cannot deduce the Gathering Lemma from it.  After all, the statement that $\depth \gamma > d$
%\end{rmk}

%%%%%%%%%%%%%%%%%%%%%%%%%%%%%%%%%%%%%%%%%%%
\section{Ramification index for symbols}

\begin{defn} \label{unram.def}
For a class $\gamma \in k_q(F)$, we put $e_\gamma = p$ if\index{$e_\gamma$}
\begin{itemize}
\item $\depth(\gamma)$ is \emph{not} divisible by $p$, or
\item $\depth(\gamma)$ is divisible by $p$ and its initial form
has nonzero projection in the second summand of $\gr^d k_q(F)$.
(Note that this condition does not depend on the choice of
uniformizer $\pi$.)
\end{itemize}
Otherwise---or if $\gamma = 0$---we put $e_\gamma = 1$.
\end{defn}

The ``ramification index" $e_\gamma$ is more subtle than in the
case of good residue characteristic, see Example \ref{ram.eg}
below.
But we do have the following positive results:

\begin{prop} \label{unramsubsym}
Let $\gamma \in k_q(F)$, $q \geq 2$, be a non-zero symbol. Then
there exist a symbol $\beta \in k_{q-1}(F)$ and an element $a \in
F^\times$ such that $\gamma = \beta\cdot\{a\}$, $e_\beta = 1$ and
one of the following holds.
\begin{itemize}
\item [(i)] $\depth\beta = 0$, $\depth\{a\} = \depth\gamma$, $a \in
\mfo^\times$ and
\[
e_\gamma = 1 \Longleftrightarrow \depth \gamma \equiv 0 \bmod p.
\]

\item [(ii)] $\depth\beta = \depth\gamma$, $0 < \lambda(a) <p$ and
$e_\gamma = p$.
\end{itemize}
\end{prop}

\begin{proof} Write $\gamma$ as in the Gathering Lemma~\ref{gathlem}, where
we may assume $a_2,\dots,a_{q-1} \in \mfo^\times$, $0 \leq
\lambda(a_q) < p$. If $\lambda(a_q) = 0$, then the Kato isomorphism
at depth zero shows $e_\alpha = 1$, so with $\beta := \alpha$ and $a
= u^{\pm 1}$ we are in Case~(i) since the Kato isomorphism at depth
$d := \depth\gamma$ gives $e_\gamma = 1$ iff $d$ is divisible by
$p$. Now suppose $\lambda(a_q) > 0$; by the final statement of the
Gathering Lemma, we may also assume $d \equiv 0 \bmod p$. Arguing as
before, in particular consulting the Kato isomorphisms again for the
determination of $e_\gamma$, $\beta := \{u,a_2,\dots,a_{q-1}\}$ has
depth $d$ and $e_\beta = 1$, so we are in Case~(ii).
\end{proof}

\begin{prop} \label{ramified}
For a symbol $\gamma \in k_{q}(F)$ for $q \ge 2$, the following
are equivalent:
\begin{enumerate}
\item $\gamma = \beta \cdot \{ a \}$ for a nonzero symbol $\beta \in k_{q-1}(F)$
with $e_\beta = 1$ and $a \in \Fx$ of value not divisible by $p$.
\item $e_\gamma = p$ and $\depth\gamma$ is divisible by $p$.
\end{enumerate}
If these conditions hold, then additionally $\depth\gamma =
\depth\beta$.
\end{prop}

\begin{proof}
The proof amounts to looking at the explicit formulas for the Kato isomorphisms.

\medskip
\noindent{\underline{(i) $\Implies$ (ii)}}: The Kato isomorphisms
send $\beta$ to a nonzero symbol $\bbeta$ (in some cohomology group
over $\Fb$), because $e_\beta = 1$.  We have $\depth(\gamma) \ge
\depth(\beta)$ by (\ref{filtration}.\ref{filt.1}), and an
examination of the isomorphisms show that the isomorphism at the
depth of $\beta$ (divisible by $p$) sends $\beta \cdot \{ a \}$ to
$\val(a) \bbeta$ in the second component of the image, which is not
zero because all of the targets of the Kato isomorphisms are abelian
groups killed by $p$; that is, $e_\gamma = p$ and
$\depth\gamma = \depth\beta$.

%We divide into cases based on the depth of $\beta$; it is divisible
%y $p$ because $e_\beta = 1$.  In each case, we content ourselves
%with describing the form of $\beta$ and $\bbeta$.

%If $\beta$ has depth 0, then $\beta$ is $\{ x_1, \ldots, x_q \}$ for
%$x_i$ of value 0 such that $\bbeta = \{ \bx_1, \ldots, \bx_q \}$
%is not zero in $k_q(\Fb)$.

%If $\beta$ has depth $d$ such that $0 < d < m$, then $\beta$ is
%$\{ 1 + x\pi^d, y_1, \ldots, y_{q-1}\}$ for $x, y_i$ of value 0
%such that $\bbeta$ is $\bx \dlog{\by_1} \wedge \cdots \wedge
%\dlog{\by_{q-1}}$ in $\Omega^{q-1}/Z^{q-1}$.

%If $\beta$ has depth $m$, then $\beta$ is $\{ 1 + x (\zeta - 1)^p, y_1, \ldots, y_{q-1}\}$ for $x, y_i$ of value 0 such that $\bbeta$ is $\bx \dlog{\by_1} \wedge \cdots \dlog{\by_{q-1}}$ in $H^1(\Fb, \nu(q - 1))$.

\medskip
\noindent{\underline{(ii) $\Implies$ (i)}}: By (ii),
alternative (i) of Prop.~\ref{unramsubsym} does not hold. Hence
alternative (ii) does.
\end{proof}
%Write $\gamma$ as in the Gathering Lemma \ref{gathlem}.  The
%hypothesis $e_\gamma = p$ implies---as can be seen by inspecting the
%Kato isomorphism at each depth---that $\val(a_2)$ is not divisible
%by $p$.  We can write $\gamma$ as in Remark \ref{preprmk}, where we
%take $\delta$ to consist of just $\{ a_2 \}$.  Inverting one of the
%slots in $\alpha$ if necessary, we can arrange for $\gamma = \alpha
%\cdot \{ a_2 \} $, as desired.

%%%%%%%%%%%%%%%%%%%%%%%%%%%%%%%%%%%%%%%%%%%%%%%%%
\section{Theorem \ref{t.CHAGORE} revisited} \label{K.CHAGORE.sec}

We will now prove a version of Theorem \ref{t.CHAGORE} for
$K$-theory mod-2 when $\Fb$ has characteristic 2; in the notation of
\S\ref{filtration.sec} we restrict to the case $p = 2$.  We replace
the anisotropic round quadratic form $P$ with $e_{P/F} = 1$ from
Th.~\ref{t.CHAGORE} with a nonzero symbol $\gamma \in k_q(F)$ with
$q \ge 1$ and $e_\gamma = 1$ (in particular, $\depth \gamma$ is
even).  We replace the hypothesis on $\texp$ with the hypothesis
$\depth \{ \mu \} + \depth \gamma = m$ (recall that $m = 2\val(2)$),
the largest possible depth for a nonzero element of $k_q(F)$.  The
extreme cases where $\depth \gamma$ is $0$ or $m$ are comparatively
easy, so we focus on the middle case.  We will use the following:

\subsection{Technique} \label{K.CHAGORE.tech} If one has a \emph{nonzero}
class $\bb \dlog{\ba_2} \wedge \cdots \wedge \dlog{\ba_q} \in
\Omega^{q-1}/Z^{q-1}$, we may apply $\der{}$ and obtain the nonzero
symbol $x_0 := \der{\bb} \wedge \dlog{\ba_2} \wedge \cdots \wedge
\dlog{\ba_q}$ in $\Omega^q$.  The $\Fb$-span of this  symbol
contains $x := \bb^{-1} x_0$, which lies in $\nu(q)$.  Indeed, it is
the unique nonzero element of $\Fb x_0$ with this property, because
$\nu(q)$ is defined to be $\ker(\gamma - 1)$ for $\gamma$ the
inverse Cartier operator \cite[pp.~123, 124]{CT:kato} and for $c \in
\Fb$ we  have $(\gamma - 1)(c x) = (c^2 - c) x \in \Omega^q$.  A
canonical isomorphism identifies $x$ with the class of the
anisotropic symmetric bilinear form $B = \pform{\bb, \ba_2, \ldots,
\ba_q}$ in the graded Witt ring \cite{Kato:Milnor}, hence with $B$
itself \cite[6.20]{MR2427530}.  The equivalence from
Prop.~\ref{t.CHARPFIBIL} takes $B$ and gives an extension $K/\Fb$ of
degree $2^q$ with a unital linear form $s \!: K \to \Fb$.

Let $\gamma \in k_q(F)$ be a nonzero symbol of even depth $d \ne 0, m$
and suppose that $e_\gamma = 1$.  Then the initial form of $\gamma$
is a nonzero symbol in $\Omega^{q-1}/Z^{q-1}$ and the technique in
the preceding paragraph gives a $(K, s)$ derived from $\gamma$.

\begin{prop} \label{K.CHAGORE}
Let $\gamma \in k_q(F)$ be a nonzero symbol of even depth $d \ne 0,
m$ with $e_\gamma = 1$. Write $(K, s)$ as in \ref{K.CHAGORE.tech},
and write $\gamma = \{ 1 + x \pi^d \} \cdot \alpha$ with $x \in
\mfo^\times$ as in the Gathering Lemma \ref{gathlem}. For $\mu = 1 -
b \pi^{m-d}$ with $b \in \mfo^\times$, we have: $\{ \mu \} \cdot
\gamma = 0$ in $k_{q+1}(F)$ if and only if the residue of $x b
\pi^m/4$ is in the image of $\wp_{K,s}$.
\end{prop}

\begin{proof}
Write $\alpha = \{ a_2, \ldots, a_q \}$ for some $a_2, \ldots, a_q
\in \mfo^\times$.  Using Lemma \ref{pumping}, we calculate:
\[
\{ \mu \} \cdot \{ 1 + x \pi^d \} = \{ 1 + x \pi^d, 1 - b \pi^{m-d}
\} = \left\{ 1 - \frac{bx}{1 + x \pi^d} \pi^m, -x (1 - b \pi^{m-d})
\right\}.
\]
We see from this that the initial form of $\{ \mu \} \cdot \gamma$
is $\bx \bb \be \dlog{\bx} \wedge \dlog{\ba_2} \wedge \cdots \wedge \dlog{\ba_q}$ where
we have set $\varepsilon := \pi^m / 4$ to simplify the notation.
This determines the class of the quadratic Pfister form
$\qpform{\bx, \ba_2, \ldots, \ba_q, \bx \bb \be}$ in the graded Witt
group of quadratic forms \cite{Kato:Milnor}.
The Arason-Pfister Hauptsatz \cite[23.7(1)]{MR2427530} implies that
this class is zero (equivalently, $\{ \mu \} \cdot \gamma$ is zero)
if and only if the Pfister form is isotropic, if and only if $\bx \bb  \be$ is in the
image of $\wp_{K,s}$ by Cor.~\ref{c.ISAS}(a), proving the claim.
\end{proof}

%%%%%%%%%%%%%%%%%%%%%%%%%%%%%%%%%%%%%%%%%%%%%%%%%
\section{Dictionary between $K$-theory and algebras and quadratic forms} \label{dict.sec}

In the cases $q = 1$, $q = 2$, or $p = 2$, we have a close
relationship between properties of symbols in $k_q(F)$ relative to
Kato's filtration and valuation-theoretic properties on the
corresponding algebras.  Specifically:

\begin{prop} \label{dict}
In cases $q = 1$ and $q = 2$
we have:\index{ramfication,
relationship between $K$-theory and valuations}
\begin{enumerate}
\item $e_\gamma = e_{D/F}$.
\item $\Db$ is a separable division algebra over $\Fb$ and distinct
from $\Fb$ if and only if $\depth(\gamma) = m$.
\end{enumerate}
\end{prop}

That is, you can determine $e_{D/F}$ and whether or not $D$ is tame
by examining the corresponding symbol in $k_q(F)$.

We prove the proposition in \S\ref{pf.dict.1.sec}  and
\ref{pf.dict.2} below.  The calculations appearing in the proof of
this result are similar to some in sections 1.1 and 2.1 of
\cite{Tignol:wildconf}.  Roughly speaking, our statements here differ
from those in \cite{Tignol:wildconf} by using the language of Kato's filtration.

In contrast with Proposition \ref{dict}, the following theorem
relies heavily on the results of Part \ref{hensel.part}. It is
proved in \S\ref{pf.texp.depth}.\index{trace exponent ($\texp$)!in
terms of $K$-theory} \index{Tignol's invariant $\wid$!in terms of
$K$-theory}

\begin{thm} \label{texp.depth}
Let $p = 2$.  Fix a nonzero symbol $\gamma \in k_q(F)$ for some $q
\ge 1$ and let $Q$ be the corresponding anisotropic $q$-Pfister
form.  Then:
\begin{enumerate}
\item \label{texp.e} $e_\gamma = e_{Q/F}$.
\item \label{texp.wild} $\overline{Q}$ is nonsingular (i.e., $Q$
is tame) if and only if $\depth(\gamma) = 2\val(2)$.
\item \label{texp.depth.eqn} $\texp(Q) = \val(2) - \floor{\frac{\depth(\gamma)}2}$.
\end{enumerate}
\end{thm}

In the statement of the theorem, $Q$ is a quadratic form, whereas
the results in \S\S\ref{s.COALHE}--\ref{s.VADAEN}  (including the
definition of tame and $\tex$) are in terms of pointed quadratic
spaces.  This is harmless in view of Prop.~\ref{p.PONO}.  The
theorem leads to:

\begin{cor} \label{texp.cor}
Fix a division quaternion or octonion algebra $C$ over $F$ and write
$\gamma$ for the corresponding symbol in $k_*(F)$.  We have:
\begin{enumerate}
\item $C$ is of primary type if and only if $e_\gamma = 2$ and
$\depth(\gamma)$ is even.  (Otherwise $C$ is of unitary type.)
\item $\wid(C) = \val(2) - \frac{\depth(\gamma)}2$
\end{enumerate}
\end{cor}

That is, one can read off properties of the composition algebra
$C$---including the invariant $\omega(C)$ studied by Saltman and
Tignol---from the properties of the corresponding symbol $\gamma$ in
Milnor $K$-theory.  We prove Theorem \ref{texp.depth} and
Cor.~\ref{texp.cor} in \S\ref{pf.texp.depth}.

\begin{rmk}
It is natural to wonder if on can extend Prop.~\ref{dict} to include
the case $p = q = 3$, where the corresponding algebraic objects are
Albert (Jordan) algebras obtained from the first Tits construction.
The answer is no, because we do not know if two such Albert algebras
corresponding to the same symbol in $k_3(F)$ are necessarily
isomorphic.  (This is a special case of open problem \#4 from
\cite{PR}.)   If this is indeed the case, then one can easily extend
Prop.~\ref{dict} to include the case $p = q = 3$ by taking
advantage of the valuation theoretic results in \cite{MR50:422}
and by imitating the proof in the cases given below.
\end{rmk}

%%%%%%%%%%%%%%%%%%%%%%%%%%%%%%%%%%%%%%%%%%%%%%%%%
\section{Proof of Proposition \ref{dict}: case $q = 1$} \label{pf.dict.1.sec}

\begin{eg}[$q = 1$ and nonzero depth divisible by $p$] \label{deg1.even}
Fix $\{ x \} \in k_1(F)$ of depth $d$ divisible by $p$, so $x = 1
+ u\pi^d$ for some $u$ of value 0. As the ``second summand"
$\Omega^{-1}/Z^{-1}$ or $H^1(\Fb, \nu(-1))$ in the Kato isomorphism
is zero, $e_{\{ x\}} = 1$.

The algebra $D$ corresponding to $\{ x \}$ is $F(\chi)$ where $\chi^p = x$.
For
\[
\alpha := \pi^{-d/p} (\chi -1) \quad \in D
\]
we have
\[
(\alpha + \pi^{-d/p})^p = \pi^{-d} x = u + \pi^{-d},
\]
so
\begin{equation} \label{y.eqn}
\sum_{i=1}^p \binom{p}{i} \alpha^i \pi^{-(p - i)d/p} - u = 0.
\end{equation}
For $1 \le i < p$, the prime $p$ divides $\binom{p}{i}$, so
\[
\val\left( \textstyle\binom{p}{i} \pi^{-(p-i)d/p} \right) \ge
\val(p) - \frac{d}{p}(p-i) = \frac{(m-d)(p-1)+d(i-1)}p.
\]
As $i \ge 1$ and $d \le m$, this is at least 0, hence $\alpha$ is
integral.  Further, the coefficient of $\alpha^i$ in \eqref{y.eqn}
has residue zero for $2 \le i < p$ in all cases and also for $i =
1$ if $d < m$.  (Clearly, $\balpha$ is not zero, so $\val(\alpha)
= 0$.)

Therefore, if $d < m $, $\Db$ contains $\balpha$ satisfying
$\balpha^p - \bu = 0$.  As $x$ has depth $d$, the Kato isomorphism
shows that $\bu$ is not a $p$-th power in $\Fb$, and we conclude
that $\Db$ is the proper extension $\Fb(\sqrt[p]{\bu})$ and
$e_{D/F} = 1$.

In case $d = m$, we set $\eta := \zeta - 1$.  As $\val(\eta) = m/p
= \val(p)/(p-1)$, we have
\[
x = 1 + b \eta^p \quad \text{for}\ b = \frac{u\pi^m}{\eta^p}.
\]
We put $\beta := \eta^{-1} (\chi - 1)$ and apply the same
reasoning as in the case $d < m$ with $\pi^{-d/p}$ replaced with
$\eta^{-1}$.  The element $\beta$ satisfies
\begin{equation} \label{deg1.1}
\sum_{i=1}^p \binom{p}{i} \beta^i \eta^{-(p-i)} - b = 0.
\end{equation}
Again, the coefficients of $\beta$ are integral because
\begin{equation} \label{deg1.2}
\val \left( \binom{p}{i} \eta^{-(p-i)} \right) \ge \val(p) - (p-i)
\frac{\val(p)}{p-1} = \val(p) \frac{i-1}{p-1} \ge 0 \quad
\text{for $1 \le i < p$}.
\end{equation}
Taking residues of \eqref{deg1.1} kills the terms with $1 < i <
p$. For the $i = 1$ term, we note that expanding the equation $(1
+ \eta)^p = 1$ gives $\sum_{i=1}^p \binom{p}{i} \eta^i = 0$, so
$p\eta = -\sum_{i=2}^p \eta^i$ and
\[
p\eta^{-(p-1)} = p\eta^{1-p} = -\sum_{i=2}^p \binom{p}{i} \eta^{-(p-i)}.
\]
Taking residues and applying \eqref{deg1.2}, only the $i = p$ term
is nonzero on the right side, so $p\eta^{-(p-1)}$ has residue
$-1$.  Taking residues of \eqref{deg1.1}, we find the equation
$\bbeta^p - \bbeta - \bb = 0$ in $\Db$.  The fact that $x$ has
depth $d$ asserts that the element $\bb$ is nonzero in $H^1(\Fb,
\nu(0)) \cong \Fb / \wp(\Fb)$, i.e., $\bbeta$ generates a proper
extension of $\Fb$ and we conclude that $D/F$ is unramified.  We
have verified Proposition \ref{dict} for $\gamma \in k_1(F)$ of
nonzero depth divisible by $p$.
\end{eg}

\begin{rmk}
It might be illuminating to compare Example \ref{deg1.even} in the
case $p = 2$ with the material in Part \ref{hensel.part}.    In case
$0 < d < m = 2\val(2) = 2\texp(F)$, we compare the example to
Prop.~\ref{p.CDUNEV} with $P = F$. The Kato isomorphism at depth $d$
sends $\{x\}$ to the image of $\bar u$ in $\bar F/\bar F^2$, and we
conclude $\bar u \notin \bar F^2$, in agreement with
Prop.~\ref{p.CDUNEV}. Moreover, identifying $F(\sqrt{x}) = \mbC(F,x)
= F \oplus Fj$ and $j = -\chi$, we obtain that $-\alpha = \Xi$ in the
sense of (\ref{p.CDUNEV}.\ref{CDPREV}); it is a normalized trace generator.

The situation with $d = 2\val(2)$ is slightly more delicate; we compare the example to
Thm.~\ref{t.CHAGORE} (resp.\ Cor.~\ref{c.CONTAM}) for
$P$ (resp.\ $C$) $= F$. First of all, we have
$\zeta = -1$, $\eta = -2$, $\beta = \frac{\chi - 1}{2}$, and $w_0 =
\frac{\pi^{\lambda(2)}}{2} \in \mfo^\times$ is the unique normalized
trace generator of $F$. Moreover, $\bar s_{w_0} = \Eins_{\bar F}$,
forcing $\wp_{\bar F,\bar s_{w_0}} = \wp$ to be the usual
Artin-Schreier map on $\bar F$. We have $u = -\beta$, $u_0 =
n_F(w_0)u = -\beta_0$ in the sense of (\ref{t.CHAGORE}.\ref{MUTEX})
and $b = u_0$. The Kato isomorphism at depth $m$ sends $\{x\}$ to
$\bar u_0 = \bar\beta_0 \in H^1(\bar F,\nu(0)) = \bar F/\wp(\bar
F)$, forcing $\bar\beta_0 \notin \Im(\wp)$ in agreement with the
equivalence (iii) $\Leftrightarrow$ (v) in Thm.~\ref{t.CHAGORE}.
Furthermore, by Cor.~\ref{c.CONTAM},
\[
\bar C = \bar F[\bft]/(\bft^2 - \bft - \bar u_0) = \bar
F[\bft]/(\bft^2 - \bft - \bar b) = \overline{F(\sqrt{x})},
\]
in agreement with the second part of Example~\ref{deg1.even}.
\end{rmk}

\subsection{} \label{pf.dict.1} \emph{Proof of Proposition
\ref{dict} for $q = 1$}: Fix a nonzero $\{ x \} \in k_1(F)$ and
write $D = F(\chi)$ as in Example \ref{deg1.even}.

Suppose first that $x$ has depth 0; we may assume that $0 \le
\val(x) < p$.  The Kato isomorphism
\[
\gr^0 k_1(F) \iso k_1(\Fb) \oplus k_0(\Fb) \cong {\Fb^\times} /
{\Fb^{\times p}} \oplus \Zm{p}
\]
is the direct sum of the specialization map $s_\pi$ and the tame
symbol $\partial$, described concretely in \cite{GilleSz}, e.g. As
$\partial(x) \equiv \val(x) \pmod{p}$,
we have $e_{\{ x \}} = p$ if and only if $\val(x)$ is not zero. As
to the extension $D/F$, the element $\chi$ satisfies $\chi^p - x =
0$, hence is integral.  If $\val(x)$ is not zero, then $\bchi^p$
is zero, hence $\bchi$ is zero, $\Db = \Fb$, and $e_{D/F} = p$.
If $\val(x)$ is zero, then since $x$ has depth 0, $\bx$ is not a
$p$-th power in $\Fb$ and $\Fb(\bchi)$ is a proper extension of
$\Fb$; in this case $\Db = \Fb(\sqrt[p]{\bx})$ and $e_{D/F} = 1$.

Suppose that $x$ has depth $d$ not divisible by $p$; we may take
$x$ of the form $1 + u \pi^d$ where $u$ has value 0.  Then
\[
N_{D/F}(\chi - 1) = (-1)^{p-1}(x - 1) = (-1)^{p-1} u \pi^d,
\]
so the value of $\chi - 1$ is not an integer and $e_{D/F} = p =
e_{\{x\}}$.

Note that in both of these cases, $\Db$ is not a proper separable
extension of $\Fb$. The remaining case where $d$ has nonzero depth
divisible by $p$ was treated in Example \ref{deg1.even}, which
concludes the proof.$\hfill\qed$

\begin{eg} \label{ram.eg}
Let $\Fb := \F_2((x))$, Laurent series over the field with 2
elements. Construct a field $F$ of characteristic zero by taking the
absolutely unramified Teichm\"uller extension of $\Fb$ as in Example
\ref{e.LAU} and adjoining $\sqrt{2}$. Then $F$ has residue field
$\Fb$ and $\sqrt{2}$ is a uniformizer. Consider the elements $a = 1+
\sqrt{2}^3$ and $b = 1 + 2u$ in $F$, where $u \in F$ is a unit with
residue $x$.  The symbols $\{a\}$ and $\{b\}$ have depth 3 and 2.
The valuation $\val$ ramifies on $F(\sqrt{a})$ and does not ramify
on $F(\sqrt{b})$.

A reader familiar with the situation of good residue characteristic
might assume that the quaternion algebra $(a, b)$ is division over
$F$, but in fact it is split.  This is easily seen because the
symbol $\{ a,  b\}$ has depth at least 5 by
(\ref{filtration}.\ref{filt.1}), which is greater than $2
\val(2)$.
\end{eg}

%%%%%%%%%%%%%%%%%%%%%%%%%%%%%%%%%%%%%%%%%%%%%%%%%
\section{Proof of Proposition \ref{dict}: cases $q = 2$} \label{deg2.sec}

We now prove Proposition \ref{dict} and Theorem \ref{texp.depth}
in the case $q = 2$.
We consider a nonzero symbol $\{ x, y \}$ in $k_2(F)$ and write
$D$ for the corresponding symbol algebra.

\begin{lem}
$e_{D/F} = p$ if and only if some subfield $L$ of $D$ has $e_{L/F} = p$.
\end{lem}

\begin{proof}
Easy.  If such an $L$ exists, then there is some $\ell \in \Lx$
such that $\val(N_{L/F}(\ell))$ is not divisible by $p$.  But
$N_{L/F}(\ell) = \Nrd_{D/F}(\ell)$ by \cite[p.~150]{Draxl}.  This
proves ``if".  For ``only if", note that a nonzero element of $D$
generates a subfield $L/F$.
\end{proof}

\subsection{} \label{pf.dict.2} \emph{Proof of Proposition \ref{dict}
for $q = 2$}: Suppose first that the depth $d$ of $\{ x, y \}$ is
not divisible by $p$ (so $e_{\{x,y\}} = p$).

By the Gathering Lemma \ref{gathlem} we can arrange that $x$ is
in $U_d$ and $y \in \Ox$ has residue that is not a $p$-th power.
As in \ref{pf.dict.1}, $e_{F(\chi)/F} = p$, hence $e_{D/F} = p$
and the dimension of $\Db/\Fb$ is $p$.  The residue field of
$F(\psi)$ is the proper extension $\Fb(\sqrt[p]{\by})$ of $\Fb$,
and by dimension count it is all of $\Db$ and $D$ is wild.  So we
may assume that the depth of $\{ x, y \}$ is divisible by $p$.

Suppose now that $e_{\{x,y\}} = p$.  By Prop.~\ref{ramified}, we
may assume that $e_{\{x\}} = 1$ and $y = u\pi^n$ for some $u \in
\mfo$ and $n$ not divisible by $p$, and $\depth \{ x, y \} =
\depth \{ x \}$.  On the one hand, $D$ contains $F(\psi)$ on which
the valuation ramifies, so $e_{D/F} = p$ as claimed.  On the other
hand, $D$ contains $F(\chi)$ whose residue algebra is a proper
extension of $\Fb$ that is inseparable (if $\depth \{ x, y \} <
m$) or separable (if $\depth \{ x, y \} = m$) by Example
\ref{deg1.even}; this proves the claim.  We are left with the case
where the depth is divisible by $p$ and $e_{\{x, y \}} = 1$.

If the depth of $\{ x, y \}$ is zero, then by the Kato isomorphisms
$\{ \bx, \by \}$
is not zero in $k_2(\Fb)$, hence $\bx, \by$ are $p$-free over
$\Fb^p$.  Following
\ref{pf.dict.1}, $\bchi$ and $\bpsi$ generate purely inseparable
extensions of $\Fb$ in $\Db$.  The value of $\chi \psi - \psi \chi
= (\zeta - 1) \psi \chi$ is $\val(\zeta - 1) = m/p > 0$, so
$\bchi$ and $\bpsi$ commute in $\Db$.  We deduce that $\Db$ is
$\Fb(\sqrt[p]{\bx}, \sqrt[p]{\by})$, $e_{D/F} = 1$, and $D$ is
wild.

If $0 < \depth \{ x, y \} < m$, then we can choose $x = 1 + a\pi^d$
where $d = \depth \{ x, y \}$ so that the initial form of $\{ x, y
\}$ in $\Omega^1 / Z^1$ is $\ba \dlog{\by}$.  As the depth is $d$,
$\ba \dlog{\by}$ is not in $Z^1$, i.e., $\mathrm{d}{\ba} \wedge
\dlog{\by}$ is not 0.  It follows that $\dim_{\Fb}
\Fb(\sqrt[p]{\ba}, \sqrt[p]{\by}) = p^2$ by \ref{Kdiff}. We claim
that this field is $\Db$.  By dimension count, it suffices to note
that for $\alpha := \pi^{-d/p}(\chi - 1)$, we have
\[
\balpha^p = \ba, \quad \bpsi^p = \by, \eand \balpha \bpsi = \bpsi \balpha.
\]
The first equation is as in Example \ref{deg1.even} and the second
is obvious.  The third follows because
\[
\alpha \psi - \psi \alpha = \pi^{-d/p}(\chi \psi - \psi \chi) =
\pi^{-d/p} (\zeta - 1) \psi \chi.
\]
But $\zeta - 1$ has value $m / p > d / p$, hence this commutator
has positive value and $\balpha, \bpsi$ commute.  This shows that
$e_{D/F} = 1$ and $D$ is wild.

Finally suppose that the depth of $\{ x, y \}$ is $m$.  Then we
write $x = 1 + b \eta^p$ and $\beta = \eta^{-1} (\chi - 1)$ as in
Example \ref{deg1.even}; again $\bbeta^p - \bbeta - \bb= 0$.
Further,
\[
\beta \psi - \psi \beta = \eta^{-1} (\chi\psi - \psi \chi) = \psi \chi,
\]
so
\[
\bbeta \bpsi = \bpsi (\bbeta + 1).
\]
The elements $\bbeta$ and $\bpsi$ generate a division algebra of
dimension $p^2$ over its center $\Fb$ \cite[p.~36]{GilleSz}, so it
must be $\Db$.  Therefore, $e_{D/F} = 1$ and $D$ is
tame.$\hfill\qed$

%\subsection{} \label{pf.dict.3} \emph{Proof of Proposition \ref{dict}
%for $q = 3$ and $p = 2$}: The last case of Prop.~\ref{dict} is
%where $p = 2$ and $\gamma$ is a nonzero symbol in $k_3(F)$, so
%$\gamma$ corresponds to an octonion division algebra $D$.

%Roughly speaking, we repeat the calculations in the case $q = 2$
%from \ref{pf.dict.2}.  For example, if the depth is odd (so
%$e_\gamma = 2$), we find that $e_{D/F} = p$ and $D$ has residue
%algebra a purely inseparable field extension of $\Fb$.  (To see
%that the residue algebra is commutative, one imitates the
%computation in the $q = 2$, depth zero, $e = 1$ case of
%\ref{pf.dict.2}.)

%In case the depth is maximal (i.e., $m$), the residue algebra is a
%composition algebra.  For example, suppose that we are in case
%(ii) of the Gathering Lemma, so that the symbol is of the form
%$\{ u, b \pi, a \}$ for some $u \in U_m$, $b \in \Ox$, and $a \in
%\Ox$ whose residue is not a $p$-th power. Recall that here $(\zeta
%- 1)^2 = 4$ has value $m$, so $u = 1 + 4x$ for some $x \in \Ox$.
%The Kato isomorphism maps
%\[
%\{ u, b\pi, a\} \mapsto \left( \bx \dlog{\bb} \wedge \dlog{\ba}, \bx
%\dlog{\ba}\right) \in H^1(\Fb, \nu(2)) \oplus H^1(\Fb, \nu(1)).
%\]
%As the depth is $m$, this image is not zero, so $\bx \dlog{\ba}$
%is not zero.  We deduce that the residue algebra $\Db$ is a
%quaternion division algebra as in the case $q = 2$, hence also
%that $e_{D/F} = 2$.$\hfill\qed$

%%%%%%%%%%%%%%%%%%%%%%%%%%%%%%%%%%%%%%%%%%%%%%%%%%%%%%%%%%%%%%%
\section{Proofs of Theorem \ref{texp.depth} and Corollary \ref{texp.cor}} \label{pf.texp.depth}

The following proofs amount to translating results of
\S\S\ref{s.VADA}, \ref{s.VADAEN} into $K$-theory.

\begin{proof}[Proof of Theorem \ref{texp.depth}] It suffices to
establish (i) and (iii) since by Prop.~\ref{p.PROTE}(a) (iii) implies (ii).

\emph{\underline{Case $q=1$}}: Suppose first that
$q = 1$, i.e., $\gamma = \{ \mu \}$ for some $\mu \in \Fx$ and $Q$
is the 1-Pfister $\pform{\mu}$.  Put $P := \qform{1}$; it is wild
because $\chr \Fb = 2$ and $\texp P = \val(2)$ by Example
\ref{e.BAFI}.  We may assume that $\mu$ has value 0 or 1.  If $\mu$
has value 1, then $\depth \gamma = 0$, $e_\gamma = 2$, and the
theorem holds by Prop.~\ref{p.CDPR}.\ref{VADAPR}.  If $\mu$ has
value 0 and $\bar{\mu}$ is a nonsquare in $\Fb$, then $\depth \gamma
= 0$, $e_\gamma = 1$, and the theorem holds by
Prop.~\ref{p.CDUNEV}.\ref{TREXPREV}.

Otherwise $\overline{\mu}$ is a square in $\Fb$, so multiplying
$\mu$ by a square in $F$ we may assume that $\overline{\mu} = 1$. If
$\depth \gamma < 2 \val(2)$, then \eqref{texp.e} and
\eqref{texp.depth.eqn} hold by Prop.~\ref{p.CDUNOD}.\ref{TREXPROD}
or Prop.~\ref{p.CDUNEV}.\ref{TREXPREV}.

Finally suppose that $\depth \gamma = 2\val(2)$.  The symbol
$\gamma$ corresponds to the quadratic extension $F(\sqrt{\mu})$ and
the residue algebra of this extension was computed in Example
\ref{deg1.even}; this verifies \eqref{texp.e} and that $Q$ is tame
(hence \eqref{texp.wild}).  Then $\texp Q = 0$ by
Prop.~\ref{p.PROTE}(d), proving \eqref{texp.depth.eqn}.

\smallskip
\emph{\underline{Case $q \ge 2$}}: We argue by induction on $q$
and decompose $\gamma$ as in Prop.~\ref{unramsubsym} (with $p = 2$).
Writing $P$ (resp.~$Q$) for the Pfister quadratic form corresponding
to $\beta$ (resp.~$\gamma$), we conclude $Q \cong \dla a \dra
\otimes P$ and put $d := \depth\gamma$. Then $e_{P/F} = 1$ and
$\texp(P) = \lambda(2) - (\depth\beta)/2$ by the induction
hypothesis. If alternative (ii) of Prop.~\ref{unramsubsym} holds,
then (\ref{p.CDPR}.\ref{VADAPR}) shows $e_\gamma = 2 = e_{Q/F}$ and
$\texp(Q) = \texp(P) = \lambda(2) - \floor{\frac{\depth\gamma}{2}}$
since $\depth\gamma = \depth\beta$ is even. Now suppose alternative
(i) of Prop.~\ref{unramsubsym} holds. Since $\beta$ has depth zero
and ramification index $1$, we can write $\beta = \{a_1,\dots,
a_{q-1}\}$, $a_i \in \mfo^\times$, $1 \leq i < q$. If $d =
2\lambda(2)$, then $e_\gamma = 1$ and $\dla a \dra$ is a \emph{tame}
$1$-Pfister quadratic subspace of $Q \cong \dla
a,a_1,\dots,a_{q-1}\dra$ with $e_{\dla a \dra/F}= 1$. Applying
Prop.~\ref{p.LATAU} $q - 1$ times yields $e_{Q/F} = 1 = e_\gamma$
and that $Q$ is tame as well, forcing $\texp(Q) = 0 = \lambda(2) -
\floor{d/2}$. We are left with the case $0 \leq d <2\lambda(2)$.
Assertions (i) and (iii) follow by combining the effect of the Kato
isomorphism at depth $d$ on $\gamma$ with \ref{Kdiff} and
(\ref{p.CDUNOD}.\ref{TREXPROD}) (for $d$ odd) or with
(\ref{p.CDUNEV}.\ref{TREXPREV}) and Cor.~\ref{c.CHARNEXOD}(a) (for
$d$ even).
\end{proof}

\begin{proof}[Proof of Corollary \ref{texp.cor}]
We now prove Corollary \ref{texp.cor}.  Claim (i) amounts to
reformulating Proposition \ref{ramified}.  For (ii), one combines
Theorem \ref{texp.depth}\eqref{texp.depth.eqn} relating $\depth
\gamma$ with $\texp C$ with the relation between $\texp C$ and
$\wid(C)$ demonstrated in the proof of Theorem \ref{t.HEWI}.
\end{proof}

%%%%%%%%%%%%%%%%%%%%%%%%%%%%%%%%%%%%%%%%%%%%%%%%%%%%%%%%%%%%%%%

\noindent\small{\textbf{Acknowledgements.} The first author was
partially supported by National Science Foundation grant
no.~DMS-0653502.  Both authors thank Eric Brussel, Detlev Hoffmann, Daniel Krashen, Ottmar
Loos, Pat Morandi, David Saltman, and Jean-Pierre Tignol for useful
comments and suggestions. Special thanks are due to Richard Weiss,
whose stimulating questions (related to his recent work \cite{We09}
on affine buildings) inspired us to look for successively improved
versions of the Local Norm Theorems~\ref{t.HENOTH} and \ref{LNTK}.}

\providecommand{\bysame}{\leavevmode\hbox to3em{\hrulefill}\thinspace}
\providecommand{\MR}{\relax\ifhmode\unskip\space\fi MR }
\providecommand{\MRhref}[2]{%
  \href{http://www.ams.org/mathscinet-getitem?mr=#1}{#2}
}
\providecommand{\href}[2]{#2}

\printindex
\end{document}